\documentclass[11pt,a4paper]{book}
\usepackage[utf8]{inputenc}
\usepackage{amsmath}
\usepackage{amsfonts}
\usepackage{amssymb}
\usepackage{amsthm}
\usepackage{mathtools}
\usepackage{mathrsfs}
\usepackage{cite}
\usepackage[titletoc]{appendix}
\usepackage{hyperref}
\usepackage{relsize}
\usepackage{stmaryrd}
\usepackage[final]{pdfpages}

\swapnumbers

\newtheorem{thm}{Theorem}
\newtheorem*{thm*}{Theorem}
\newtheorem*{theorem*}{Théorème}

\theoremstyle{definition}
\newtheorem{defn}[thm]{Definition}

\theoremstyle{remark}

\newcommand{\R}{\mathbb{R}}
\newcommand{\Z}{\mathbb{Z}}
\newcommand{\h}{\mathcal{H}}
\newcommand{\T}{\mathbf{T}}
\newcommand{\Y}{\mathbf{Y}}
\newcommand{\W}{\mathbf{W}}
\newcommand{\sk}{\mathbf{sk}}
\newcommand{\D}{\Delta}

\renewcommand{\H}{\mathcal{H}}

\begin{document}







\begin{figure}[h]
\vspace*{-4,4cm}
\hspace*{-3,39cm}
\includegraphics[scale=0.3595]{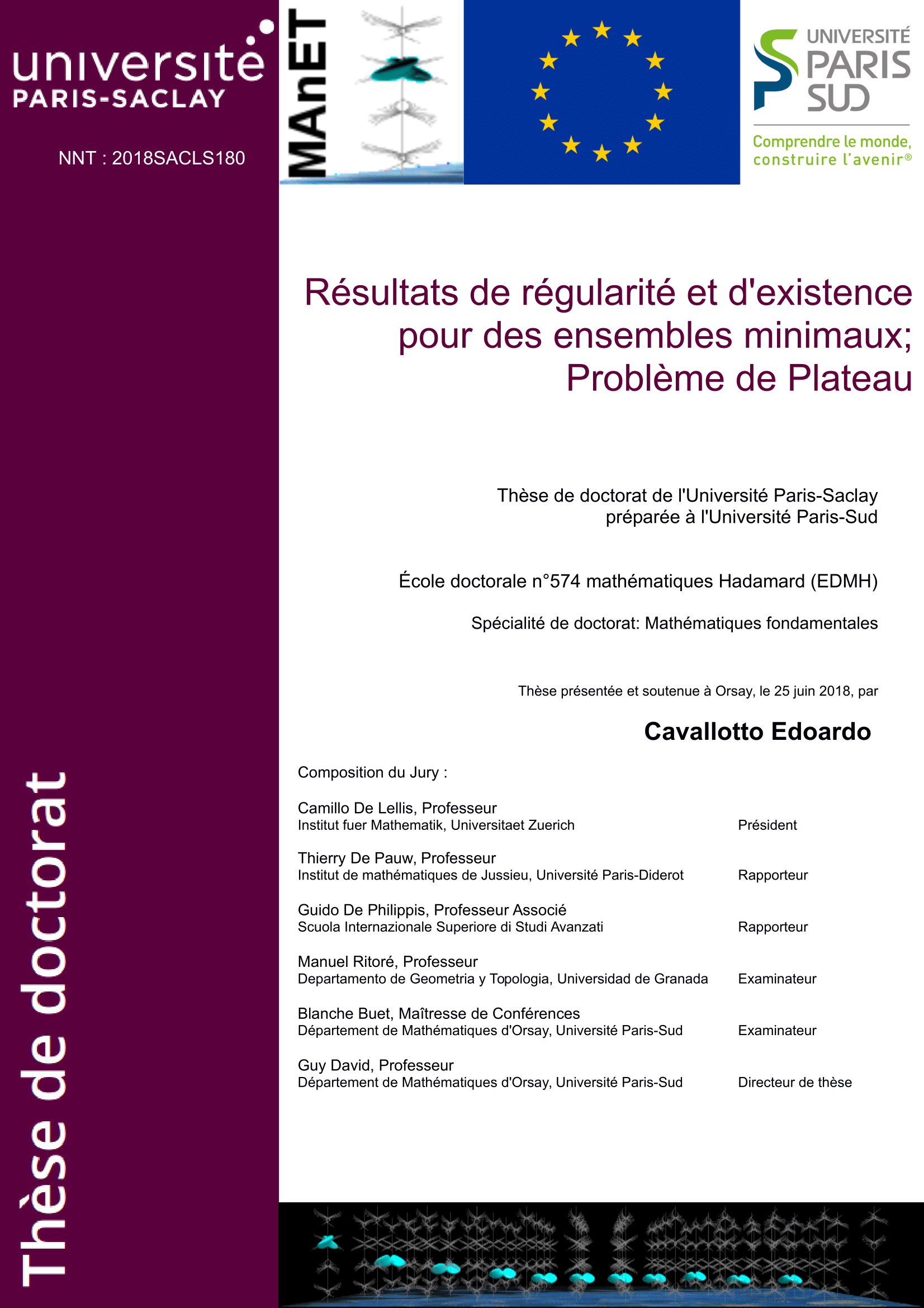}
\end{figure}

\clearpage

\phantom{pain}

\newpage

\begin{figure}[h]
\vspace*{-4cm}
\hspace*{-4cm}
\includegraphics[scale=0.35]{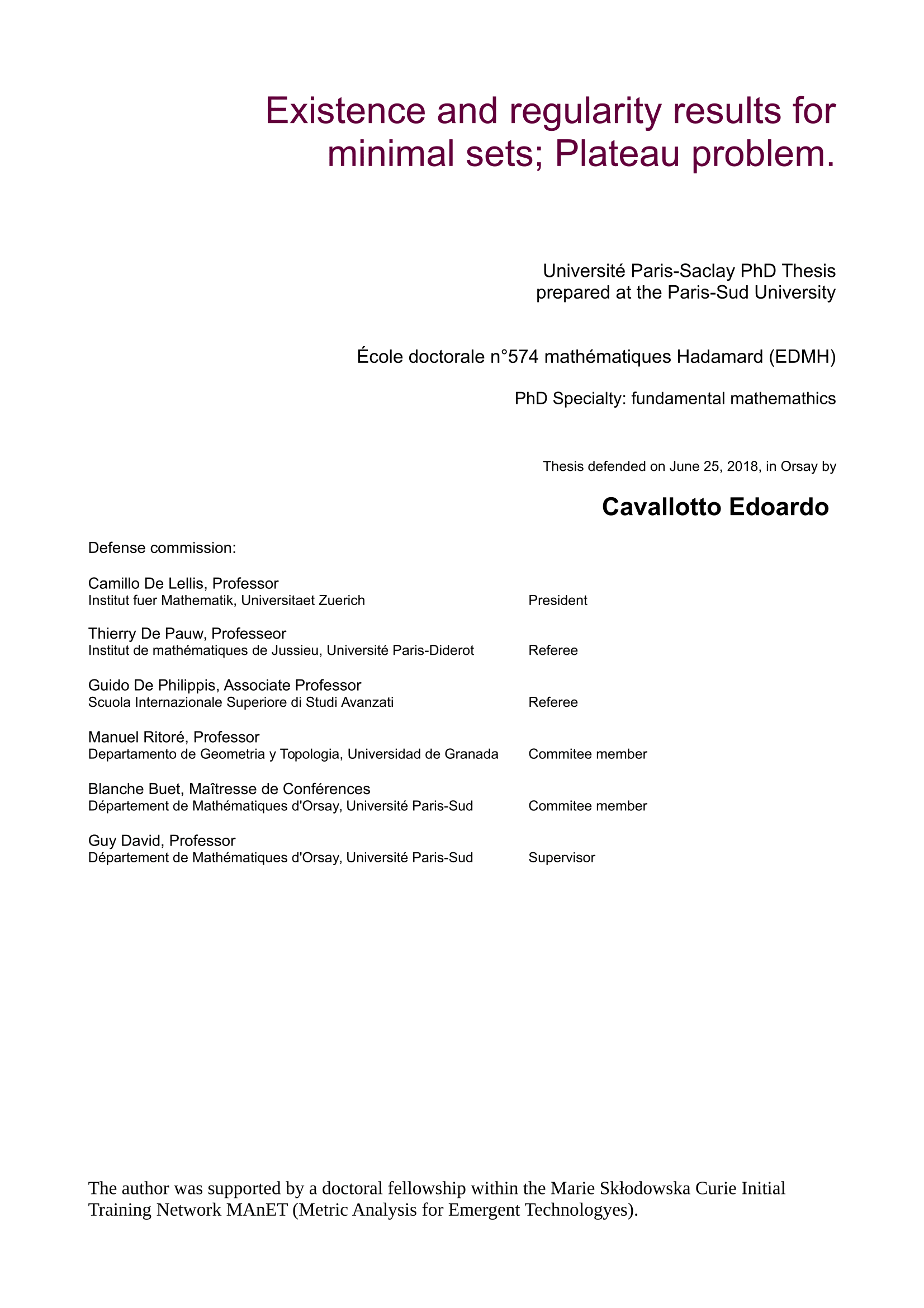}
\end{figure}

\newpage

\tableofcontents

\chapter*{Introduction}

\addcontentsline{toc}{chapter}{Introduction (Français)}

\paragraph{Contexte historique}
Le Problème de Plateau tire son origine de la physi-que et il porte le nom du physicien belge Joseph Plateau (1801-1883), qui a étudié le comportement des bulles de savon et a décrit les singularités typiques qu'elles produisent (voire \cite{plateau1873statique}). Résoudre le Problème de Plateau signifie trouver la surface ayant l'aire minimale parmi toutes les surfaces avec un bord donné. Une partie du problème réside dans le fait de donner des définitions appropriées aux concepts de ``surface'', ``aire'' et ``bord''. Au fil du temps plusieurs formulations différentes ont été données au Problème de Plateau (pour une vue d'ensemble voir \cite{david2013regularity}) et, en fait, ce problème a fourni la motivation principale pour le développement de la théorie géométrique de la mesure.

Un premier essai de résolution du Problème de Plateau s'est fait en termes du problème de Dirichlet, dans lequel on cherche à minimiser l'énergie de Dirichlet parmi les fonctions avec un valeur donnée au bord. C'est à dire
\begin{displaymath}
min\left\{\int_\Omega|\nabla u|^2\,:\, u_{|\partial\Omega}=\gamma\right\}
\end{displaymath}
où $\Omega\subset\mathbb{R}^3$ et $\gamma:\partial\Omega\to\mathbb{R}^3$ est fixée.

Une autre formulation possible consiste à injecter un disque de dimension 2 dans $\mathbb{R}^3$ de sorte que son bord soit envoyé sur un chemin fermé et fixé en minimisant la fonctionnelle d'aire associée à l'image du disque (voire \cite{douglas1931solution}). Autrement dit: soit $D$ un disque de dimension 2, et $\Gamma$ une courbe simple et fermée dans $\R^3$ paramétrée par $g:\partial D\to\mathbb{R}^n$; alors une solution pour le Problème de Plateau est une fonction $f:D\to\mathbb{R}^3$ telle que $f_{|\partial D}=g$ et $f$ minimise 
\begin{displaymath}
A(f):=\int_D J_f(x)dx,
\end{displaymath}
où $J_f$ est le Jacobien de $f$.

Dans les deux cadres précédents le Problème de Plateau est bien défini et il existe bien des solution, par contre on ne veut pas forcement supposer que notre bulle de savon soit le graphe d'une fonction lisse ou bien l'image d'un disque. De plus on voudrait pouvoir définir notre problème dans un espace euclidien des dimension arbitraire $n>0$, et pour de surfaces de dimension arbitraire $d>0$ telle que $0<d<n$.

Dans \cite{reifenberg1960solution} la condition de bord a été définie par Reifenberg en utilisant l'homologie de Čech et il a montré ensuite qu'il est possible minimiser la mesure de Hausdorff de dimension $d$ par rapport à cette condition. Les solutions de Reifenberg nous donnent un bonne description de plusieurs films de savon différentes, mais il y a encore des exemples physiques de films de savon qui ne peuvent pas être obtenus dans ce contexte.

Probablement la formulation du Problème de Plateau la plus célèbre et qui a eu le plus de succès est celle donnée par Federer et Fleming en termes des courants (voire \cite{federer2014geometric} et \cite{federer1960normal}). Les courants sont définis comme des distributions sur l'espace des formes différentielles, et pour les courants on peut définir de façon assez naturelle des notions de bord et d'aire. En particulier un courant intégral de dimension $d$ peut être vu, de façon très libre, comme une surface $(\H^d,d)$ rectifiable équipée d'une orientation et d'une multiplicité, et telle que son bord est $(\H^{d-1},d-1)$ rectifiable. Malheureusement l'orientation d'un courant dépend du groupe des coefficients choisi, et en autre la fonctionnelle d'aire doit tenir compte de la multiplicité de la surface. A cause de ce deux raisons on est pas entièrement satisfait de ce modèle de film de savon.

\begin{figure}
\centering
\includegraphics[scale=0.3]{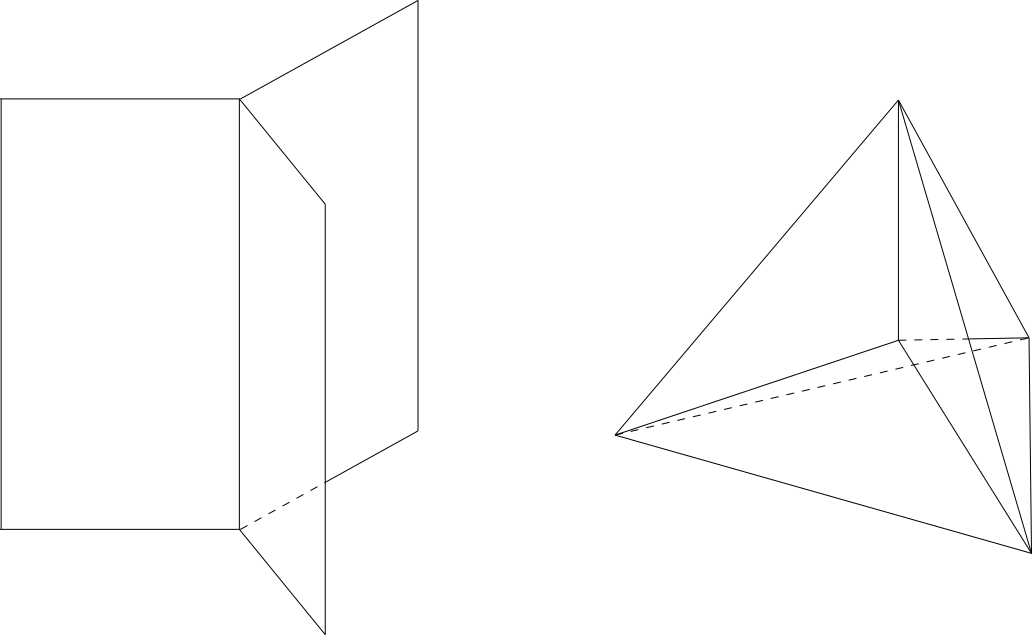}
\caption{Les cônes de type $\mathbf{Y}$ et $\mathbf{T}$.}
\label{YeTintro.pngfr}
\end{figure}
\paragraph{Almgren}
Dans \cite{almgren1968existence} Almgren a introduit un description plus ``naturelle'' des bulles de savon qui ensuite a été étudiée par David et Semmes dans \cite{david2000uniform}. Les objets considérés dans ce contexte sont des ensembles dont la mesure de Hausdorff $\H^d$ est localement finie, la fonctionnelle à minimiser appartient à la classe des intégrants elliptiques (laquelle contient la mesure de Hausdorff $\H^d$), et la condition de bord est donnée par rapport à une famille à un paramètre de déformations compactes; finalement on dit qu'un ensemble est \emph{minimal} si son énergie ne peut pas être réduite par aucune déformation continue qui a lieu dans une bulle qui n'intersecte pas le bord. Almgren a montré que, en dehors d'un ensemble négligeable, les ensembles minimaux sont des sous-variétés différentielles d'ordre $C^{1,\alpha}$ plongées dans $\R^n$ et que le cône tangent à un ensemble minimal en n'importe quel point est un cône minimal (c'est à dire un ensemble invariant par dilatations qui est aussi un ensemble minimal). Donc pour comprendre complètement le comportement local d'un ensemble minimal il est vital d'étudier les cônes minimaux.  

Dans $\mathbb{R}^2$ il n'y a que deux types de cônes minimaux: les droites, et les cônes de type $Y$, formés par trois demi-droites qui se rencontrent avec des angles de 120 degrés. La liste complète des cônes minimaux de dimension deux dans $\mathbb{R}^3$ a été proposé à l'origine par Plateau, et ensuit a été démontrée par Taylor en \cite{taylor1976structure}; il y en a de 3 types: les plans,  $\mathbf{Y}:=Y\times\mathbb{R}$, et $\mathbf{T}$, le cône sur les arêtes d'un tétraèdre régulier (voire la Figure \ref{YeTintro.pngfr}). Pour $n\ge4$ on ne connaît que des exemples isolés de cônes minimaux, et la liste semble loin d'être complète (voire \cite{brakke1991minimal} \cite{lawlor1994paired} \cite{liang2012almgren} \cite{liang2014almgren}).

\paragraph{Bord glissant}
Une nouvelle notion de bord, le \emph{bord glissant}, a été introduite par David dans \cite{david2014local} afin d'investiguer la régularité prés du bord des ensembles minimaux. La partie d'une surface soumise à cette condition qui touche le bord n'est plus fixée comme dans la définition d'Almgern mais libre de se délacer le le long d'un domaine fermé ayant la fonction d'une coulisse. Un exemple physique ou le bord glissant s'applique à une surface est celui d'une pellicule de savon contenue dans un tube; dans ce cas la pellicule peut se déplacer à l'intérieur du tube de façon à ce que son bord reste toujours en contact avec la surface intérieure du tube sans jamais la laisser. David a montré que les ensemble minimaux glissants (c'est à dire les ensemble minimaux dans ce cadre) sont uniformément rectifiables; il a entre autre montré que, sous des hypothèses raisonnables de régularité pour le bord, le cône tangent à un ensemble minimal glissant par rapport à un point sur son bord est un cône minimal glissant par rapport à un bord conique. Comme en supposant plus de régularité pour le bord (par exemple $C^1$ ou rectifiable) on obtient que son tangent est plat (respectivement partout ou presque partout), il suit que l'investigation des cônes minimaux par rapport à un bord plat est une étape important vers une compréhension complète du comportement des ensembles minimaux près du bord.

Dans cette thèse nous nous sommes intéressés au cônes minimaux glissants de dimension $d$ contenus dans le demi-espace $\mathbb{R}^n_+:=\{x_n\ge0\}$, avec $0<d<n$, où le domaine du bord glissant est l'hyperplan ``horizontal'' $\Gamma:=\{x_n=0\}$ qui borne le demi-espace. La fonctionnelle que on a considéré n'est pas la mesure de Hausdorff tout court, mais, en vue de fonctionnelles plus générales (les intégrands elliptiques) nous nous sommes concentrés sur une petite modification de celle-ci définie comme il suit
\begin{displaymath}
J_\alpha(E):=\mathcal{H}^d(E\setminus\Gamma)+\alpha\mathcal{H}^d(E\cap\Gamma)
\end{displaymath}
pour $E\subset\mathbb{R}^n_+$ (où $\alpha\in[0,1]$ pour des raisons de semi-continuité).

\begin{figure}
\centering
\includegraphics[scale=0.3]{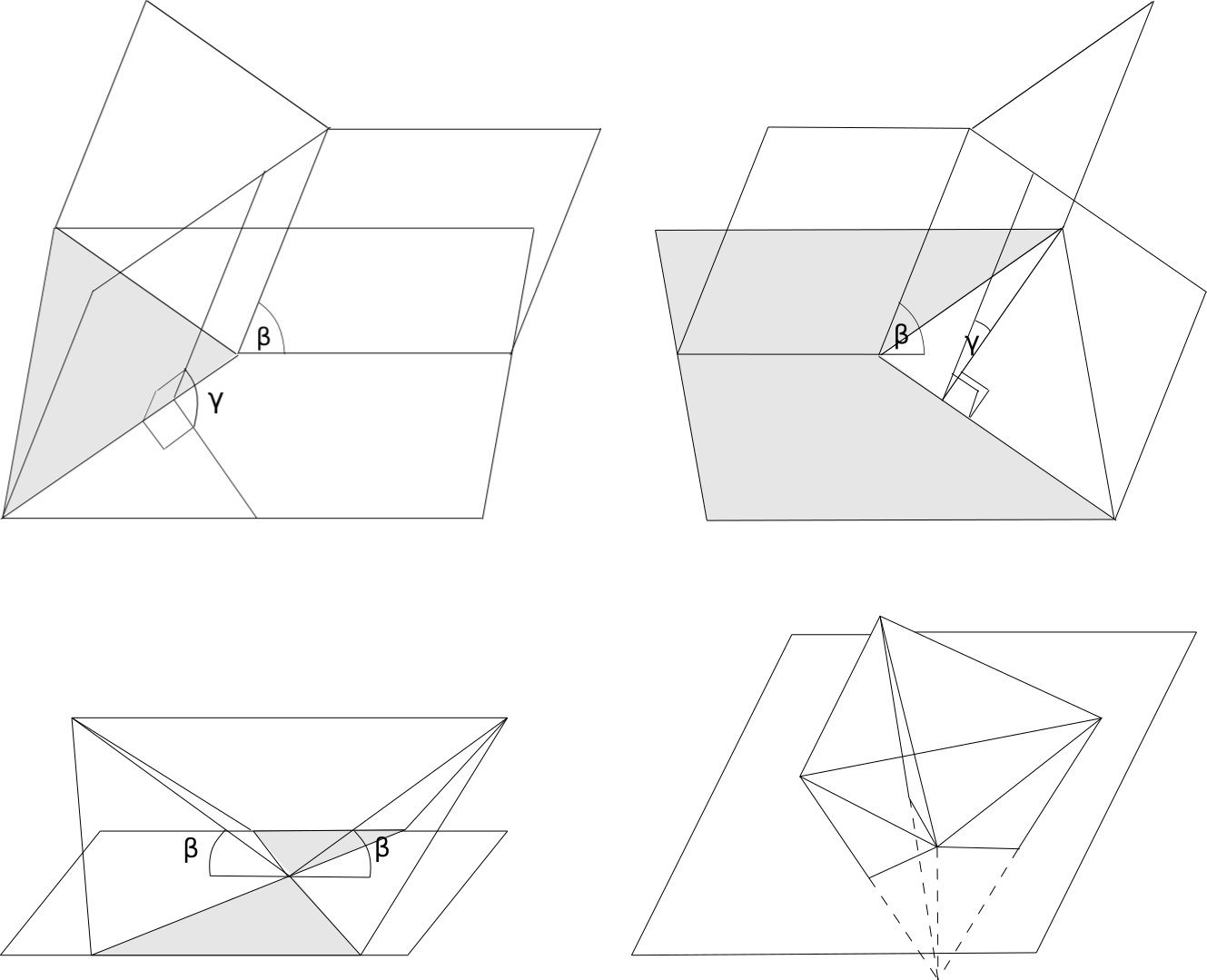}
\caption{Les cônes minimaux de dimension 2 dans le demi-espace de dimension 3, dans l'ordre: $\Y_\beta$, $\overline{\Y}_\beta$, $\W_\beta$, $\T_+$. Pour chacun la région grise représente l'intersection entre le cône et le plan horizontal.}
\label{coniscorrevolifr}
\end{figure}

On a montré la minimalité de 4 nouveau types des cônes de dimension 2 dans le demi-espace de dimension 3 (autres que ceux qui peuvent être obtenus comme produit cartésien avec $\R$ d'un cône minimal de dimension 1 contenu dans le demi-plan) ou pour chaque type de cône on a en fait trouvé une famille à un paramètre de cônes minimaux glissants dépendant du paramètre $\alpha$ (voire Figure \ref{coniscorrevolifr}).

L'argument principal utilisé pour montrer la minimalité de ces cônes est l'utilisation des \emph{calibrations couplées}, un outil employé par Lawlor et Morgan dans \cite{lawlor1994paired} pour montrer la minimalité du cône sur les arêtes de dimension $n-2$ d'un simplexe régulier de dimension $n$, qui consiste à appliquer le théorème de la divergence à chacune des composants connexes du complémentaire du cône par rapport à une famille de champs vecteurs à divergence nulle proprement choisie.

Les deux premiers exemples de cônes minimaux glissants peuvent se construire avec la procédure suivante. On se donne un cône de type $\mathbf{Y}$ contenu dans $\R^3$ et penché de sorte que l'intersection entre les trois demi-plans rencontre le plan horizontale $\Gamma$ en formant un angle $\beta\in[0,\pi/2]$ et qu'une de trois nappes rencontre $\Gamma$ orthogonalement. Le cône obtenu en intersectant le précédent avec le demi-espace $\R^3$ par construction est formé de deux nappes inclinées, qui forment un même angle $\gamma$ (qui dépend de $\beta$) avec $\Gamma$, et d'une nappe verticale ayant la forme d'une secteur de plan décrit par l'angle $\beta$ ou bien $\pi-\beta$. Dans le premier cas on peut définir le cône $\Y_\beta$ comme la réunion de l'ensemble obtenu avec la construction précédente avec la région de $\Gamma$ contenue entre le deux nappes inclinées, et dans le deuxième cas on peut définir le cône $\overline{\Y}_\beta$ comme la réunion de l'ensemble obtenu avec la construction précédente avec la région de $\Gamma$ qui n'est pas contenue entre le deux nappes inclinées. Les deux cônes définis de cette manière sont minimaux si et seulement si $\cos\gamma=\alpha$.

On peut construire le troisième cône, appelé $\mathbf{W}_\beta$, comme la réunion entre l'intersection de $\overline{\mathbf{Y}}_\beta$ avec un demi-espace borné par un plan vertical $P$, orthogonal à la nappe verticale, et sa réflexion par rapport à $P$. Le cône $\mathbf{W}_\beta$ est minimal si et seulement si $\beta\le30^\circ$ et $\cos\gamma=\alpha$.

Le quatrième cône minimal s'appelle $\mathbf{T}_+$ et il est tout simplement l'intersection d'un cône de type $\T$ proprement placé avec le demi-espace $\R^3_+$. Par un argument de calibration on montrera que $\T_+$ est un cône minimal glissant pour tout $\alpha\ge\sqrt{\frac{2}{3}}$, en outre, pour tout $\alpha<\sqrt{\frac{2}{3}}$ on peut obtenir un compétiteur dont l'énergie est strictement plus petite par rapport à l'énergie du cône en poussant une partie de cône sur le plan horizontal de sorte qu'un petit triangle horizontal est produit et en rangeant proprement les nappes inclinées. L'argument de calibration précédent peut s'étendre en grand dimension comme suit
\begin{theorem*}
Soit $n\ge3$ et $\alpha_n=\sqrt{\frac{n+1}{2n}}$. Si $\alpha\ge\alpha_n$ alors le cône de dimension $(n-1)$ sur les arêtes de dimension $(n-2)$ d'un simplex régulier de dimension $n$ proprement placé et intersecté avec le demi-espace $\mathbb{R}^n_+$ est un cône minimal glissant.
\end{theorem*}

Dans le Chapitre \ref{minimal sets} nous introduirons les ensembles minimaux d'Almgren, en décrivant leurs propriétés et les cône minimaux produits. Ensuite nous introduirons la notion de bord glissant et finalement nous décrirons le cadre auquel nous sommes intéressés.

Le Chapitre \ref{One-dimensional} sera dédié à l'étude des cônes minimaux de dimension 1 dans le demi-plan.

Dans le Chapitre \ref{two-dimensional cones} nous nous concentrerons sur les cônes minimaux de dimension 2 dans le demi-espace, en essayant de les classifier par rapport au graphe engendré par leur intersection avec la sphere unitaire. Des première partie du Chapitre nous discuterons les cônes minimaux glissants, à partir de ceux qui peuvent s'obtenir comme un produit cartésien et des conditions nécessaires qu'un cône doit satisfaire pour pouvoir être minimal. Dans la deuxième partie du Chapitre nous présenterons une vue d'ensemble, malheureusement pas totalement exhaustive à ce jour, sur les cônes pas minimaux qui satisfont aux conditions nécessaires pour la minimalité, tout en essayent de les classifier.

Dans le Chapitre \ref{higher dimension} nous démontrerons le Théorème énoncé ci-dessus, qui généralise en haute dimension l'argument par calibration employé pour montrer la minimalité du cône $\T_+$ dans le demi-espace de dimension 3.

\chapter*{Introduction}

\addcontentsline{toc}{chapter}{Introduction (English)}

\paragraph{Historical background}
The Plateau problem arises from physics, and it is named after the Belgian physicist Joseph Plateau (1801-1883), who studied the behaviour of soap films and described the typical singularities that they produce (see \cite{plateau1873statique}). Solving the Plateau problem means finding the surface with minimal area among all surfaces with a given boundary. Part of the problem actually consists in giving a suitable definition to the notions of ``surface'', ``area'' and ``boundary''. Over time a number of different settings for the Plateau problem have been studied (for an overview on the topic see \cite{david2013regularity}), and solving this problem actually provided the main motivation for the development of geometric measure theory.

A first possible setting for the Plateau problem is the Dirichlet problem, in which one wants to minimise the Dirichlet energy among all the functions with a given boundary data, i.e.
\begin{displaymath}
min\left\{\int_\Omega|\nabla u|^2\,:\, u_{|\partial\Omega}=\gamma\right\}
\end{displaymath}
where $\Omega\subset\mathbb{R}^3$ and $\gamma:\partial\Omega\to\mathbb{R}^3$ is fixed.

Another possible formulation consists of mapping a 2-dimensional disc into $\mathbb{R}^3$ in such a way that the boundary is mapped onto a given closed path, and then minimising the area functional (see \cite{douglas1931solution}). That is to say: let $D$ be the two-dimensional unit disc, and $\Gamma$ a simple closed curve in $\mathbb{R}^3$ parameterised by $g:\partial D\to\mathbb{R}^n$; a solution to the Plateau problem is a function $f:D\to\mathbb{R}^3$ such that $f_{|\partial D}=g$ and $f$ minimises the following quantity
\begin{displaymath}
A(f):=\int_D J_f(x)dx,
\end{displaymath}
where $J_f$ is the Jacobian of $f$.

The two previous settings provide positive solutions to the Plateau problem. However we do not necessarily want to assume our surfaces to be graphs of smooth functions or images of discs. Moreover we would like to state our problem in a space with arbitrary dimension $n>0$ and for surfaces of any dimension $d$, with positive integers such that $d<n$.

In \cite{reifenberg1960solution} Reifenberg stated the boundary condition in terms of Čech homology, and then proved that the $d$-dimensional Hausdorff measure can be minimised under this constraint. Reifenberg’s solutions are nice and seem to give a good description of many soap films, but still there are some physical examples of soap films spanned by a curve that cannot be obtained in this framework.

Perhaps the most celebrated and successful model is the description of films in terms of currents, given by Federer and Fleming (see e.g. \cite{federer2014geometric} and \cite{federer1960normal}). Currents are defined as distributions on differential forms and they are naturally endowed with a notion of boundary and area. In particular a $d$-dimensional integral current is, loosely speaking, a $(\H^d,d)$ rectifiable surface endowed with an orientation and a multiplicity, and such that its boundary is $(\H^{d-1},d-1)$ rectifiable. However the orientation of a current depends on the group of coefficient chosen, and the area functional takes into account the multiplicity of a surface. For these reasons currents are not entirely satisfying as a model for soap films.

\paragraph{Almgren setting}
A more ``natural'' setting was introduced by Almgren in \cite{almgren1968existence} and later studied by David and Semmes in \cite{david2000uniform}. In this framework the considered objects are sets with locally finite $d$-dimensional Hausdorff measure, the functional to be minimised belongs to the class of elliptic integrands (which contains the Hausdorff measure itself), and the boundary condition is given in terms of a one parameter family of compact deformations. In this setting a set is said to be a \emph{minimiser} if its energy cannot be decreased by any continuous deformation acting in a ball that does not intersect the boundary. Almgren also proved that minimisers are $C^{1,\alpha}$ embedded submanifolds of $\mathbb{R}^n$ up to a negligible set, and that the tangent cone to any point of such a minimiser is a minimal cone. Therefore in order to completely understand the local behaviour of minimal surfaces one has to know what minimal cones look like. 
\begin{figure}
\centering
\includegraphics[scale=0.3]{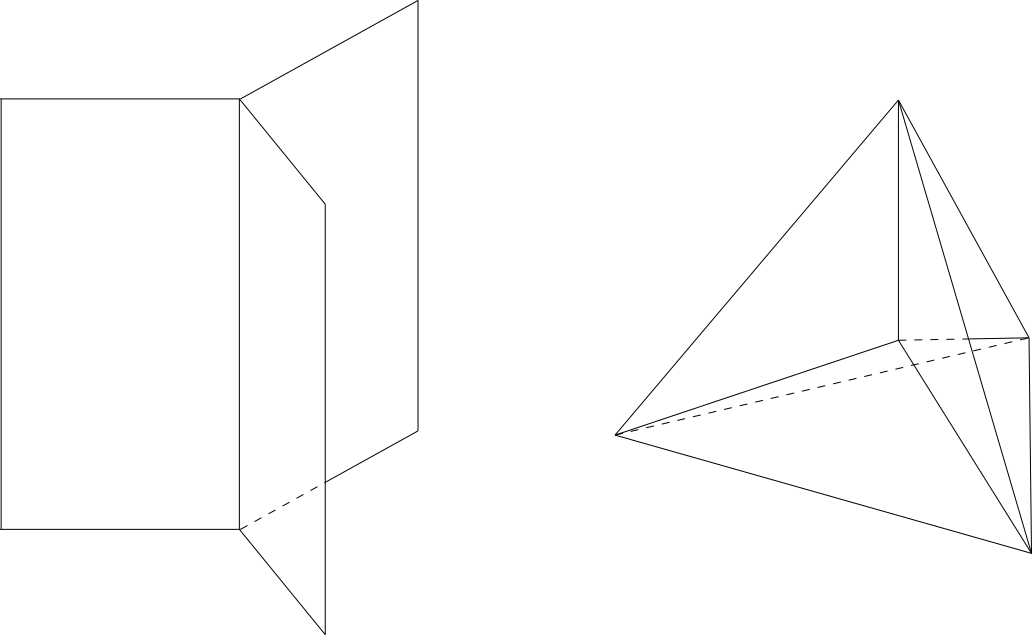}
\caption{Cones of type $\mathbf{Y}$ and $\mathbf{T}$.}
\label{YeTintro.png}
\end{figure}

In $\mathbb{R}^2$ there are only two types of minimal cones: straight lines, and the cone called $Y$, formed by three half-lines meeting with equal angle of 120 degrees. The complete list of minimal cones of dimension 2 in $\mathbb{R}^3$ was first conjectured by Plateau and then proved by Taylor in \cite{taylor1976structure}; it has three types of cones: planes,  $\mathbf{Y}:=Y\times\mathbb{R}$, and $\mathbf{T}$, the cone over the edges of a regular tetrahedron (see Figure \ref{YeTintro.png}). When $n\ge4$ the list of minimal cones looks far from complete and we only know few examples (see e.g. \cite{brakke1991minimal} \cite{lawlor1994paired} \cite{liang2012almgren} \cite{liang2014almgren}).


\paragraph{Sliding boundary}
In order to study boundary regularity, in \cite{david2014local} David introduced a new notion of boundary, called \emph{sliding boundary}. Loosely speaking the boundary of a surface subject to this condition is not fixed but is allowed to move in a closed set. A physical example where this condition applies is a soap film contained in a tube: the boundary of the film can move along the inner surface of the tube without leaving it. David proved that sliding minimisers (i.e. minimal surfaces in this new setting) are uniformly rectifiable; moreover he proved that under some mild regularity condition of the boundary, the blow-up limit of a sliding minimal set at a boundary point is a sliding minimal cone with respect to a conical boundary. The assumption of more regularity of the boundary (like $C^1$ or rectifiable) provides flatness (everywhere or almost-everywhere) of its blow-up. Therefore an important step toward understanding the behaviour of sliding minimisers close to the boundary is to list the sliding minimal cones with respect to flat boundaries.

The subject of this thesis is the behaviour of $d$-dimensional sliding minimal cones contained in the n-dimensional half-space $\mathbb{R}^n_+:=\{x_n\ge0\}$, for $0<d<n$, where the domain of the sliding boundary is the bounding hyperplane $\Gamma:=\{x_n=0\}$. In view of more general functionals to minimise (namely elliptic integrands), the considered one is not the $d$-dimensional Hausdroff measure $\mathcal{H}^d$ itself, but a small modification of it, defined by
\begin{displaymath}
J_\alpha(E):=\mathcal{H}^d(E\setminus\Gamma)+\alpha\mathcal{H}^d(E\cap\Gamma)
\end{displaymath}
for $E\subset\mathbb{R}^n_+$ (where $\alpha\in[0,1]$ for semicontinuity reasons). This energy is also related to functionals appearing in capillarity theory and free boundary problems (see e.g. \cite{finn1974capillarity} \cite{giusti1976boundary} \cite{taylor1977boundary} \cite{dephilippis2015regularity}).

\begin{figure}[h]
\centering
\includegraphics[scale=0.3]{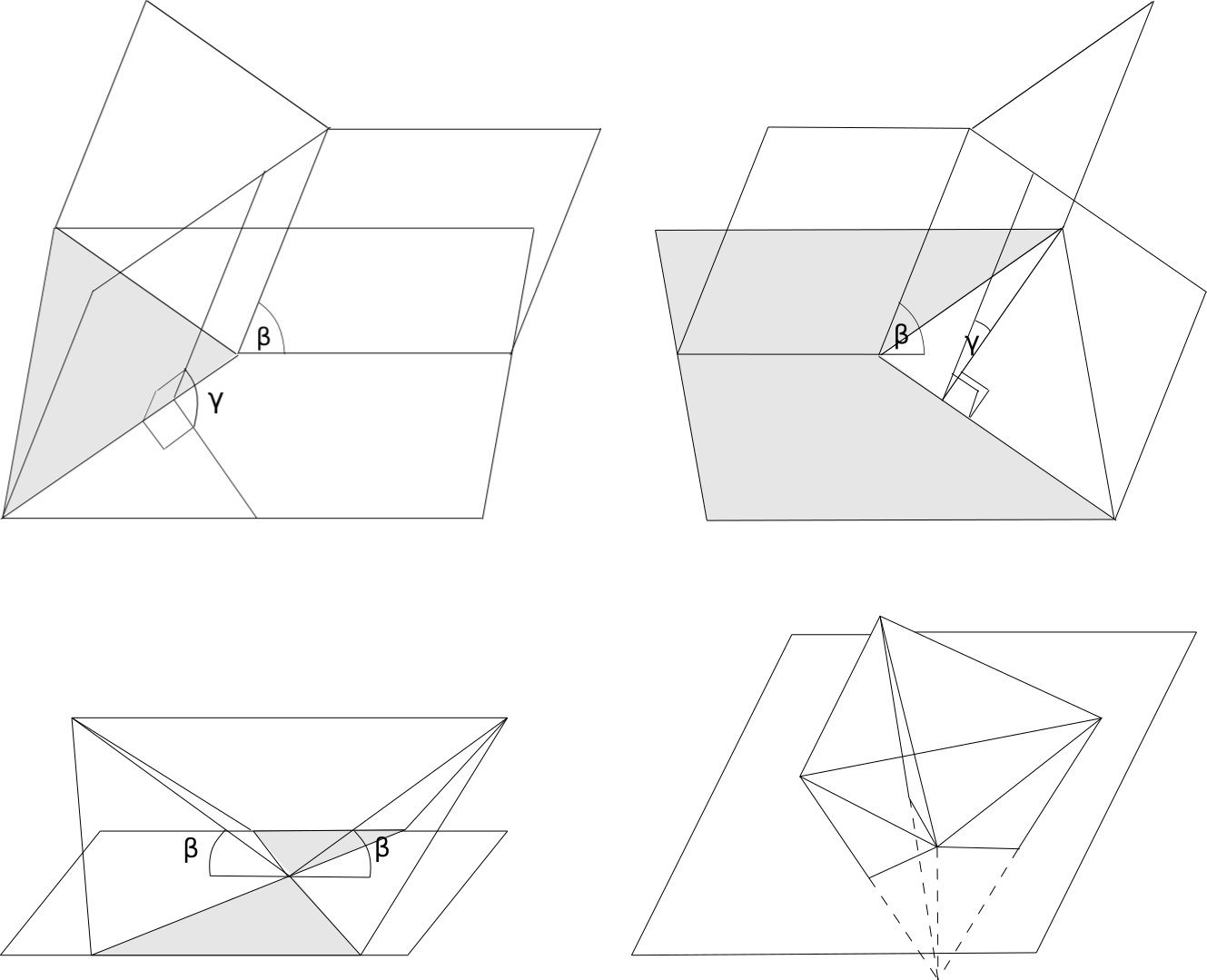}
\caption{The 2-dimensional sliding minimal cones in 3-dimensional half space, respectively: $\Y_\beta$, $\overline{\Y}_\beta$, $\W_\beta$, $\T_+$. For each one of them the grey region represent the intersection of the cone with the horizontal plane.}
\label{coniscorrevoli}
\end{figure}

Beside the cones obtained as the Cartesian product of $\mathbb{R}$ with a one-dimensional minimal cone in the half-plane, we will prove the sliding minimality of 4 new types of 2-dimensional cones in the 3-dimensional half-space (see Figure \ref{coniscorrevoli}) each of them being indeed a one-parameter family of sliding minimal cones depending on the parameter $\alpha$. In order to prove the sliding minimality of these cones we will use \emph{paired calibrations}. This tool has been employed by Lawlor and Morgan in the proof \cite{lawlor1994paired} of the minimality of the cone over the $(n-2)$-dimensional skeleton of an $n$-dimensional regular simplex; and has recently been generalised to currents with coefficients in a group by Marchese and Massaccesi in \cite{marchese2016steiner} and \cite{massaccesi2014currents}. The technique of paired calibrations consists in applying the divergence theorem to each of the connected components of the complement of a cone using a suitable family of divergence-free vector fields.  

The first two examples of minimal cones are obtained with the following procedure. First embed in $\mathbb{R}^3$ a cone of type $\mathbf{Y}$. Then tilt it in such a way that the intersection between the three half-planes meets the horizontal plane with an angle $\beta\in[0,\pi/2]$ and one of the three folds meets orthogonally the horizontal plane. Take now the intersection of this cone with the upper half-space $\mathbb{R}^3_+$. By construction the two sloping folds are forced to meet the horizontal plane with equal angle $\gamma$ (depending on $\beta$), and the vertical fold has the shape of a planar sector whose angle can be either $\beta$ or $\pi-\beta$. In the first case we define $\mathbf{Y}_\beta$ as the union of the cone obtained with the previous construction and the sector of the horizontal plane contained in between the two sloping folds. In the second case we define $\overline{\mathbf{Y}}_\beta$ as the union of the cone obtained with the previous construction and the sector of the horizontal plane not contained in between the two sloping folds. Both $\mathbf{Y}_\beta$ and $\overline{\mathbf{Y}}_\beta$ are minimal if and only if $\cos\gamma=\alpha$.

The third minimal cone is called $\mathbf{W}_\beta$ and is obtained  by taking the union between the intersection of $\overline{\mathbf{Y}}_\beta$ with a half-space bounded by a vertical plane $P$ orthogonal to the vertical fold and its reflection with respect to the plane $P$ itself. The cone $\mathbf{W}_\beta$ is minimal if and only if $\beta\le30^\circ$ and $\cos\gamma=\alpha$.

The fourth minimal cone is called $\mathbf{T}_+$ and is obtained by taking a cone of type $\mathbf{T}$ as in the first picture, flipping it upside down, placing its barycentre at the origin, and finally intersecting it with the half-space $\mathbb{R}^3_+$. Using paired calibrations it is possible to prove the minimality of $\mathbf{T}_+$ for every $\alpha\ge\sqrt{\frac{2}{3}}$. Moreover for every $\alpha<\sqrt{\frac{2}{3}}$ a better competitor to the cone can be found by pinching it down on the horizontal plane in such a way to produce a little triangle and then by connecting it to the boundary in a proper way. The previous calibration argument can be replicated in any dimension, therefore we have the following

\begin{thm*}
Given $n\ge3$ and $\alpha_n=\sqrt{\frac{n+1}{2n}}$. If $\alpha\ge\alpha_n$ then the $(n-1)$-dimensional cone over the $(n-2)$-dimensional skeleton of a regular and properly placed $n$-simplex intersected with the upper half-space $\mathbb{R}^n_+$ is a sliding minimal cone.
\end{thm*}

In Chapter \ref{minimal sets} we will introduce Almgren minimal sets, we will describe their properties and the minimal cones they produce. Then we will introduce the notion of sliding boundary and finally we will describe the setting in which we will work.

Chapter \ref{One-dimensional} will be devoted to the study of one-dimensional sliding minimal cones in the half plane.

In Chapter \ref{two-dimensional cones} we will study 2-dimensional cones in the 3-dimensional half-space and we will try to classify them with respect to their intersection with the unit sphere. In the first part of the Chapter we will discuss sliding minimal cones, starting from the ones that can be obtained as a Cartesian product and from the necessary condition that a cone has to satisfy in order to be a sliding minimal one. In the second part of the Chapter we will present a non exhaustive list of non minimal cones that satisfy the necessary condition for sliding minimality, in an attempt to classify them.

In Chapter \ref{higher dimension} we will prove the theorem stated above, generalising to higher dimension the calibration argument used to prove the sliding minimality of the cone $\T_+$ in the 3-dimensional half-space.

\chapter{Minimal sets}\label{minimal sets}

Throughout this chapter, unless otherwise specified, we will denote with $n>0$ the dimension of the ambient space, which will always be Euclidean, and with $0<d<n$ the measure theoretic dimension of the considered objects. Both $n$ and $d$ will always be integer.

In the following we will refer to the works of Almgren \cite{almgren1968existence}, David and Semmes \cite{david2000uniform}, and David \cite{david2013regularity} \cite{david2014local}, in which more general classes of minimal set are studied. The notions of restricted sets, almost minimal sets and quasi minimal sets are beyond our purpose here; however most of the results we are about to present can be extended to these objects.

\section{Area minimisers}\label{area_minimisers}

We start our discussion by presenting the notion and main properties of area minimisers, first introduced by Almgren in \cite{almgren1968existence} and then studied by David and Semmes in \cite{david2000uniform}, which represents the main reference for this section.

Let $U$ be an open subset of $\R^n$, and $E\neq\emptyset$ be a subset of $U$ which is relatively closed in $U$, that is to say $\overline{E}\setminus E\subset\R^n\setminus U$, and in particular it implies that $\overline{E}\setminus E\subset\partial U$. Assume also that $E$ has locally finite d-dimensional Hausdorff measure, that is to say $\h^d(E\cap K)<+\infty$ for every compact subset of $U$.
\begin{defn}[Admissible competitors]\label{competitori_ammissibili}
Let $E$ satisfy the previous assumptions. We say that $F$ belongs to the family of admissible competitors to $E$ in $U$ if there exists a continuous function $\phi:[0,1]\times\R^n\to\R^n$ such that (setting $\phi_t(x)=\phi(t,x)$):
\begin{enumerate}
\item $\phi_0$ is the identity;
\item $F=\phi_1(E)$;
\item $\phi_1$ is Lipschitz;
\item Set $W_t:=\{x\in\R^n:\phi(t,x)\neq x\}$ for $t\in[0,1]$ and $W=\cup_{t\in[0,1]}W_t$; then there exists a compact set $K$ such that $\phi(W)\subset K\subset U$.
\end{enumerate}
Such a function $\phi$ will also be called an admissible deformation.
\end{defn}
Since $W_t\subset\phi(W_t)$ for every $t\in[0,1]$, an admissible deformation will only affect a region of $U$ bounded away from $\R^n\setminus U$. This property of admissible deformation will play the role of a boundary condition. Another important remark is that no bound on the Lipschitz constant of $\phi_1$ is required.

\begin{defn}[Area minimiser]\label{area_minimiser}
A set $E$ is an area minimiser in $U$ if it satisfies the previous assumptions and for every admissible competitor $F$ we have that $\h^d(E\setminus F)\le\h^d(F\setminus E)$, where $\h^d$ denotes the d-dimensional Hausdorff measure (for definition and properties of the Hausdorff measure see e.g. \cite{evans1991measure}, \cite{mattila1999geometry} or \cite{ambrosio2000functions}).
\end{defn}

This means that the Hausdorff measure of $E$ cannot be decreased by any deformation that takes place ``away from the boundary''. In this sense $E$ is a solution to the Plateau problem in $U$ with respect to the boundary $S:=\overline{E}\setminus E$ (we recall that our assumptions imply $S\subset\partial U$).

First of all let us remark that, by definition, $E$ and $F$ have only locally finite Hausdorff measure, therefore it would not make sense to change the inequality in the previous definition with $\h^d(E)\le\h^d(F)$, because both sides could be $+\infty$. On the other hand, since any admissible deformation takes place in a compact set, we have that both $E\setminus F$ and $F\setminus E$ are relatively compact in $U$, therefore they have finite Hausdorff measure.

Another useful remark is that the notion of area minimiser can trivially be localised. Let $E$ be an area minimiser in $U$, and let $V$ be an open subset of $U$ such that $E\cap V\neq\emptyset$, then $E$ is an area minimiser in $V$. This follows from the fact that, being $V$ a subset of $U$, the family of admissible deformation with respect to $V$ is contained in the family of admissible deformations with respect to $U$.

Finally we remark that the set $E$ is not required to have any regularity property, like rectifiability; but, as we are about to see in next theorem, it turns out that minimality implies a much stronger regularity condition.

\begin{thm}[Almgren]\label{Almgren}
Let $E$ be an area minimiser in $U$. Then there exists a set $N$ such that $\h^d(N)=0$ and $E\setminus N$ is a $C^{1,\alpha}$ embedded submanifold of $U$ for every $0<\alpha<1$.
\end{thm}

The previous theorem was stated by Almgren in \cite{almgren1968existence} with a slightly different definition of area minimiser. In particular he required $E$ to have finite diameter and Hausdorff measure, and he didn't ask for an homotopy between the identity map and $\phi_1$. However the two notions are equivalent for our purpose, and this is because, up to a localisation, they coincide. To see this it is enough to assume $U$ to be bounded and convex: the boundedness condition forces $E$ to have finite diameter and Hausdorff area, and the convexity allows us to choose the affine homotopy between the identity map and $\phi_1$, that is to say $\phi(t,x)=(1-t)x+t\phi_1(x)$.

\section{Minimal cones}\label{Minimal_cones}

Theorem \ref{Almgren} is a strong regularity result, characterising minimal sets almost everywhere. Because of this result we have that minimal sets admit almost everywhere a flat tangent space in the classical sense, therefore, in order to get a complete local description of minimal sets one still needs to know what happens in the points of the exceptional set $N$.

It can be proved that any blow-up of a minimal set at any of its point is a dilation invariant set which is itself an area minimiser on the whole space $\R^n$ (that is to say $U=\R^n$ in Definition \ref{area_minimiser}) see \cite{david2000uniform}. Dilation invariant sets will be called cones, and will always be supposed to be centred at the origin. Cones that arise from blow-ups will be called tangent cones. Therefore knowing minimal cones is an important step in order to understand the local behaviour of minimal sets.

Since cones are dilation invariant unbounded set, and admissible deformation only take place in a compact set, in order to prove that a cone is minimal in $\R^n$ one only has to prove that the cone is minimal in a given open and bounded neighbourhood of the origin. This remark will allow us to pick each time the most suitable open set in which to prove the minimality of a cone.

There are only two types of 1-dimensional minimal cones in $\R^2$, and they are: straight lines, and the cone called Y, formed by three half lines meeting with equal angles of $120^\circ$. These two are the only 1-dimensional minimal cones also in any higher dimensional ambient space.

In $\R^3$ there are three types of 2-dimensional minimal cones: the planes, the sets $\mathbf{Y}:=Y\times\R$, and $\mathbf{T}$, the cone over the 1-dimensional skeleton of a regular tetrahedron centred in the origin. Given an n-dimensional polyhedron $P$ we will refer to the union of all its m-dimensional faces as the m-dimensional skeleton of $P$, denoted by $\sk_m(P)$. Given any subset $A$ of $\R^n$ we define the cone over $A$, sometimes written as $cone(A)$, as the smallest cone containing the set $A$, or, equivalently, as the union of all the half lines starting from the origin and intersecting $A$.

\begin{figure}
\centering
\includegraphics[scale=0.37]{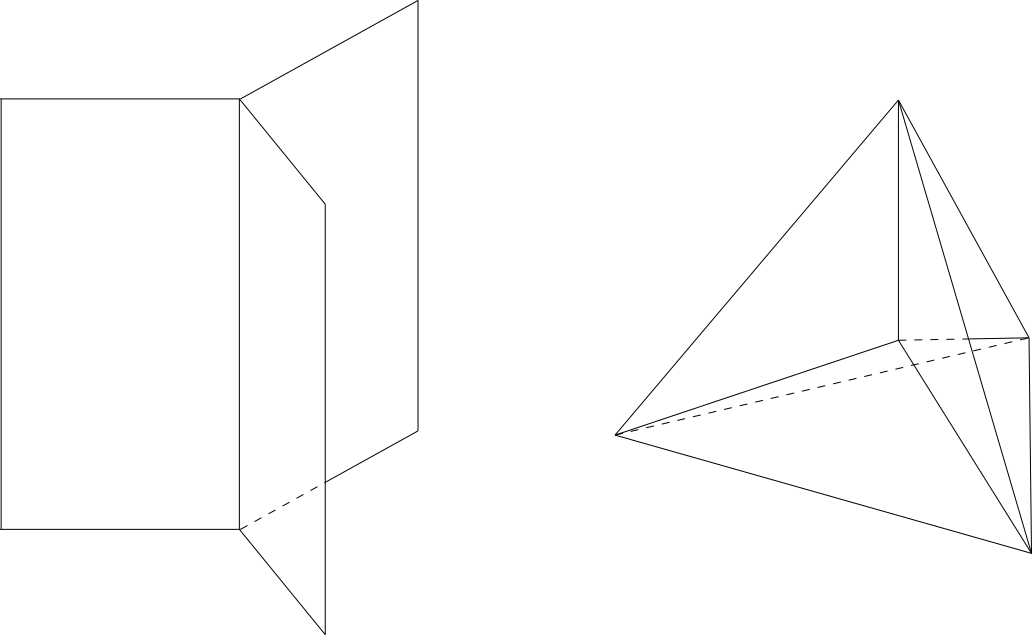}
\caption{Cones of type $\mathbf{Y}$ and $\mathbf{T}$.}
\label{YeT.png}
\end{figure}

These minimal cones are known since a long time ago, but in \cite{taylor1976structure} Taylor proved that there are no others. In the same paper she also proved that minimal sets of dimension 2 in $\R^3$ are locally equivalent, through $C^1$ diffeomorphisms, to minimal cones. That is to say, given $x\in E$ such that $\mathbf{C}$ is the blow-up limit of $E$ at $x$ (where $\mathbf{C}$ could be a plane, $\mathbf{Y}$, or $\mathbf{T}$), there exist a radius $r>0$ and a $C^1$ diffeomorphism $\varphi:B^n(0,1)\to B^n(x,r)$ such that $\varphi(\mathbf{C}\cap B^n(0,1))=E\cap B^n(x,r)$.

In higher dimension we only know some examples of minimal cones, and the list is incomplete. Even in dimension 4 we don't know whether the known examples provide the complete list of minimal cones. In \cite{brakke1991minimal} Brakke proved that the cone over the (n-2)-dimensional skeleton of a top-dimensional cube in $\R^n$ is a minimal cone (of codimension 1) if and only if $n\ge4$.

\begin{figure}
\centering
\includegraphics[scale=0.4]{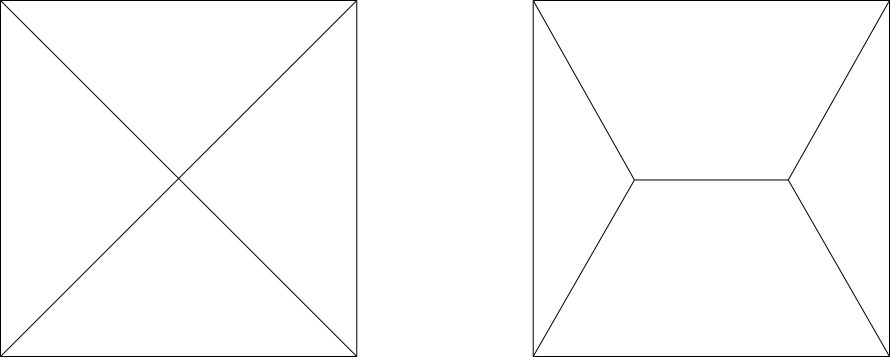}
\caption{Cone over the square on the left and a better competitor on the right.}
\label{conoquadrato}
\end{figure}

When $n=2$, a better competitor for the cone over the vertices of a square can be easily provided by pinching together the four branches in such a way as to produce two branching point whose tangent cone is a $Y$, as in picture \ref{conoquadrato}.

\begin{figure}
\centering
\includegraphics[scale=0.25]{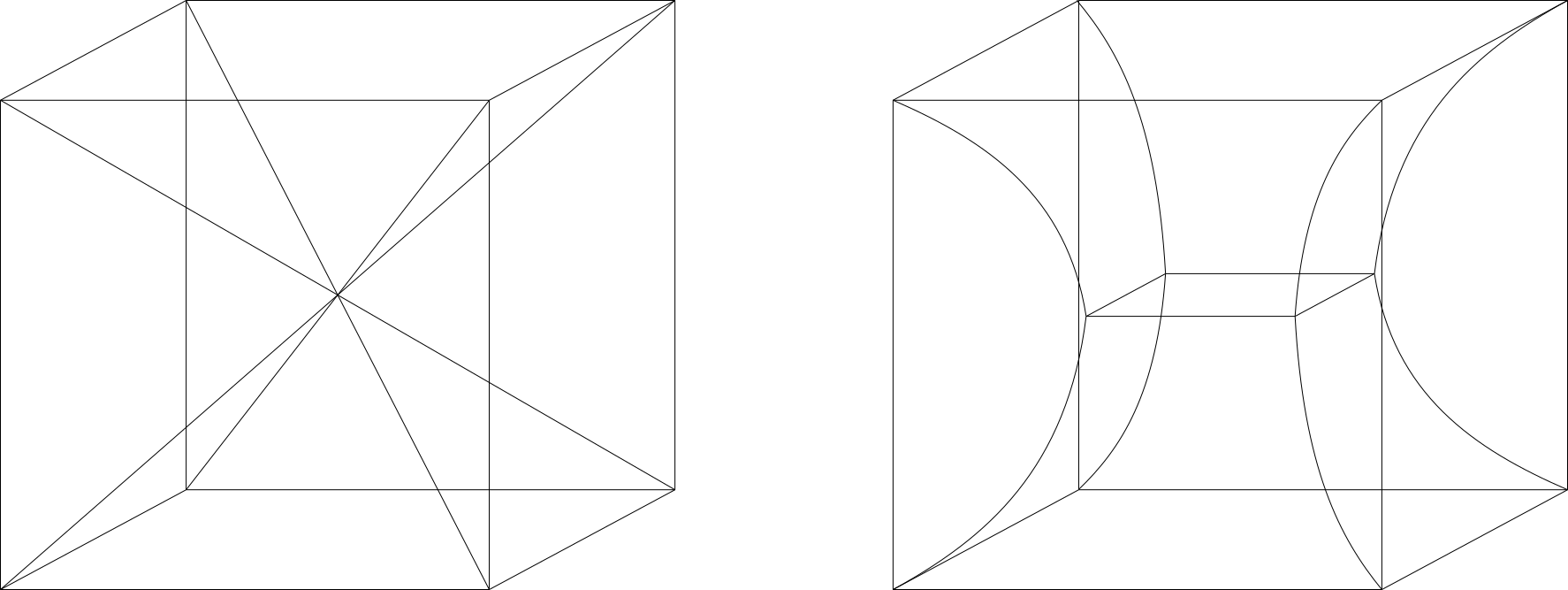}
\caption{Cone over the edges of a cube on the left and a better competitor on the right.}
\label{conocubo}
\end{figure}

The same can be done with the cone over the edges of a cube in $\R^3$. The sloping faces can be pinched together at the center of the cube in such a way as to produce a little square. Brakke showed that the sloping faces can be bent along a suitable curve in order to have a better competitor than the cone, see picture \ref{conocubo}. Moreover he proved that this kind of competitor is still better even when the square interface is more expensive, in term of area, by a parameter $1\le T\le\sqrt{2}$. The critical value of $T$, below which the cone is minimal, is smaller than one in dimension 4, and decreases for higher dimensions. From dimension 7 and higher the cone is still minimal even when the cost of the interface between two opposite regions is 0.

In \cite{lawlor1994paired} Lawlor and Morgan proved that the cone over the $(n-2)$-dimensional skeleton of a top-dimensional regular simplex in $\R^n$ is a minimal cone (again the codimension is 1) for any $n\ge2$. These cones will be called $\mathbf{\Delta}^n$, let us stress the fact that $n$ is the dimension of the ambient space, while the dimension of the cone is $n-1$. It follows that $\mathbf{\Delta}^2=Y$ and $\mathbf{\Delta}^3=\mathbf{T}$. The proof of this result relies on a technique called `paired calibration', and we report it here below because, in the author's opinion, it is very elegant and also because it will be useful to our purpose later.

Let $\Delta^n:=[p_0,...,p_n]$ be the n-dimensional simplex whose n+1 vertices are $p_0,...p_n\in\R^n$; assume it to be regular and centred on the origin. For $0\le k\le n$, we set $\Delta^n_k$ to be the $(n-1)$-dimensional simplex whose vertices are the previous ones except for $p_k$. Analogously, for $0\le k\neq h\le n$, we will denote with $\Delta^n_{k,h}$ the $(n-2)$-dimensional simplex whose vertices are the previous ones except for $p_k$ and $p_h$. The key observation for the proof is that, for $0\le i<j\le n$, the vector $p_i-p_j$ is orthogonal to the fold of the cone over the facet $\Delta^n_{i,j}$; more compactly $(p_i-p_j)\perp cone(\Delta^n_{i,j})$, a proof of this fact can be found in Appendix \ref{proprieta_simplessi}. As we are about to show, the calibration argument relies on this remark and on the divergence theorem.

As we mentioned before, it is enough to prove that the cone is minimal with respect to admissible deformations acting in a fixed bounded neighbourhood of the origin. We choose this neighbourhood to be the open simplex $\Delta^n$ itself. For every $i=0,\dots,n$ we set $F_i:=\Delta^n_i$, which is the $(n-1)$-dimensional face of the simplex opposed to the vertex $p_i$. Let $M$ be an admissible competitor to $\mathbf{\Delta}^n$ in $\R^n$, obtained by one of the deformations mentioned before, in such a way that the symmetric difference between the cone and the competitor is compactly contained in the simplex $\Delta^n$. The competitor is countably $(\h^{n-1},n-1)$ rectifiable since it is the image of $\mathbf{\Delta}^n$ via a Lipschitz map, therefore its intersection with the chosen simplex is $(\h^{n-1},n-1)$ rectifiable because it has finite $(n-1)$-dimensional Hausdorfff measure.

Hence $\R^n\setminus M$ has $n+1$ unbounded connected components (cf. \cite[XVII 4.3]{dugundji1970topology}), each one of them having locally finite perimeter and containing one face of the simplex $\Delta^n$. For $i=0,...,n$ we name $V_i$ the connected component of $\R^n\setminus M$ containing $F_i$ and in case $\R^n\setminus M$ also has bounded connected components we just include them in $V_0$. Now let us set $U_i:=V_i\cap\Delta^n$. The $U_i$ are finite perimeter sets, hence we can define $M_i:=\partial^*U_i\setminus F_i$ (where $\partial^*U_i$ denotes the reduced boundary of $U_i$) and $n_i$ as the exterior unit normal to $\partial^*U_i$. Then $M_i\subset M\cap\Delta^n$, and, by definition of reduced boundary, it follows that the blow-up limit of $M_i$ at any of its point is a hyperplane separating its complement in exactly two connected component, each one of them being a half-space (respectively the blow-up limit of the set itself and the blow-up limit of it complement).

For $0\le i\neq j\le n$ let us define $M_{ij}:=M_i\cap M_j=\partial^*U_i\cap\partial^*U_j$ and $n_{ij}$ as the unit normal to $M_{ij}$ pointing in direction of $U_j$. The sets $M_{ij}$ are contained in $M$ and in particular $M_{ij}$ is contained in the interface between the regions $U_i$ and $U_j$. Let us now remark that for every $i=0,...,n$, $\H^{n-1}$-almost every point of $M_i$ lies on the interface between exactly two regions $U_i$ and $U_j$. Therefore the interfaces between different couples of regions are essentially disjoint with respect to $\H^{n-1}$ and
\begin{equation}
\H^{n-1}(M_i)=\H^{n-1}\left(\bigcup_{j\neq i}M_{ij}\right)=\sum_{j\neq i}\H^{n-1}\left(M_{ij}\right).
\end{equation}
In order to show it let us first define, for $i=0,...,n$, the exceptional set
\begin{equation}
E_i:=\R^n\setminus\left(U_i^0\cup\partial^*U_i\cup U_i^1\right),
\end{equation}
where, for an $\mathcal{L}^n$-measurable set $A\subset\R^n$ and $t\in[0,1]$, we define
\begin{equation}
A^t:=\left\{x\in\R^n:\lim_{r\to0}\frac{\mathcal{L}^n(A\cap B_r(x))}{\mathcal{L}^n( B_r(x))}=t\right\}.
\end{equation}
Since $U_i$ is a finite perimeter set for every $i=0,...,n$, by Federer's theorem $\H^{n-1}(E_i)=0$ (see \cite[Theorem 3.61]{ambrosio2000functions}). Therefore the exceptional set $E:=\cup_iE_i$ is negligible with respect to $\H^{n-1}$. Let us now assume that a point $x$ belongs to the common boundaries of at least three sets $U_i$, $U_j$ and $U_k$. Clearly the point $x$ cannot belong to the reduced boundaries of the three of them because in this case the blow-up limit of each one of them in the point $x$ would be a half-space and that is a contradiction. Let us assume $x\notin\partial^*U_i$, then $x\in E_i\subset E$. Since $E$ is $\H^{n-1}$-negligible, the same holds true for its intersection with $M_i$ for any $i=0,...,n$. 

Let us now define $\widetilde{M}:=\cup_iM_i$. It follows that $\widetilde{M}$ is contained in $M\cap\Delta^n$ and $\h^{n-1}$-almost every point in $\widetilde{M}$ lies on the interface between exactly two connected components of its complement. Now we define the vectors of the calibration as $w_i:=\frac{-p_i}{|p_i-p_j|}$, for $0\le i\le n$ and $j\neq i$. Recalling that $w_i\perp F_i$ we can compute as follows
\begin{displaymath}
\begin{aligned}
|w_i|\h^{n-1}(\partial\Delta^n) & =\sum_{i=0}^n\int_{F_i}w_i\cdot n_id\h^{n-1}\stackrel{(i)}{=}-\sum_{i=0}^n\int_{M_i}w_i\cdot n_id\h^{n-1}\\
&\stackrel{(ii)}{=}-\sum_{i\neq j}\int_{M_{ij}}\!\!w_i\cdot n_{ij}d\h^{n-1}=\sum_{i<j}\int_{M_{ij}}\!\!(w_j-w_i)\cdot n_{ij}d\h^{n-1}\\
&\stackrel{(iii)}{\le} \sum_{i<j}\h^{n-1}(M_{ij})=\h^{n-1}(\widetilde{M})\le\h^{n-1}(M),
\end{aligned}
\end{displaymath}
where in $(i)$ we applied the divergence theorem to the constant vectors $w_i$ on the finite perimeter sets $U_i$, in $(ii)$ we used the definitions and the properties of $M_{ij}$ and $n_{ij}$, and in the inequality $(iii)$ we employed the fact that $|w_i-w_j|=1$. Since the first term on the left is a constant, what we obtained is a lower bound for the area of a competitor. By choosing the cone itself as a competitor we get a chain of equalities since in this case $(w_j-w_i)=n_{ij}$ (because $(w_i-w_j)\perp cone(\Delta^n_{ij})$ as we mentioned above) and $M=\widetilde{M}$. Hence the lower bound is attained by the cone, and this means the cone is an area minimiser.

In \cite{massaccesi2014currents} Massaccesi proved that the cone $\mathbf{\Delta}^n$ supports a current whose coefficients are carefully chosen in an appropriate group. Moreover that current is mass minimising among all the currents with the same boundary, and its mass coincides with the Hausdorff area of the cone. In the author's opinion it is very interesting how such a current was constructed, and how the proof of the mass minimality turned out to be very similar to Morgan's paired calibration (as Massaccesi herself remarked).

The other known examples of higher dimensional minimal cones are due to Liang. In \cite{liang2012almgren} she proved that the union of two almost orthogonal planes in $\R^4$ is minimal, providing for the first time an example of a one-parameter family of minimal cones that are not isometric to each other. The particular case where the planes are orthogonal can be proved with calibrations and relies on the fact that, given $(x_1,x_2,x_3,x_4)$ an orthonormal base of $\R^4$, the mass of the two-vector $x_1\wedge x_2\pm x_3\wedge x_4$ is exactly one (see \cite{liang2014almgren}).

In \cite{liang2014almgren} Liang proved the minimality of $Y\times Y\subset\R^4$, the two-dimensional cone obtained by the Cartesian product of two one-dimensional cones of type $Y$, each one of them contained in a different copy of $\R^2$. The proof is by calibration and in this case the calibration is the family of two-vectors $v_i\wedge w_j$, for $i,j=0,1,2$, where $v_i$ and $w_j$ are the calibrations of the two different copies of $Y$. The main difference with $\mathbf{\Delta}^n$ is that we are in codimension higher than one. In codimension one $\h^d$-almost every point on an admissible competitor lies on the interface between exactly two connected components of its complement, and that means that exactly two vectors of the calibration act on it. In higher codimension different folds of a competitor are allowed, in principle, to overlap multiple times, therefore, in order to get an inequality as in $(iii)$ above, one has to check all possible cases that may arise, taking care of all the possible multiplicities.

Finally in \cite{liang2015topological} Liang proved a very general result about the minimality of the union of subspaces of arbitrary codimension, generalising her result about almost orthogonal planes in $\R^4$. The statement is the following. For each $d\ge2$ and $m\ge2$, there exists an angle $\theta_{m,d}\in(0,\frac{\pi}{2})$, such that, given a family $P_1,...,P_m$ of d-dimensional subspaces of $\R^{dm}$, their union $\cup_{i=1}^mP_i$ is minimal in $\R^{dm}$ if all the characteristic angles between any two subspaces are bigger than $\theta_{m,d}$.

\section{Sliding boundary}

In this section we will discuss the notion of sliding boundary, introduced by David in \cite{david2013regularity} and \cite{david2014should} with the purpose to extend onto the boundary the regularity results stated above. As we did before, given a set $E\subset\R^n$ with locally finite $d$-dimensional Hausdorff measure, we will first define a new class of competitors, here called sliding competitors, an then we will say that $E$ is a sliding minimiser if it minimises $\h^d$ among the sliding competitors. Let $\Omega$ be a closed subset of $\R^n$, which may coincide with $\R^n$ itself; and let $\Gamma_i\subset\Omega$, for $0\le i\le I$, be finite family of closed sets. For the sake of notation we will set $\Gamma_0:=\Omega$.

\begin{defn}[Sliding competitor]\label{sliding_competitor}
We say that $F$ is an admissible competitor to $E$ in $\Omega$, with respect to the sliding boundary domains $\{\Gamma_i\}_{0\le i\le I}$, if there exists a continuous function $\phi:[0,1]\times E\to\mathbb{R}^n$ such that (setting $\phi_t(x)=\phi(t,x)$):
\begin{enumerate}
\item $\phi_0$ is the identity;
\item $F=\phi_1(E)$;
\item $\phi_1$ is Lipschitz;
\item Set $W_t:=\{x\in E:\phi(t,x)\neq x\}$ for $t\in{0,1}$ and $W=\cup_{t\in[0,1]}W_t$, there exists a compact set $K$ such that $\phi(W)\subset K\subset \R^n$.
\item for any $0\le i\le I$, if $x\in E\cap\Gamma_i$ then $\varphi_t(x)\in \Gamma_i$ $\forall t\in[0,1]$.
\end{enumerate}
\end{defn}

For short we will refer to this class of competitor simply as sliding competitor, and to the deformation satisfying the definition as sliding deformations. The requirement 5 is new with respect to Definition \ref{competitori_ammissibili} and encodes the sliding boundary condition, but there are two more differences that we want to talk about. First, the deformation in Definition \ref{sliding_competitor} is only defined on $E$, and, in general, it is a non-trivial problem to extend it onto $\R^n$ preserving the sliding boundary condition. However we will show later that this can be done in the cases we are concerned with. Second, $\Omega$ is closed while $U$ in Definition \ref{competitori_ammissibili} is open; as we will see later, this little difference will make the definition more suitable for the case in which $E$ stays on one side of the sliding boundary. Finally let us remark that this notion can be localised to an open set $U$ as we did for Almgren minimal sets.

\begin{defn}[Sliding minimal set]\label{minimo_scorrevole}
A set $E$ is a minimiser in $\Omega$ with respect to the sliding boundary $\{\Gamma_i\}_{0\le i\le I}$, if, for any sliding competitor $F$, we have that $\h^d(E\setminus F)\le \h^d(F\setminus E)$.
\end{defn}

In \cite{david2014local} David proved that under some mild regularity condition of the boundary, the blow-up limit of a sliding minimal set at a boundary point is a sliding minimal cone with respect to a conical boundary. Since the assumption of more regularity of the boundary (like $C^1$ or rectifiable) provides the flatness (everywhere or almost-everywhere) of its blow-up, an important step, in order to understand the behaviour of sliding minimisers close to the boundary, is to know the list of sliding minimal cones with respect to flat boundaries.

In \cite{fang2016holder} Fang proved that given $\Gamma$ a two-dimensional $C^1$ submanifold of $\R^3$, and a two-dimensional set $E$ which is sliding minimal with respect to $\Gamma$, then $E$ is locally biHölder equivalent to a sliding minimal cone under the assumption that $E$ contains $\Gamma$ and stays on one side of it. As a first step of the proof he provided a complete list of one-dimensional sliding minimal cones in the half-plane with respect to the bounding straight line, and a partial list of two-dimensional minimal cones in the three-dimensional half-space with respect to the bounding plane. In the recent paper \cite{fang2017local} Fang improved this result to $C^{1,\alpha}$-regularity at the boundary of sliding minimal sets. We expect a similar result (at least the Hölder part) to hold with the corresponding analogue of our basic problem when $0\le\alpha\le1$ and follow from our description of minimal cones.

\section{General integrands}

Up to now we used the Hausdorff measure to weight the size of our sets, but one could try to minimise more general functionals.

Let $f:\R^n\times G(d,n)\to(0,+\infty)$ be a Borel-measurable positive function, where $G(d,n)$ denotes the Grassmannian manifold of unoriented d-planes in $\R^n$. Such a function $f$ will be referred to as a positive integrand. Given a rectifiable set $E$ we define the following functional
\begin{displaymath}
J_f(E):=\int_Ef(x,T_xE)d\h^d(x),
\end{displaymath}
where $T_xE$ denotes the non oriented d-plane which gives the approximate tangent
plane to $E$ at $x$; thus $T_xE$ is defined $\h^d$-almost everywhere on $E$ because $E$ is rectifiable.

\begin{defn}[$f$-minimiser]\label{f-minimo}
Under the same assumptions as in Definition \ref{area_minimiser}, and given a positive integrand $f$, a rectifiable set $E$ is said to be an $f$-minimiser in $U$ if for every admissible competitor $F$ we have that $J_f(E\setminus F)\le J_f(F\setminus E)$.
\end{defn}

\begin{defn}[Sliding $f$-minimal set]\label{f-minimo_scorrevole}
Under the same assumptions as in Definition \ref{minimo_scorrevole}, and given a positive integrand $f$, a rectifiable set $E$ is an $f$-minimiser in $\Omega$ with respect to the sliding boundary $\{\Gamma_i\}_{0\le i\le I}$, if, for any sliding competitor $F$, we have that $J_f(E\setminus F)\le J_f(F\setminus E)$.
\end{defn}

The aforementioned regularity results provided by Almgren and David for area minimisers and sliding minimal sets have actually been proven for $f$-minimisers and sliding $f$-minimal sets when $J_f$ belongs to the class of \textit{elliptic integrands} (which includes the area functional $f\equiv 1$, see \cite{almgren1968existence}).

\section{Our setting}

In the following the subject of our investigation will be blow-up limits of sliding minimal sets arising at boundary points. In particular we will focus on the case where the cones are contained in a closed half-space $\R^n_+:=\{(x_1,\dots,x_n)\in\R^n:x_n\ge0\}$ and the domain of the sliding boundary is the hyperplane $\Gamma:=\{(x_1,\dots,x_n)\in\R^n:x_n=0\}$ bounding the half-space. Therefore, using the notation introduced in the previous section we have $\Gamma_0=\Omega=\R^n_+$ and $\Gamma_1=\partial\Omega=\Gamma$.

In view of the general class of elliptic integrands, the functional we will try to minimise will not be the $d$-dimensional Hausdroff measure $\mathcal{H}^d$ but a small modification of it, defined as follows

\begin{defn}[Modified functional]
Let $\Gamma$ be the hyperplane defined above, and $\alpha\in[0,1]$ be a constant. For $E\subset\mathbb{R}^n_+$ we define
\begin{displaymath}
J_\alpha(E):=\mathcal{H}^d(E\setminus\Gamma)+\alpha\mathcal{H}^d(E\cap\Gamma).
\end{displaymath}
\end{defn}

This functional is indeed a weighted Hausdorff measure and the constant  $\alpha$ is chosen to take values in the interval $[0,1]$ for semicontinuity reasons. In the following we may refer to this functional simply as `energy' or `cost' and we will call a set $E$ negligible if $J_\alpha(E)=0$. The effect produced by this modification will be that for a surface, in order to minimise the cost, it will often be more convenient to lie on $\Gamma$ rather than in the rest of $\R^n_+$. Since $\Gamma$ is both the domain of the sliding boundary and the region where the Hausdorff area has a different weight, this setting can be viewed as a particular case of Plateau's Problem, a case where there is a physical interaction between the soap bubble and the material which the sliding boundary is made of.

In the following, unless otherwise specified, we will refer to a sliding $J_\alpha$-minimal set as $\alpha$-minimal, or simply as minimal in case the value of the constant $\alpha$ is clear from the context.

\chapter{One-dimensional cones}\label{One-dimensional}

\section{Half-plane}\label{Half-plane}
In this section we will discuss one-dimensional minimal cones in the half-plane, with respect to the bounding line. Using the notation introduced in the previous chapter we have that $\Omega=\mathbb{R}^2_+=\{(x,y)\in\mathbb{R}^2:y\ge0\}$ and $\Gamma=\{y=0\}$.

Given $0\le\theta\le\frac{\pi}{2}$, let $P_\theta$ be a half-line meeting $\Gamma$ at the origin with angle $\theta\in[0,\frac{\pi}{2}]$, and let $\theta_\alpha$ be such that $\alpha=\cos\theta_\alpha$. In particular let us remark that $\theta_\alpha\to\pi/2$ when $\alpha\to0$ and $\theta_\alpha\to0$ when $\alpha\to1$. Given $0\le\alpha\le1$ the minimal cones are the following (see Figure \ref{coni1in2}):

\begin{figure}
\centering
\includegraphics[scale=0.35]{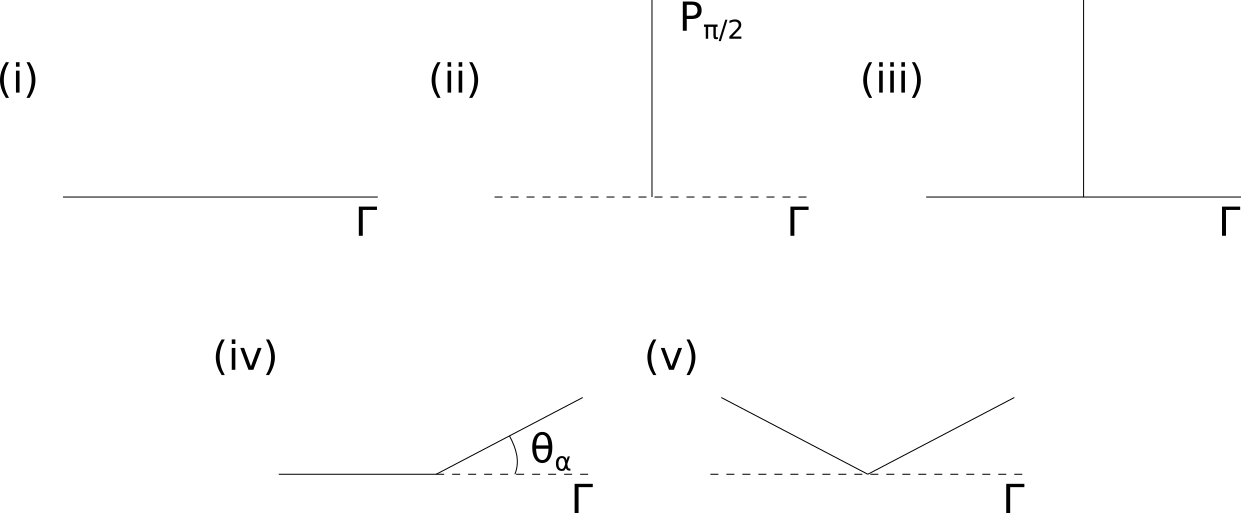}
\caption{One dimensional sliding $\alpha$-minimal cones in the half plane.}
\label{coni1in2}
\end{figure}

\begin{figure}
\centering
\includegraphics[scale=0.35]{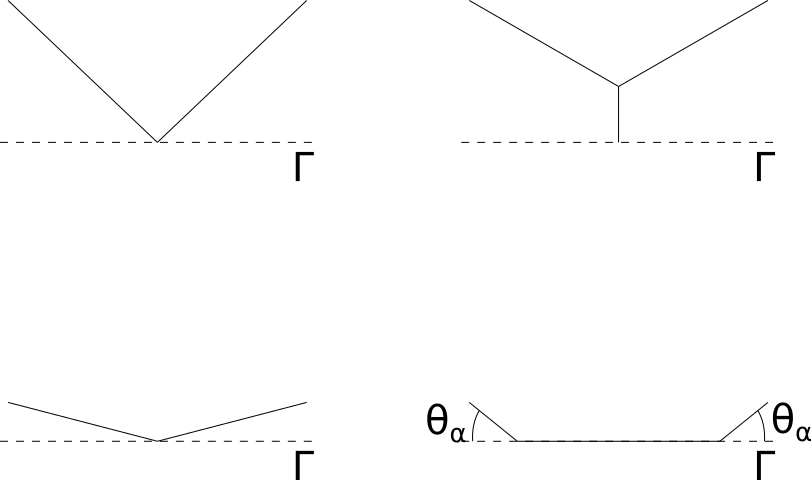}
\caption{Better competitors to the cone (v) in case $\theta>\pi/6$ above, and in case $\theta\le\theta_\alpha$ below.}
\label{competitori_cono_v}
\end{figure}

\begin{itemize}
\item[(i)] $\Gamma$; this cone is trivially minimal because the only member in the family of sliding competitors is $\Gamma$ itself.
\item[(ii)] $P_{\frac{\pi}{2}}$; given a compact set $K$ and $x_0\in P_{\frac{\pi}{2}}\setminus K$, any sliding competitor to the vertical half-line would be (or at least contain) a path connecting $x_0$ to $\Gamma$, therefore it would be longer than the vertical line segment connecting $x_0$ to the origin.
\item[(iii)] $\Gamma\cup P_{\frac{\pi}{2}}$; we can show the minimality of this cone combining the two previous arguments. In fact any competitor to this cone still contains $\Gamma$, and given $x_0$ as before we have that any competitor contains a path connecting $x_0$ to $\Gamma$.
\item[(iv)] the union of $P_{\theta_\alpha}$ with a horizontal half-line laying on $\Gamma$; let $B_1(0)$ be the unit ball centred at the origin, and $\theta\in[0,\frac{\pi}{2}]$. We define $A=(-1,0)$ and $B=(\cos(\theta),\sin(\theta))$ to be the two endpoints of the cone intersected with $B_1(0)$. In order to prove the minimality of the cone it is sufficient to consider all competitors obtained as the union of the two segments $\overline{AC}$ and $\overline{CB}$ where $C=(x,0)$. Hence we have to minimise $J_\alpha$ among a one-parameter family of competitors. Let $E_x$ be one of such competitors, than
\begin{equation}
\begin{aligned}
J_\alpha(E_x) & =\alpha(1+x)+\sqrt{(x-\cos(\theta))^2+\sin^2(\theta)}\\
\left.\frac{\partial J_\alpha(E_x)}{\partial x}\right|_{x=0} & =\left[\alpha+\frac{x-\cos(\theta)}{\sqrt{(x-\cos(\theta))^2+\sin^2(\theta)}}\right]_{x=0}= \alpha-\cos(\theta)\\
\left.\frac{\partial^2 J_\alpha(E_x)}{\partial x^2}\right|_{x=0} & = \left.\frac{\sin^2(\theta)}{\left[(x-\cos(\theta))^2+\sin^2(\theta)\right]^{3/2}}\right|_{x=0}=\sin^2(\theta).
\end{aligned}
\end{equation}
Therefore $x=0$ is a critical point if and only if $\cos(\theta)=\alpha$, and the second derivative is always positive.
\item[(v)] $V_\theta$: the union of $P_\theta$ and its symmetric with respect to the vertical axis, for $\theta_\alpha\le\theta\le\pi/6$. As before let us call $A$ and $B$ the endpoints of $V_\theta$ intersected with the ball $B_1(0)$. Any admissible competitor for $V_\theta$ has to connect $A$, $B$ and $\Gamma$. It is easily seen that $V_\theta$ is minimal among all the competitors obtained as the union of the two segments $\overline{AC}$ and $\overline{CB}$ where $C=(x,0)$. However other types of competitor may occur (see picture \ref{competitori_cono_v}). Pinching together the two segments $\overline{AC}$ and $\overline{CB}$ we can produce a triple junction, and in case $\theta>\pi/6$ it is possible to arrange it in the shape of a $Y$ cone. In this case the obtained competitor is actually the minimiser. Otherwise we can push down the two segments onto $\Gamma$ producing a segment $\overline{C'C''}$ in $\Gamma$. In case $\theta<\theta_\alpha$ we can keep pushing down up to the point when the angles formed by $\overline{AC'}$ and $\overline{C''B}$ with $\Gamma$ turn into $\theta_\alpha$, and again we obtain a minimiser.
\end{itemize}

We remark that cones (i), (ii) and (iii) are independent from $\alpha$ while cone (iv) forms a one parameter family depending on this value. Cones of type (v) form a one-parameter family independent on $\alpha$ but whose endpoints depend on it. When $\alpha=1$ we have that $J_1=\h^1\llcorner\Omega$ and the cone (iv) collapses to the cone (i). On the other side, when $\alpha=0$, we have that $J_0=\h^1\llcorner(\Omega\setminus\Gamma)$. In this case (i) turns into an (even more) trivial minimal cone, and the cones (ii) and (iii) become equivalent with respect to $J_0$ because they only differ from a ``null-measure'' set. Moreover also the cone (iv) collapses to the type (ii)-(iii).

The minimality of the cones (i) and (iv) is only due to the definition $J_\alpha$ and it would still be minimal without imposing the sliding boundary condition. The cones (ii) and (v) are in the opposite situation: they would not be minimal without the sliding boundary condition, regardless to the coefficient $\alpha$. Cone (iii) in an exceptional case, its minimality relies on the cost functional when $\alpha\le\sqrt{3}/2$ and on the sliding boundary condition when $\alpha>\sqrt{3}/2$. Indeed, in the latter case, what could happen if we drop the sliding boundary condition is that the branching point of the cone could move upwards assuming the $Y$ configuration. Therefore the two branches of this $Y$ would meet the boundary with an angle of $30^\circ$, whose cosine is $\sqrt{3}/2$.

\begin{figure}
\centering
\includegraphics[scale=0.35]{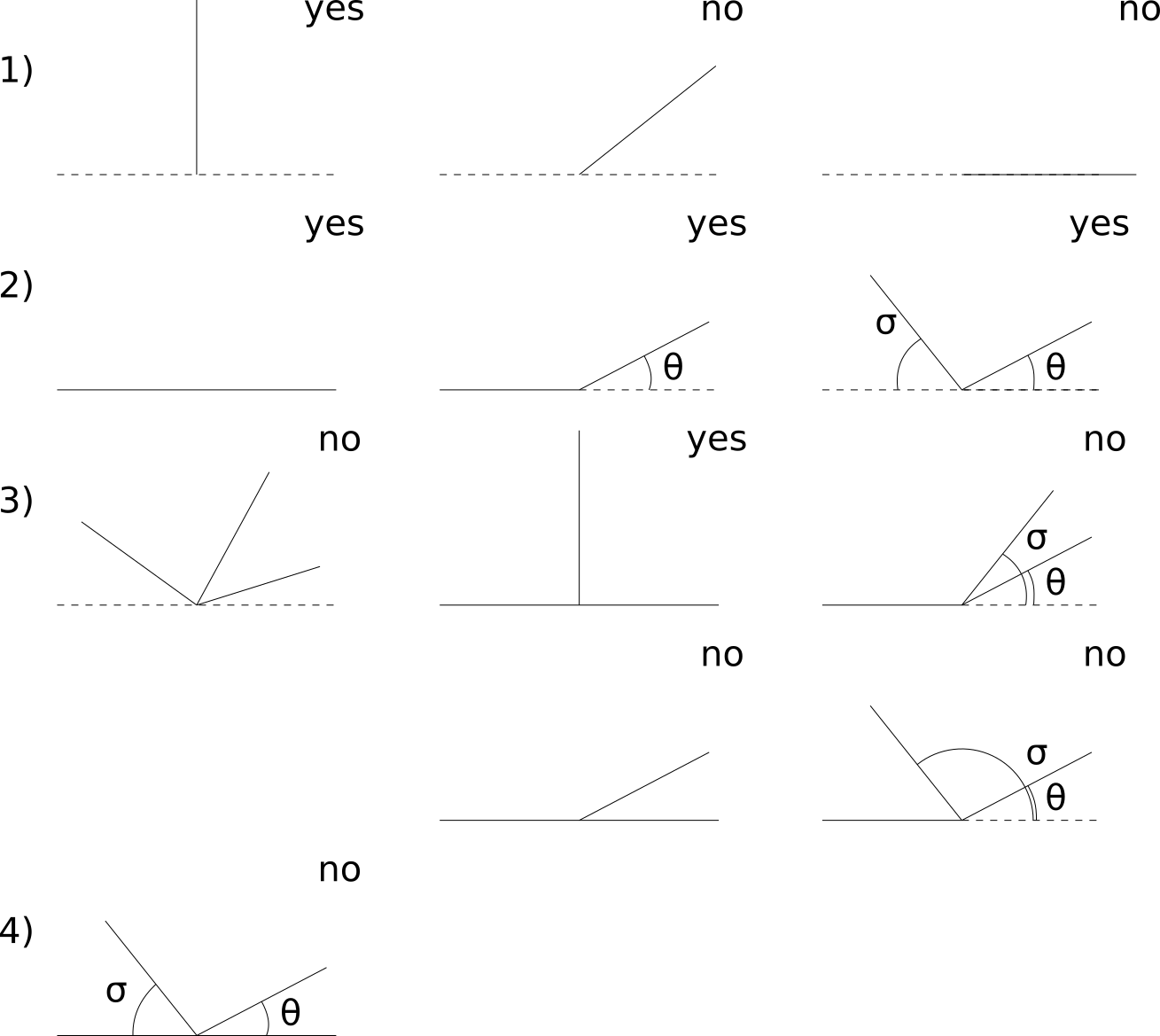}
\caption{Classification of one-dimensional cones by the number of their branches, and whether or not there exists a minimal cone of a given type.}
\label{coni1in2-nonminini}
\end{figure}

In order to show that the aforementioned list of sliding minimal cones is complete we can classify all the one-dimensional cones by the number of distinct half lines (or branches) they are composed by (see Figure \ref{coni1in2-nonminini}) and by their position with respect to $\Gamma$. In case the cone is composed by only one branch than it is sliding minimal only if the branch is vertical (cone of type (ii)); otherwise it is very easy to find a better competitor. In case the cone has two branches than we have three sub-cases, depending on the number of branches contained in $\Gamma$. If they are both contained in $\Gamma$ we find again the cone of type (i); if only one is contained then the cone if minimal if and only if it is of type (iv); and it both the branches are not contained in $\Gamma$ then the cone is minimal if and only if it is of type (v). Let us now discuss the case of three branches. If the three of them are not contained in $\Gamma$ than at least two of them form an angle smaller then $120^\circ$ therefore they can be pinched together decreasing the total length. If only one branch is not contained in $\Gamma$ than the cone is minimal if and only if it is of type (iii) otherwise the sloping branch can be projected onto $\Gamma$ decreasing the total energy. If exactly one branch is contained in $\Gamma$ we can call $\sigma$ and $\theta$ the angles formed by the sloping branches with $\Gamma$ (as in Figure \ref{coni1in2-nonminini}). Than we have two sub-sub-cases. If both the angles are less than a right angle than the angle between them is less than $120^\circ$ and they can be pinched together. If at least one angle is bigger than a right angle than the corresponding branch can be pushed down onto $\Gamma$ in such a way as to obtain a better competitor. The case of four branches is rather simple to rule out. Since there cannot be three branches outside $\Gamma$ (because otherwise we could pinch together two of them as before) the only possibility is that exactly two of them are contained in $\Gamma$, but that means that the other two can be projected onto $\Gamma$ in such a way as to decrease the energy of the set. Therefore there are no four-branched minimal cones. By the same argument no cone with more than 4 branches can be sliding minimal.

\chapter{Two-dimensional cones}\label{two-dimensional cones}

In this Chapter we will discuss two-dimensional cones in the half-plane $\R^3:=\{(x,y,z)\in\R^3:z\ge0\}$. The domain of the sliding boundary will be the horizontal plane $\Gamma=\{(x,y,z)\in\R^3:z=0\}$.

\section{Minimal cones}

\subsection{Cartesian products}\label{Cartesian products}

\begin{figure}[h]
\centering
\includegraphics[scale=0.3]{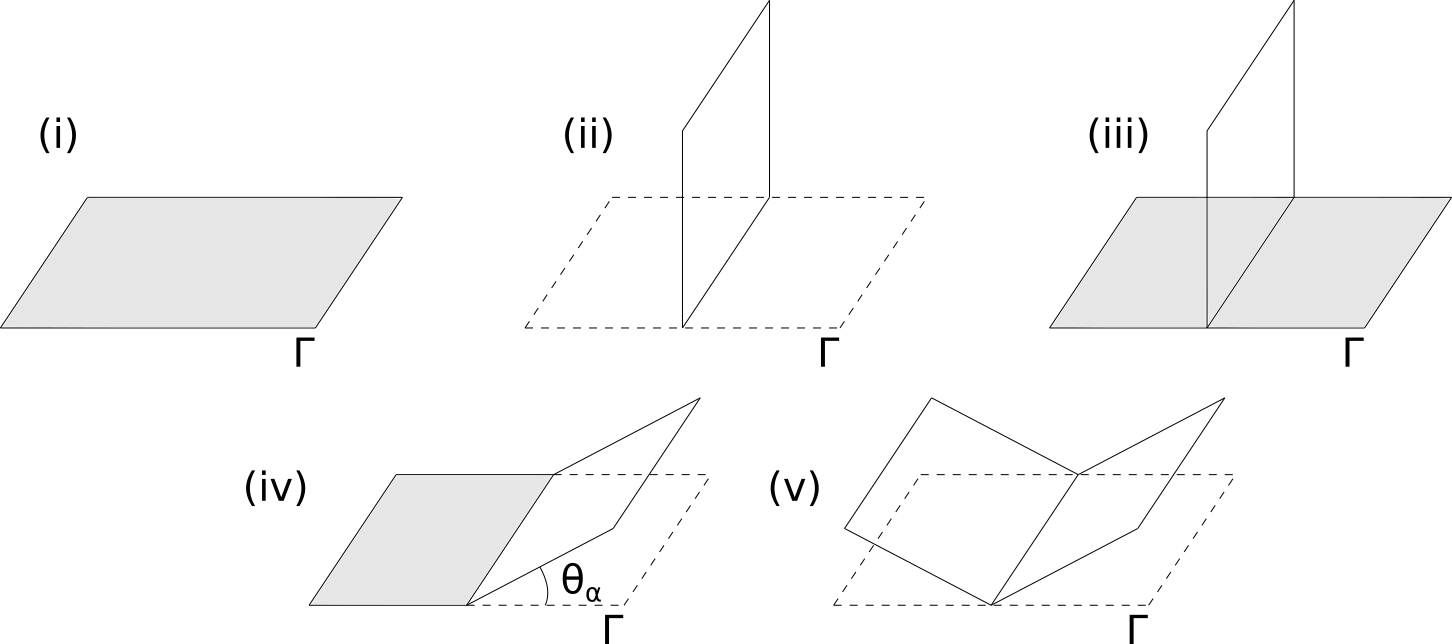}
\caption{Cones obtained as the Cartesian product of $\R$ with a 1-dimensional minimal cone in $\R^2_+$ (the gray region is the intersection between the cones and $\Gamma$).}
\label{coni2in3}
\end{figure}

Let us start our discussion with the 2-dimensional that can be obtained as the Cartesian product of $\R$ with one of the 1-dimensional minimal cones in the previous Chapter (see Figure \ref{coni2in3}).

The minimality of this kind of cones can be proved by a slicing argument. In the following we are going to provide a proof for a cone of type (iv), the minimality of the other cones can be proved in the same way.

Let us denote with $B$ the unit ball of $\R^2$ centred at the origin, and let $A$ be the intersection of the 1-dimensional cone of type (iv) with the ball B. Let $C:=B\times[0,1]\subset \R^2\times\R$ be a cylinder and $D:=A\times[0,1]\subset \R^2\times\R$, then $D$ is the intersection of a 2-dimensional cone of type (iv) with the cylinder $C$. Let us now identify $\R^2\times\R$ with $\R^3$, abusing the notation we will denote with $J_\alpha$ both the functional defined on $\R^2_+$ and the corresponding functional defined on $\R^3_+$. Let $\phi:D\to\R^3$ be a sliding deformation acting in the interior of the cylinder $C$, that is to say $\phi$ is a Lipschitz function and $\phi(W)\subset \textrm{Int}(C)$ where $W:=\{p\in D:\phi(p)\neq p\}$. Therefore $M:=\phi(D)$ is a sliding competitor in the cylinder $C$. Let $\llbracket A\rrbracket\in \mathscr{P}_1(\R^2,\Z_2)$ be the 1-dimensional polyhedral chain with coefficient in $\Z_2$ whose support is $A$, and let $\llbracket D\rrbracket:=\llbracket A\rrbracket\times\llbracket0,1\rrbracket\in \mathscr{P}_2(\R^3,\Z_2)$ be the 2-dimensional polyhedral with coefficient in $\Z_2$ whose support is $D$, it follows that $f_\sharp(\llbracket D\rrbracket)$ is supported in $M$. Let us now define the following orthogonal projection
\begin{equation}
\begin{aligned}
\pi:&\R^2\times\R\to\R\\
&(x,y,t)\mapsto t.
\end{aligned}
\end{equation}
Since $\H^2(M)<+\infty$ we have that $\H^1(M\cap\pi^{-1}(t))<+\infty$ for almost every $t\in[0,1]$, therefore the slice $\langle \phi_\sharp\llbracket D\rrbracket,\pi,t\rangle\in \mathscr{F}_1(\R^2,\Z_2)$ exists for almost every $t\in[0,1]$ since its support is contained in $M\cap\pi^{-1}(t)$.
Moreover the boundary of this slice is the same as the boundary of $\langle \llbracket D\rrbracket,\pi,t\rangle$, and it can be shown using the properties of the slices and the definition of $\phi$ as follows
\begin{equation}
\begin{aligned}
\partial\langle \phi_\sharp\llbracket D\rrbracket,\pi,t\rangle &=\langle \partial\phi_\sharp\llbracket D\rrbracket,\pi,t\rangle\\
&=\langle \phi_\sharp\partial\llbracket D\rrbracket,\pi,t\rangle\\
&=\phi_\sharp\langle\partial\llbracket D\rrbracket,\pi\circ\phi,t\rangle\\
&=\langle\partial\llbracket D\rrbracket,\pi,t\rangle\\
&=\partial\llbracket A\rrbracket.
\end{aligned}
\end{equation}
Hence for almost every $t\in[0,1]$ the slice $M\cap\pi^{-1}(t)$ contains a curve whose length is finite and whose endpoints are the same as the endpoints of $A$, and by the argument in the previous chapter we have that $J_\alpha(M\cap\pi^{-1}(t))\ge J_\alpha(A)$. Therefore we can now compute as follows
\begin{equation}
\begin{aligned}
J_\alpha(M)&\ge\int_0^1J_\alpha(M\cap\pi^{-1}(t))d\H^1(t)\ge\int_0^1J_\alpha(A)d\H^1(t)=J_\alpha(D).
\end{aligned}
\end{equation}

\subsection{Characterisation of minimal cones}

Let $S:=\{(x,y,z)\in\R^3:x^2+y^2+z^2=1\}$ be the unit sphere of $\R^3$ centred at the origin, and let $C\subset\R^3_+$ be a sliding minimal cone. We remark that, since a cone is invariant by dilations, it is completely characterised by its intersection with $S$. In our case $C$ is contained in a half-space, hence it is completely characterised by the graph obtained as its intersection with the upper hemisphere $S$, that is to say $S_+:=S\cap\R^3_+$.

In \cite{taylor1973regularity} and \cite{taylor1976structure} Taylor proved that the graph obtained as the intersection of a minimal cone with the unit sphere must consist of arcs of great circles intersecting three at a time at a finite number of points, and the angles of intersection must be $120^\circ$. Moreover she proved that if $A$ is one of the region in which the sphere is divided by the cone, then $A$ is a spherical polygon having at most 5 sides and the lengths of the arcs of these nets can be computed (in terms of the angle at the origin they subtend) using the following formulae:
\begin{itemize}
\item if $A$ is bounded by only one edge then $A$ is a hemisphere bounded by a great circle;
\item if $A$ is bounded by 2 edges then it is a gore whose side length is $\pi$;
\item If $A$ is a spherical equiangular triangle (all angles $120^\circ$), its side
length  is $\arccos(-\frac{1}{3})$;
\item If $A$ is a spherical quadrilateral with $120^\circ$ angles at its vertices, then
it is ``rectangular'' in the sense that opposite sides are of equal length, and
the lengths $\alpha$ and $\beta$ of its adjacent sides are related by the formula
\begin{equation}\label{lati_rettangolo}
\cos(\beta)=\frac{3-5\cos(\alpha)}{5-3\cos(\alpha)}
\end{equation}
which in terms of half angles becomes
\begin{equation}\label{lati_rettangolo/2}
\cos(\beta/2)=2\sin(\alpha/2)\sqrt{1+3\sin^2(\alpha/2)}
\end{equation}
\item If $A$ is a spherical pentagon ($120^\circ$ angles), and if $\alpha$ and $\beta$ are the lengths of adjacent sides, then the length $\gamma$ of the side adjacent to
neither is given by
\begin{equation}\label{lati_pentagono}
2\cos(\gamma)=\frac{1}{3}+\cos(\alpha)+\cos(\beta)+\cos(\alpha)\cos(\beta)-\sin(\alpha)\sin(\beta).
\end{equation}
\end{itemize}
Finally, applying the slicing argument of the previous subsection to the blow-up of a cone in a point on $\Gamma$, we have that the arcs can meet the equator only with one of the 1-dimensional optimal profiles of Chapter \ref{One-dimensional}. In the following we will use this characterisation of cones trying to classify them in terms of the number of triple junctions and in terms of the spherical polygons they produce.

\begin{figure}[h]
\centering
\includegraphics[scale=0.5]{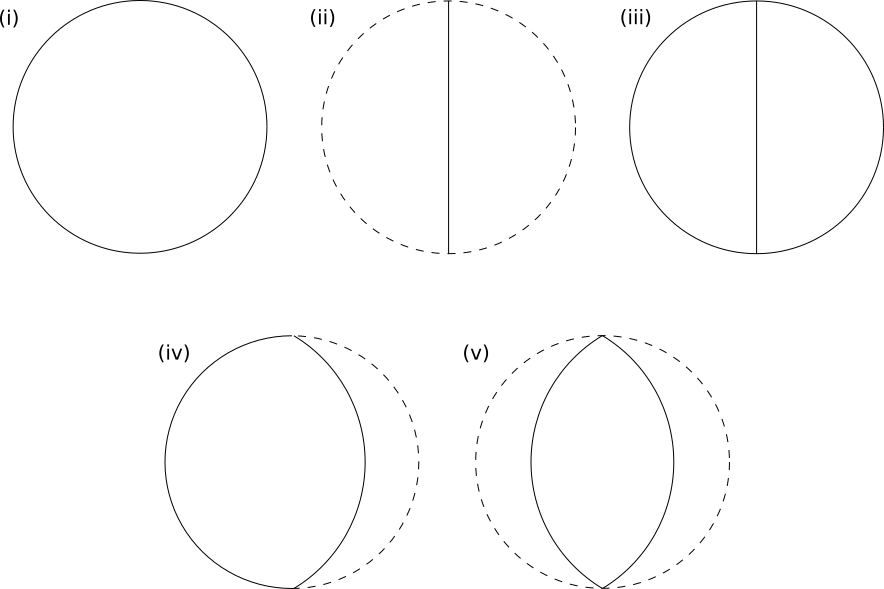}
\caption{Representation of the graphs on the hemisphere characterising the cones in Figure \ref{coni2in3}. The half sphere is seen from above.}
\label{calibrazioneT}
\end{figure}

\subsection{Half $\T$}\label{Half_T}

The first minimal cone we are going to discuss is called $\T_+$ (or half $\T$) and is obtained by taking a cone of type $\T$ as in Figure \ref{YeT.png}, flipping it upside down, placing its barycentre at the origin, and finally intersecting it with the half-space $\R^3_+$ (see Figure \ref{calibrazioneT}).
\begin{figure}[h]
\centering
\includegraphics[scale=0.4]{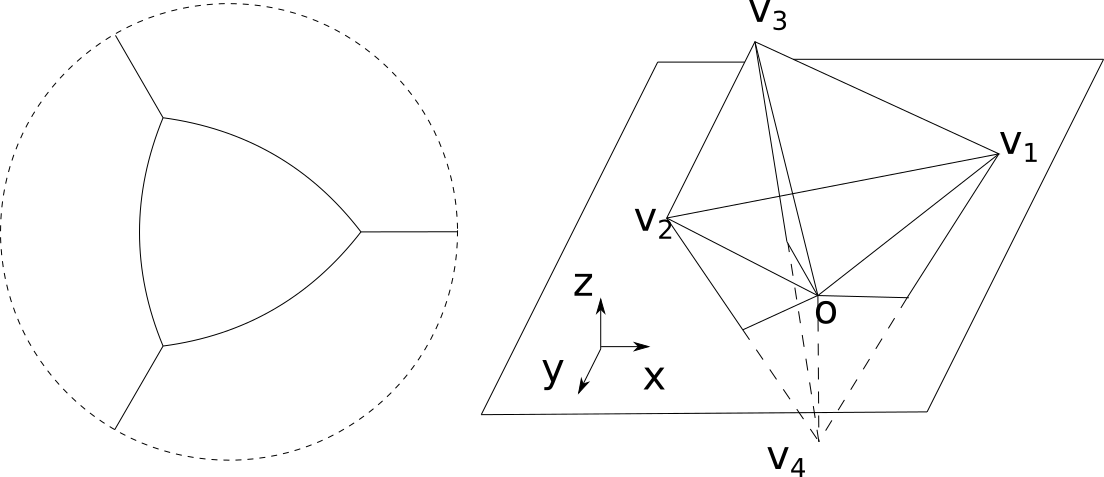}
\caption{On the right the cone $\mathbf{T}_+$ or ``half $\mathbf{T}$'', and on the left its intersection with the hemisphere.}
\label{calibrazioneT}
\end{figure}
 In this section we are going to prove the following theorem:

\begin{thm}\label{teoremaT+}
The cone $\T_+$ is an $\alpha$-sliding minimiser in the half-space $\R^3_+$ with respect to $\Gamma=\partial\R^3_+$ if and only if $\alpha\ge\sqrt{\frac{2}{3}}$.
\end{thm}

Let us begin with a formal description of $\T_+$, which will be similar to the description of the regular simplex $\Delta^n$ given in Section \ref{Minimal_cones}. Given the following unitary vectors
\begin{equation}
\begin{array}{lrrrl}
v_1=\Big( & 2\frac{\sqrt{2}}{3}, &0, &\frac{1}{3} &\Big)\\
v_2=\Big( & -\frac{\sqrt{2}}{3}, &\sqrt{\frac{2}{3}}, &\frac{1}{3} &\Big)\\
v_3=\Big( & -\frac{\sqrt{2}}{3}, &-\sqrt{\frac{2}{3}}, &\frac{1}{3} &\Big)\\
v_4=\Big( & 0, &0, &-1 &\Big)\\
\end{array}
\end{equation}
let $\Delta^3:=[v_1,v_2,v_3,v_4]$ be the 3-dimensional simplex whose vertices are $v_1,v_2,v_3,v_4$, and let $\mathbf{\Delta}^3:=cone(\sk_{1}(\Delta^3))$; then $\T_+:=\mathbf{\Delta}^3\cap\R^3_+$.

We set $\Delta^3_+:=\Delta^3\cap\R^3_+$; for $1=1,2,3,4$ let $F_i:=\Delta^3_i$  be the two-dimensional face of $\Delta^3$ opposed to the vertex $v_i$  and $F_i^+:=F_i\cap\R^3_+$. Let $M\subset\R^3_+$ be a sliding competitor for $\T_+$ such that the symmetric difference between the two is contained in $\Delta^3_+\cap\textrm{Int}(\Delta^3)$. It follows that $\R^3\setminus M$ has 2 unbounded connected components: one of them contains $F_4^+$, and the other one contains the other three faces $F_1^+$, $F_2^+$, $F_3^+$, as well as the lower half-space $\R^3\setminus\R^3_+$. However $\R^3_+\setminus M$ has 4 unbounded connected components each one of them containing one of the faces $F_i^+$.

For $i=1,2,3,4$ we name $V_i$ the connected component of $\R^3_+\setminus M$ containing $F_i^+$ and we set $V_0:=\R^3\setminus\R^3_+$. In case $\R^3_+\setminus M$ also has some bounded connected components we just include them in $V_1$. By the definition of sliding competitor we have that $M$ is a Lipschitz image of $\T_+$, therefore $M$ has locally finite $2$-dimensional Hausdorff measure. This means that any of the $V_i$ is a set whose perimeter is locally finite. In particular the sets $U_i:=V_i\cap\Delta^3$ have finite perimeter, these are the set we will use in order to apply the divergence theorem. Let us now introduce the following notation:
\begin{eqnarray}
M_{ij} &:=& \partial^* U_i\cap\partial^* U_j\\
M_i &:=& \bigcup_{j=0,\, j\neq i}^4M_{ij}\\
M_i^+ &:=& \bigcup_{j=1}^4M_{ij}, \textrm{ for } i=1,2,3,4\\
\widetilde{M} &:=& M_1^+\cup M_2^+\cup M_3^+\cup M_4
\end{eqnarray}
where $\partial^*E$ denotes the reduced boundary of $E$ and $i,j=0,1,2,3,4$ unless otherwise specified. As we have seen in Section \ref{Minimal_cones} it follows that $\widetilde{M}\subset M\cap\Delta^3_+$ and $\H^2$-almost every point in $\widetilde{M}$ lies on the interface between exactly two regions of its complement (taking into account also $U_0$). Moreover the interfaces between different couples of regions are essentially disjoint with respect to $\H^2$.

 Let us now remark the following useful facts
\begin{eqnarray}
M_i &=& M_i^+\cup M_{i0}\\
\widetilde{M}\setminus\Gamma &=& \bigcup_{i,j\neq0}M_{ij}\quad \textrm{ (up to $\H^2$-negligible sets)}\\
\widetilde{M}\cap\Gamma &=& M_{40}\quad \textrm{ (up to $\H^2$-negligible sets)}\\
\partial^*U_i &=& F_i^+\cup M_i=F_i^+\cup M_i^+\cup M_{i0}\\
\partial^*U_0 &=& \partial^*(\Delta^3\setminus\R^3_+)=(\partial^*\Delta^3\setminus\R^3_+)\cup M_0.
\end{eqnarray}
Finally we denote with $n_i$ the exterior unit normal to $\partial U_i$, and $n_{ij}$ will denote the unit normal to $M_{ij}$ pointing in direction of $U_j$. The vectors of the calibration we will use are $w_i:=\frac{-v_i}{|v_i-v_j|}=-\sqrt{\frac{3}{8}}v_i$ for $i=1,2,3,4$; whose components are the following:
\begin{equation}\label{calibrazione_mezzoT}
\begin{array}{lrrrl}
w_1=\Big( & -\frac{1}{\sqrt{3}}, &0, &-\frac{1}{2\sqrt{6}} &\Big)\\
w_2=\Big( & \frac{1}{2\sqrt{3}}, &-\frac{1}{2}, &-\frac{1}{2\sqrt{6}} &\Big)\\
w_3=\Big( & \frac{1}{2\sqrt{3}}, &\frac{1}{2}, &-\frac{1}{2\sqrt{6}} &\Big)\\
w_4=\Big( & 0, &0, &\frac{1}{2}\sqrt{\frac{3}{2}} &\Big).\\
\end{array}
\end{equation}

We are now ready to start the paired calibration machinery. After applying the divergence theorem to the sets $U_i$ with the vectors $w_i$ we can isolate the interface with the negative half-space as follows
\begin{equation}\label{calcoli_calibrazione_mezzoT}
\begin{aligned}
\sqrt{\frac{3}{8}}\H^2\left(\cup_iF^+_i\right)
&=\sum_{i=1}^4\int_{F_i^+}w_i\cdot n_id\mathcal{H}^2\\
&=-\sum_{i=1}^4\int_{M_i}w_i\cdot n_id\mathcal{H}^2\\
&=-\sum_{i=1}^4\int_{M_i^+}w_i\cdot n_id\mathcal{H}^2-\sum_{i=1}^4\int_{M_{i0}}w_i\cdot n_{i0}d\mathcal{H}^2\\
&=\sum_{1\le i<j\le4}\int_{M_{ij}}(w_j-w_i)\cdot n_{ij}d\mathcal{H}^2+\sum_{i=1}^4\int_{M_{i0}}w_i\cdot \hat{z}d\mathcal{H}^2,
\end{aligned}
\end{equation}
where we denoted with $\hat{z}$ the vector $(0,0,1)$. We could replace $n_{i0}$ with $-\hat{z}$ because all the interfaces of kind $M_{i0}$ are contained in the horizontal plane $\Gamma$ and their exterior normal points downward. Let us now focus on the second of the two sums. Using \eqref{calibrazione_mezzoT} and the definition of $M_{ij}$ we get
\begin{equation}
\begin{aligned}\label{calcoli_interfaccia}
\sum_{i=1}^4\int_{M_{i0}}w_i\cdot \hat{z}d\mathcal{H}^2 &=-\sum_{i=1}^3\frac{1}{2\sqrt{6}}\mathcal{H}^2(M_{i0})+\sqrt{\frac{3}{8}}\mathcal{H}^2(M_{40})\\
&=-\frac{1}{2\sqrt{6}}\mathcal{H}^2(M_{10}\cup M_{20}\cup M_{30})+\sqrt{\frac{3}{8}}\mathcal{H}^2(M_{40}).
\end{aligned}
\end{equation}
Plugging \eqref{calcoli_interfaccia} in \eqref{calcoli_calibrazione_mezzoT} and using the fact that $\mathcal{H}^2(M_{10}\cup M_{20}\cup M_{30})=\mathcal{H}^2(M_0)-\mathcal{H}^2(M_{40})$ we obtain

\begin{equation}
\begin{aligned}
&\sqrt{\frac{3}{8}}\mathcal{H}^2(\cup_iF_i^+)+\frac{1}{2\sqrt{6}}\mathcal{H}^2(M_0)=\\
&=\sum_{1\le i<j \le4}\int_{M_{ij}}(w_j-w_i)\cdot n_{ij}d\mathcal{H}^2+\sqrt{\frac{2}{3}}\mathcal{H}^2(M_{40})\\
&\le\H^2\left(\widetilde{M}\setminus\Gamma\right)+\sqrt{\frac{2}{3}}\H^2\left(\widetilde{M}\cap\Gamma\right).
\end{aligned}
\end{equation}

Let us now assume $\alpha=\sqrt{\frac{2}{3}}$. Recalling that $\widetilde{M}\subset M\cap\Delta^n$ the previous inequality becomes
\begin{equation}\label{calibrazione_quasi_fatto}
\sqrt{\frac{3}{8}}\mathcal{H}^2(\cup_iF_i^+)+\frac{1}{2\sqrt{6}}\mathcal{H}^2(M_0)\le J_\alpha(\widetilde{M})\le J_\alpha(M).
\end{equation}
Since the left-hand side of \eqref{calibrazione_quasi_fatto} is a constant, and the chain of inequalities turns into a chain of equalities for $M=\mathbf{T}_+$, we proved that this cone is minimal when $\alpha=\sqrt{\frac{2}{3}}$. Moreover is is also minimal for $\alpha'\ge\sqrt{\frac{2}{3}}$ and it is due to the fact that $J_{\alpha'}(\T_+)=J_{\alpha}(\T_+)$ because $\H^2(\T_+\cap\Gamma)=0$. To show the $\alpha'$-minimality of $\T_+$ we can compute as follows:
\begin{equation}
J_{\alpha'}(\T_+)=J_{\alpha}(\T_+)\le J_{\alpha}(M)\le J_{\alpha'}(M)
\end{equation}
for every sliding competitor $M$.

\begin{figure}
\centering
\includegraphics[scale=0.5]{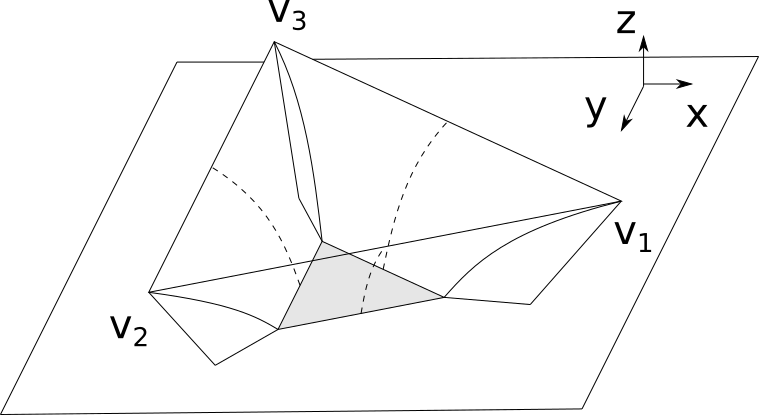}
\caption{Competitor $M$. The dotted line represent the profile given by the function \eqref{profilo_competitore}. The grey region is the intersection between $M$ and $\Gamma$.}
\label{competitoreT}
\end{figure}

Let us now check that $\mathbf{T}_+$ is not a sliding minimiser for $\alpha<\sqrt{\frac{2}{3}}$, we will do it by providing a better competitor, similar to the competitor to the cone over the 1-skeleton of a cube in \ref{conocubo} (see Figure \ref{competitoreT}).

\begin{figure}
\centering
\includegraphics[scale=0.45]{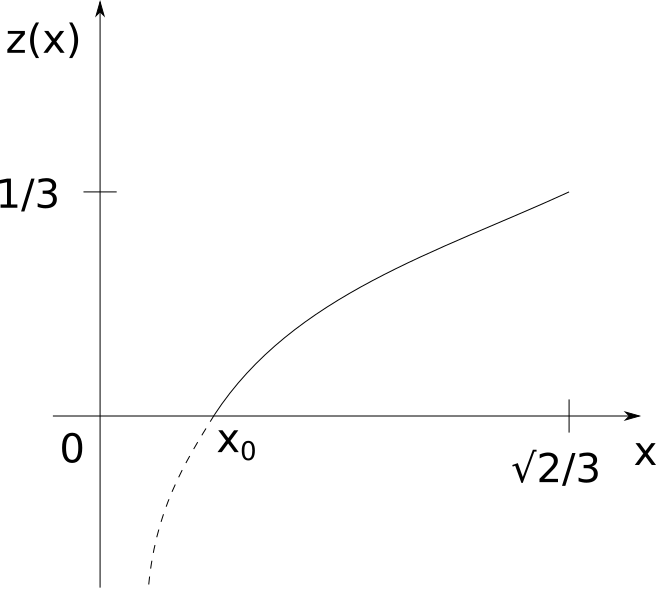}
\caption{Graph of the profile function $z$, the dotted line is its negative part.}
\label{funzionez}
\end{figure}

This is a modification of $\mathbf{T}_+$ obtained by pushing it down on $\Gamma$ in such a way to create a little horizontal equilateral triangle centred at the origin and bending the sloping folds along a profile defined by the positive part of the following function (see Figure \ref{funzionez})

\begin{equation}\label{profilo_competitore}
z(x)=\frac{x}{\sqrt{2}}+c\log\left(\frac{3}{\sqrt{2}}x\right),\quad z'(x)=\frac{1}{\sqrt{2}}+\frac{c}{x}.
\end{equation}
The function has been chosen in such a way that $z(\sqrt{2}/3)=1/3$ for any $c>0$ (which will be fixed later). Moreover we set $x_0\in(0,\sqrt{2}/3)$ as the unique solution of $z(x_0)=0$
\begin{equation}\label{eqx0}
0=\frac{x_0}{\sqrt{2}}+c\log\left(\frac{3}{\sqrt{2}}x_0\right).
\end{equation}
What we just defined is a one-parameter family of admissible competitors for $\T_+$. It is possible to use either $c$ or $x_0$ as parameter of the family, and in particular we can use \eqref{eqx0} in order to write $c$ in terms of $x_0$ as follows
\begin{equation}\label{cx0}
c=-\frac{x_0}{\sqrt{2}\log\left(\frac{3}{\sqrt{2}}x_0\right)}.
\end{equation}

Let us now show that for every $\alpha<\sqrt{\frac{2}{3}}$ there exist a competitor in this family, denoted by $M_c$, with less energy than the cone. As usual we define $M_c^+:=M_c\cap\Gamma$, and we have

\begin{equation}\label{JalphaMc}
J_\alpha(M_c)=\H^2(M_c^+)+\alpha\H^2(M\cap \Gamma).
\end{equation}
By construction the part of $M_c$ laying on $\Gamma$ is just an equilateral triangle whose apothem is $x_0$ therefore its area is $3\sqrt{3}x_0^2$. On the other hand $M_c^+$ is composed by three equal vertical folds and three equal curved folds. Let us call $B$ and $V$ respectively any of the bended or vertical folds, hence
\begin{equation}
\H^2(M_c^+)=3\H^2(B)+3\H^2(V),
\end{equation}
and \eqref{JalphaMc} becomes
\begin{equation}\label{JalphaMc2}
J_\alpha(M_c)=3\H^2(B)+3\H^2(V)+\alpha\H^2(M\cap \Gamma).
\end{equation}
Since the profile of both $B$ and $V$ can be described in terms of the function $z(x)$, their area can be computed by slicing along the direction of the $x$ axis and then integrating on the interval $[x_0,\sqrt{2}/3]$. To be precise we should integrate on the interval $[-x_0,-\sqrt{2}/3]$, what we are actually doing is to compute the area of $M_c^+$ after a symmetry with respect to the $yz$ plane. We have that
\begin{equation}
\H^2(B)=\int_{x_0}^{\frac{\sqrt{2}}{3}}2\sqrt{3}x\sqrt{1+(z'(x))^2}dx,
\end{equation}
where $2\sqrt{3}x$ is the length of the slice and $\sqrt{1+(z'(x))^2}$ is the Jacobian of the function $Z(x):=(x,z(x))$; and
\begin{equation}
\H^2(V)=\int_{x_0}^{\frac{\sqrt{2}}{3}}2z(x)dx
\end{equation}
where we have to multiply by two because the fold $V$ that we are slicing makes and angle $\pi/3$ with the direction of the $x$ axis.
Using the inequality
\begin{equation}
\sqrt{\frac{3}{2}x^2+\sqrt{2}cx+c^2}\le\sqrt{\frac{3}{2}}x+\frac{c}{\sqrt{3}}+\sqrt{\frac{2}{3}}\frac{c^2}{x}
\end{equation}
we can now compute as follows
\begin{equation}\label{faldaB}
\begin{aligned}
\H^2(B) &=2\sqrt{3}\int_{x_0}^{\frac{\sqrt{2}}{3}}x\sqrt{1+\frac{1}{2}+\sqrt{2}\frac{c}{x}+\frac{c^2}{x^2}}dx\\
&=2\sqrt{3}\int_{x_0}^{\frac{\sqrt{2}}{3}}\sqrt{\frac{3}{2}x^2+\sqrt{2}cx+c^2}dx\\
&\le2\sqrt{3}\int_{x_0}^{\frac{\sqrt{2}}{3}}\left[\sqrt{\frac{3}{2}}x+\frac{c}{\sqrt{3}}+\sqrt{\frac{2}{3}}\frac{c^2}{x}\right]dx\\
&=3\sqrt{3}\left[\frac{1}{2}\sqrt{\frac{3}{2}}x^2+\frac{c}{\sqrt{3}}x+\sqrt{\frac{2}{3}}c^2\log\left(\frac{3}{\sqrt{2}}x\right)\right]_{x_0}^\frac{\sqrt{2}}{3}
\end{aligned}
\end{equation}
and
\begin{equation}\label{faldaV}
\begin{aligned}
\H^2(V) &=2\int_{x_0}^{\frac{\sqrt{2}}{3}}\left[\frac{x}{\sqrt{2}}+c\log\left(\frac{3}{\sqrt{2}}x\right)\right]dx\\
&=2\left[\frac{x^2}{2\sqrt{2}}+cx\log\left(\frac{3}{\sqrt{2}}x\right)-cx\right]_{x_0}^\frac{\sqrt{2}}{3}.
\end{aligned}
\end{equation}
Therefore, using \eqref{JalphaMc2}, \eqref{faldaB} and \eqref{faldaV} we obtain
\begin{equation}
J_\alpha(M_c)\le6\left[\sqrt{2}x^2+\sqrt{2}c^2\log\left(\frac{3}{\sqrt{2}}x\right)+cx\log\left(\frac{3}{\sqrt{2}}x\right)\right]_{x_0}^\frac{\sqrt{2}}{3}+\alpha3\sqrt{3}x_0^2.
\end{equation}
Let us now compute the energy of the cone $\mathbf{T}_+$. Since $\mathbf{T}_+=M_0$ we can just use the previous computation with $c=0$ and $x_0=0$ and we get
\begin{equation}
\begin{aligned}
z(x) &=\frac{x}{\sqrt{2}},\quad z'(x)=\frac{1}{\sqrt{2}},\\
J_\alpha(\mathbf{T}_+) &=3\int_{0}^{\frac{\sqrt{2}}{3}}2\sqrt{3}x\sqrt{1+(z'(x))^2}dx+3\int_{0}^{\frac{\sqrt{2}}{3}}2z(x)dx\\
&=12\sqrt{2}\int_{0}^{\frac{\sqrt{2}}{3}}xdx=\frac{4}{3}\sqrt{2}.
\end{aligned}
\end{equation}
Now we can compare the energy of the competitor with the energy of the cone and, using  \eqref{cx0}, we have
\begin{equation}
\begin{aligned}
&J_\alpha(M_c)-J_\alpha(M_0) =\\
=&6\left[\sqrt{2}\frac{2}{9}-\!\sqrt{2}x_0^2-\!\sqrt{2}c^2\log\left(\frac{3}{\sqrt{2}}x_0\right)-\!cx_0\log\left(\frac{3}{\sqrt{2}}x_0\right)\right]+\!\alpha3\sqrt{3}x_0^2-\!\frac{4}{3}\sqrt{2}\\
=&-6\sqrt{2}x_0^2-6\sqrt{2}c^2\log\left(\frac{3}{\sqrt{2}}x_0\right)-6cx_0\log\left(\frac{3}{\sqrt{2}}x_0\right)+\alpha3\sqrt{3}x_0^2\\
=&-6\sqrt{2}x_0^2-3\sqrt{2}\frac{x_0^2}{\log\left(\frac{3}{\sqrt{2}}x_0\right)}+3\sqrt{2}x^2_0+\alpha3\sqrt{3}x_0^2\\
=&3x_0^2\left[-\sqrt{2}-\frac{\sqrt{2}}{\log\left(\frac{3}{\sqrt{2}}x_0\right)}+\alpha\sqrt{3}\right].
\end{aligned}
\end{equation}
Therefore the competitor has less energy than the cone if
\begin{equation}
\alpha\le\sqrt{\frac{2}{3}}\left[1+\frac{1}{\log\left(\frac{3}{\sqrt{2}}x_0\right)}\right].
\end{equation}
Since
\begin{equation}
\lim_{x_0\to0^+}\left[1+\frac{1}{\log\left(\frac{3}{\sqrt{2}}x_0\right)}\right]=1^-
\end{equation}
it follows that for every $\alpha<\sqrt{\frac{2}{3}}$ there exists an $x_0$ (and hence a $c$) such that $J_\alpha(M_c)\le J_\alpha(\mathbf{T}_+)$. This completes the proof of \ref{teoremaT+}.

\begin{figure}
\centering
\includegraphics[scale=0.5]{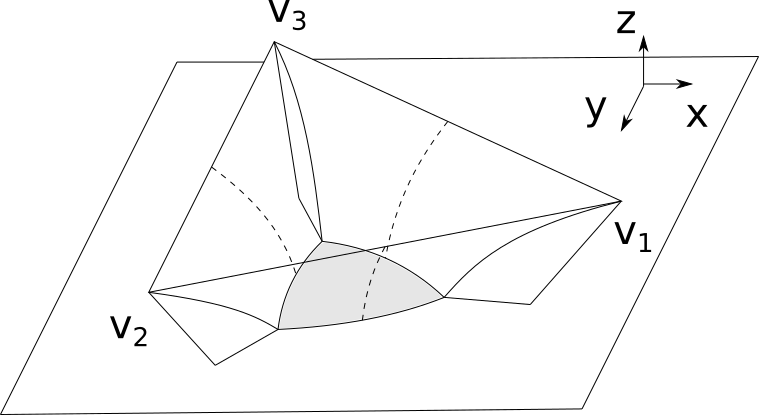}
\caption{Minimiser for $\alpha<\alpha_3$. The grey region (fat triangle) is the intersection between the set and $\Gamma$.}
\label{minimoT}
\end{figure}

For a reason that will be clear later let us call $\alpha_3:=\sqrt{\frac{2}{3}}$. In order to understand why $\alpha_3$ is the threshold between minimality and non-minimality of the cone let us recover how it shows up in the calibration argument as the difference of two scalar products
\begin{equation}\label{valore_soglia}
\alpha_3=\sqrt{\frac{2}{3}}=\frac{1}{2}\sqrt{\frac{3}{2}}+\frac{1}{2\sqrt{6}}=w_4\cdot\hat{z}-w_i\cdot\hat{z}.
\end{equation}
The last term in \eqref{valore_soglia} can be reformulated as follows
\begin{equation}
w_4\cdot\hat{z}-w_i\cdot\hat{z}=(w_4-w_i)\cdot\hat{z}=n_{i4}\cdot\hat{z},
\end{equation}
where $n_{i4}$ is the unit normal to one of the sloping folds of $\T_+$. Hence $\alpha_3$ turns out to be the cosine of the angle between the two unit vectors $n_{i4}$ and $\hat{z}$, which is the same as the cosine of the angle between the plane containing the interface $M_{i4}$ and $\Gamma$ since the previous vectors are the unit normals to these planes. Therefore when $\alpha=\alpha_3$ the sloping folds of $\T_+$ satisfy the optimal profile condition $\cos\theta_{\alpha_3}=\alpha_3$, stated in Chapter \ref{One-dimensional}. In case $\alpha<\alpha_3$ the corresponding optimal profile angle $\theta_\alpha$ is bigger than $\theta_{\alpha_3}$ and a competitor, in order to minimise its energy, would try to attain such optimal angle with the sloping folds; resulting in a shape similar to $M_c$. Numerical simulations with Brakke's Surface Evolver show that in this case the minimiser is very similar to $M_c$, it particular the part of it laying on $\Gamma$ is a ``fat'' triangle (see Figure \ref{minimoT}).

On the other hand if $\alpha>\alpha_3$ the optimal profile angle $\theta_\alpha$ is smaller than $\theta_{\alpha_3}$ and for a competitor is impossible to attain it with its sloping folds minimising the energy at the same time.

\subsection{$\Y_\beta$}

\begin{figure}
\centering
\includegraphics[scale=0.35]{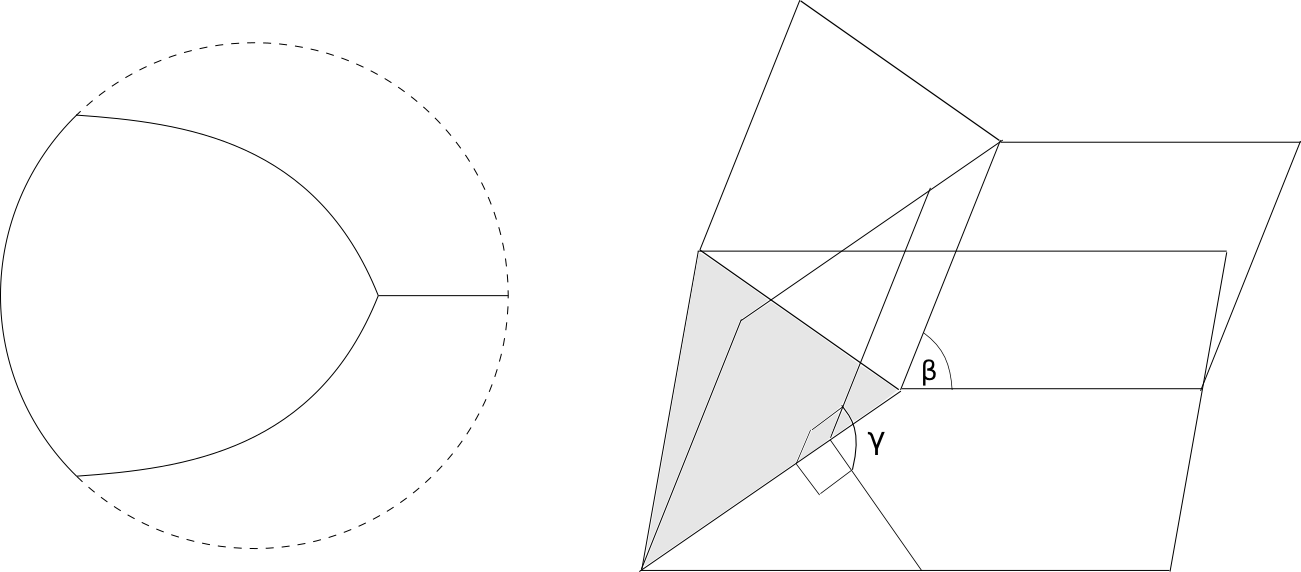}
\caption{The cone $\Y_\beta$ (the grey region is the intersection between the cone and $\Gamma$) and on the left its intersection with the hemisphere.}
\label{Y_tagliato.png}
\end{figure}

Let us introduce a new kind of cone that we will call $\Y_\beta$, where $\beta\in[0,\pi/2]$. It can be obtained with the following procedure: first take the cone $\Y\subset\mathbb{R}^3$, tilt it in a proper way, then intersect it with $\R^3_+$, and finally join it with a section of the horizontal plane $\Gamma$ (see Figure \ref{Y_tagliato.png}).

\begin{figure}
\centering
\includegraphics[scale=0.35]{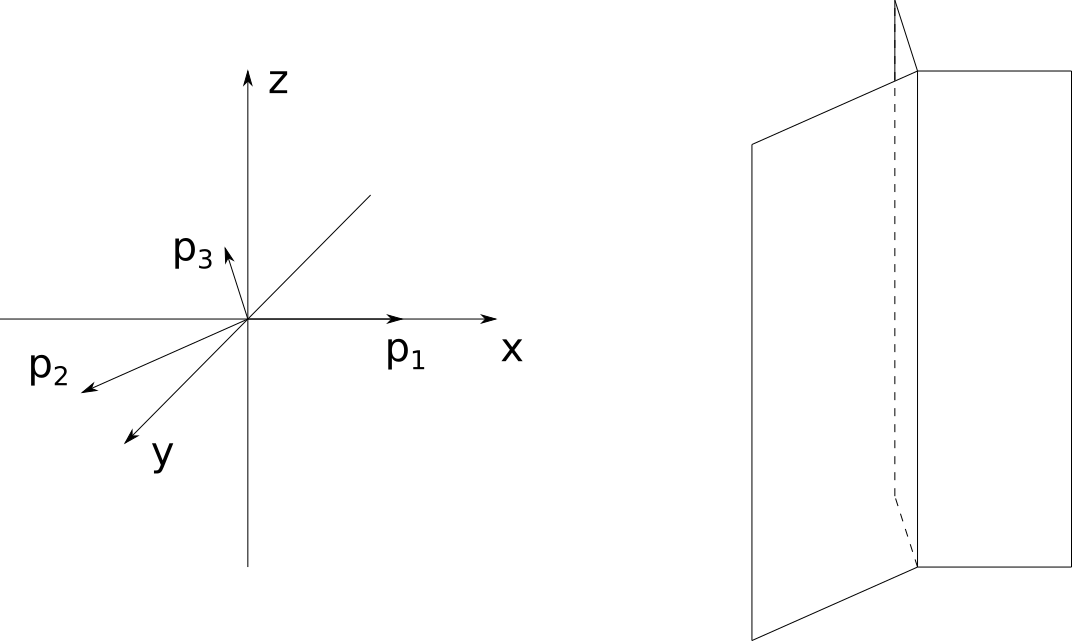}
\caption{The three points generating the cone $Y$ on the left, and the cone $\Y$ on the right.}
\label{yeY.png}
\end{figure}

The cone $\Y_\beta$ can be constructed as follows. Let $r$ be the straight line spanned by the vector $\hat{z}=(0,0,1)$ and $Y\subset\R^3$ be the 1-dimensional cone over the three points
\begin{equation}
\begin{aligned}
p_1 &=(1,0,0)\\
p_2 &=\left(-\frac{1}{2},\frac{\sqrt{3}}{2},0\right)\\
p_3 &=\left(-\frac{1}{2},-\frac{\sqrt{3}}{2},0\right).
\end{aligned}
\end{equation}
We define the 2-dimensional cone $\Y$ as the Cartesian product $Y\times r$ (see Figure \ref{yeY.png}).
Now we can rotate $\Y$ around the $y$ axis using the following rotation
\begin{equation}
R_\beta:=\left(\begin{array}{ccc}
\sin\beta &0 &\cos\beta\\
0 &1 &0\\
-\cos\beta &0 &\sin\beta
\end{array}\right)
\end{equation}
and we obtain a cone, $R_\beta(\Y)$, composed of one vertical fold and two sloping ones meeting $\Gamma$ with the same angle $\gamma$ (by symmetry). Let us now introduce the vectors
\begin{equation}
\begin{aligned}
n_2 &:=\left(-\frac{\sqrt{3}}{2},-\frac{1}{2},0\right)\\
n_3 &:=\left(-\frac{\sqrt{3}}{2},\frac{1}{2},0\right),
\end{aligned}
\end{equation}
which are respectively orthogonal to $p_2$ and $p_3$ and are contained in $\Gamma$. It follows that the normal vectors to the sloping folds of $R_\beta(\Y)$ can be obtained by rotating $n_2$ and $n_3$ with $R_\beta$; that is to say:
\begin{equation}
\begin{aligned}
m_2:=R_\beta(n_2)=&\left(-\frac{\sqrt{3}}{2}\sin\beta,-\frac{1}{2},\frac{\sqrt{3}}{2}\cos\beta\right)\\
m_3:=R_\beta(n_3)=&\left(-\frac{\sqrt{3}}{2}\sin\beta,\frac{1}{2},\frac{\sqrt{3}}{2}\cos\beta\right).
\end{aligned}
\end{equation}
Since $\cos\gamma=(m_2,\hat{z})$ we find the following relation between $\beta$ and $\gamma$:
\begin{equation}\label{condizione_Yb}
\frac{\sqrt{3}}{2}\cos\beta=\cos\gamma.
\end{equation}
The intersection between $R_\beta(\Y)$ and $\Gamma$ is the union of three half-lines meeting at the origin, each one being the intersection of one of the folds with $\Gamma$. We name $q_1$, $q_2$ and $q_3$ these half-lines, and using $m_2$ and $m_3$ we find that
\begin{equation}\label{retteq123}
\begin{aligned}
q_1 &=\{(t,0,0): t\ge0\}\\
q_2 &=\left\{\left(t,-\sqrt{3}\sin\beta \;t,0\right): t\le0\right\}\\
q_3 &=\left\{\left(t,\sqrt{3}\sin\beta \;t,0\right): t\le0\right\}.
\end{aligned}
\end{equation}
Let us call $S$ the convex subset of of $\Gamma$ bounded by $q_2\cup q_3$; we can now define our cone as
\begin{equation}
\Y_\beta:=(R_\beta(\Y)\cap\R^3_+)\cup S.
\end{equation}

In the the rest of this section we will prove the following theorem.
\begin{thm}\label{teoremaYbeta}
The cone $\Y_\beta$ is sliding minimal if and only if
\begin{equation}\label{condizioneYb}
\alpha=\frac{\sqrt{3}}{2}\cos\beta.
\end{equation}
\end{thm}
In order for $\Y_\beta$ to be a minimal set a necessary condition is that for every $P\in\Y_\beta$ the blow-up of $\Y_\beta$ at $P$ has to be a minimal cone (here we are assuming $P\neq0$ because otherwise the necessary condition would turn into a tautology since the blow-up of $\Y_\beta$ at the origin is $\Y_\beta$ itself). In particular if $P\in q_1\cup q_2\cup q_3$ the blow-up of $\Y_\beta$ has to assume one of the optimal profiles described in Section \ref{Half-plane}. In case $P\in q_1$ this condition is satisfied because the cone assumes the profile (ii). When $P\in q_2\cup q_3$ the profile taken by the cone is of type (iv), and in this case the blow-up is a minimal cone if and only if
\begin{equation}
\cos\gamma=\alpha.
\end{equation}
Hence, combining the equality \eqref{condizione_Yb} with the previous one we obtain that the condition \eqref{condizioneYb} expressed in Theorem \ref{teoremaYbeta} is necessary for the minimality of $\Y_\beta$.

We are now going to prove with a calibration argument that condition \eqref{condizioneYb} is also sufficient for the minimality of $\Y_\beta$. The calibration we will use is obtained by rotating with $R_\beta$ a calibration for the cone $Y\subset\Gamma$. Let
\begin{equation}
\begin{aligned}
v_1 &:=\left(-\frac{1}{\sqrt{3}},0,0\right)\\
v_2 &:=\left(\frac{1}{2\sqrt{3}},-\frac{1}{2},0\right)\\
v_3 &:=\left(\frac{1}{2\sqrt{3}},\frac{1}{2},0\right)\\
\end{aligned}
\end{equation}
be a calibration for $Y$ in $\Gamma$ (hence for $\Y$ in $\R^3$), we define the calibration for $\Y_\beta$ as $w_i:=R_\beta(v_i)$ and we get
\begin{equation}
\begin{aligned}
w_1 &=\left(-\frac{\sin\beta}{\sqrt{3}},0,\frac{\cos\beta}{\sqrt{3}}\right)\\
w_2 &=\left(\frac{\sin\beta}{2\sqrt{3}},-\frac{1}{2},-\frac{\cos\beta}{2\sqrt{3}}\right)\\
w_3 &=\left(\frac{\sin\beta}{2\sqrt{3}},\frac{1}{2},-\frac{\cos\beta}{2\sqrt{3}}\right).
\end{aligned}
\end{equation}

\begin{figure}
\centering
\includegraphics[scale=0.35]{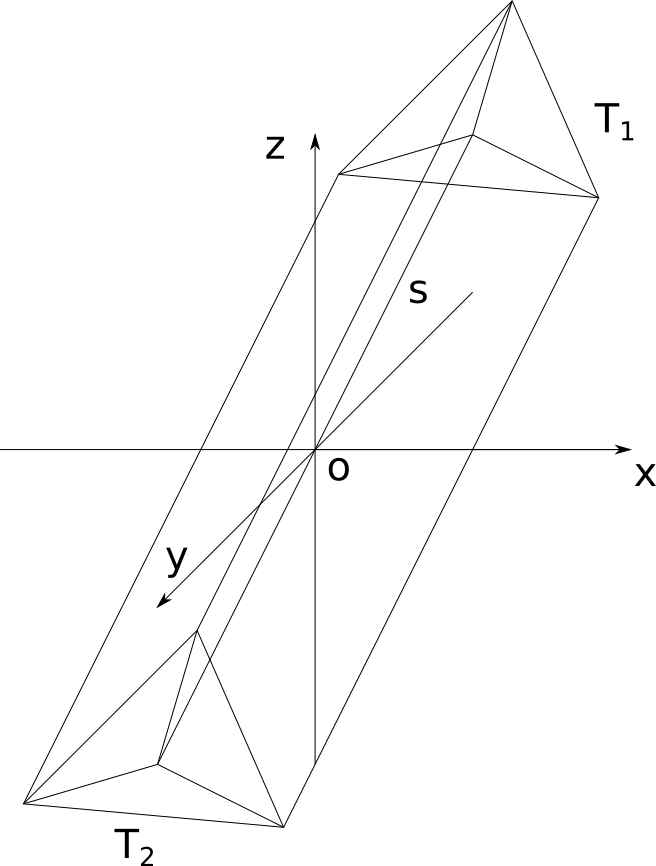}
\caption{The right prism $P$ enclosing part of the cone $R_\beta(\Y)$.}
\label{prisma.png}
\end{figure}

Let $s$ be the straight line spanned by $R_\beta(\hat{z})=(\cos\beta,0,\sin\beta)$, in the following we will refer to it as the spine of $\Y_\beta$. We fix the compact set in which the sliding deformation takes place as the right prism $P$ whose bases are two equilateral triangles, $T_1$ and $T_2$, orthogonal to the spine  $s$ and centred on it, such that their vertices lie in $\Y_\beta$. We assume that the barycentre of $P$ is the origin and its height is large enough such that the two triangular bases do not intersect $\Gamma$ (see Figure \ref{prisma.png}).

By definition the vectors $w_i$ are orthogonal to the lateral faces of the prism $P$, therefore, for $i=1,2,3$ we name $F_i$ the face orthogonal to $w_i$. We set $P_+:=P\cap\R^3_+$ and $F_i^+:=F_i\cap\R^3_+$. Let $M\subset\R^3_+$ be a sliding competitor for $\Y_\beta$ such that $M\triangle\Y_\beta\subset P_+$. It follows that $\R^3_+\setminus M$ has 3 unbounded connected components. For $i=1,2,3$ we name $V_i$ the connected component of $\R^3_+\setminus M$ containing $F_i^+$, and we set $V_0:=\R^3\setminus\R^3_+$. In case $\R^3_+\setminus M$ also has bounded connected components we can just include them in $V_1$. For $i=0,1,2,3$ the sets $V_i$ have locally finite perimeter, hence the sets $U_i:=V_i\cap P$ are finite perimeter sets. We can now introduce the following notation
\begin{eqnarray}
M_{ij} &:=& \partial^* U_i\cap\partial^* U_j\\
M_i &:=& \bigcup_{j=0, \, j\neq i}^3M_{ij}\\
M_i^+ &:=& \bigcup_{j=1}^3M_{ij}, \textrm{ for } 1=1,2,3\\
\widetilde{M} &:=& M_1^+\cup M_2^+\cup M_3.
\end{eqnarray}
As in the previous section it follows that $\widetilde{M}\subset M\cap P$, $\H^2$-almost every point in $\widetilde{M}$ lies on the interface between exactly two of the $U_i$, and the sets $M_{ij}$ are essentially disjoint when $i<j$. Finally we call $n_i$ the outer normal to $\partial U_i$, and $n_{ij}$ the unit normal to $M_{ij}$ pointing in direction of $U_j$.

We are now ready for the calibration argument. In particular, when applying the divergence theorem to the sets $U_i$ with respect to the vectors $w_i$, we can ignore the contribution given by the upper base of the prism since by definition the vectors $w_i$ are orthogonal to the normal vector to $T_1$. Thus we can compute as follows:

\begin{equation}\label{calcoli_calibrazioneYb}
\begin{aligned}
\frac{1}{\sqrt{3}}\sum_{i=1}^3\H^2(F_i^+) =& \sum_{i=1}^3\int_{F_i^+}w_i\cdot n_id\mathcal{H}^2\\
=&-\sum_{i=1}^3\int_{M_i}w_i\cdot n_id\mathcal{H}^2\\
=&-\sum_{i=1}^3\int_{M_i^+}w_i\cdot n_id\mathcal{H}^2-\sum_{i=1}^3\int_{M_{i0}}w_i\cdot n_id\mathcal{H}^2\\
=&\sum_{1\le i<j\le3}\int_{M_{ij}}(w_j-w_i)\cdot n_{ij}d\mathcal{H}^2+\sum_{i=1}^3\int_{M_{i0}}w_i\cdot \hat{z}d\mathcal{H}^2.
\end{aligned}
\end{equation}
In the last line we isolated the contribution given by the interface with the negative half-space, let us now focus on this term. Using the definition of the vectors $w_i$ we get

\begin{equation}\label{calcoli_interfacciaYb}
\sum_{i=1}^3\int_{M_{i0}}w_i\cdot \hat{z}d\mathcal{H}^2 =\frac{\cos\beta}{\sqrt{3}}\mathcal{H}^2(M_{10})-\frac{\cos\beta}{2\sqrt{3}}\left(\mathcal{H}^2(M_{20})+\mathcal{H}^2(M_{30})\right).
\end{equation}
Plugging \eqref{calcoli_interfacciaYb} in \eqref{calcoli_calibrazioneYb} and using the fact that $\mathcal{H}^2(M_{20})+\mathcal{H}^2(M_{30})=\mathcal{H}^2(M_0)-\mathcal{H}^2(M_{10})$ we obtain
\begin{equation}
\begin{aligned}
\frac{1}{\sqrt{3}}\mathcal{H}^2(\cup_iF_i^+)+\frac{\cos\beta}{2\sqrt{3}}\mathcal{H}^2(M_0) &=\sum_{1\le i<j \le3}\int_{M_{ij}}(w_j-w_i)\cdot n_{ij}d\mathcal{H}^2\\
&+\frac{\sqrt{3}}{2}\cos\beta\mathcal{H}^2(M_{10})\\
&\le\H^2\left(\widetilde{M}\setminus\Gamma\right)+\frac{\sqrt{3}}{2}\cos\beta\left(\widetilde{M}\cap\Gamma\right).
\end{aligned}
\end{equation}
Assuming $\alpha=\frac{\sqrt{3}}{2}\cos\beta$ we have
\begin{equation}\label{calibrazione_quasi_fattoYb}
C(\alpha,P)\le J_\alpha(\widetilde{M})\le J_\alpha(M)
\end{equation}
where $C(\alpha,K)$ is a constant that depends only on the parameter $\alpha$ and on the compact set $K$ containing $M\triangle \Y_\beta$, in our case $K=P$.
Since the left-hand side of \eqref{calibrazione_quasi_fattoYb} is a constant, and the chain of inequalities turns into a chain of equalities for $M=\Y_\beta$, we proved that this cone is minimal when $\alpha=\frac{\sqrt{3}}{2}\cos\beta$. Therefore condition \eqref{condizioneYb} is both necessary and sufficient for the cone $\Y_\beta$ to be minimal.

I the same way it has happened in the previous Section, it might look surprising how a necessary condition for minimality turns into a sufficient one. However, as we did before, let us remark how the constant $\frac{\sqrt{3}}{2}\cos\beta$ showed up from the computation. It appears as the scalar product between the normal vector to the sloping folds and the normal vector to the domain of the sliding boundary. Since the two vectors have unitary norm their scalar product simply is the cosine of the angle between them, which is the same as the cosine of the angle between the planes they are orthogonal to,
\begin{equation}
\frac{\sqrt{3}}{2}\cos\beta=\frac{\cos\beta}{2\sqrt{3}}+\frac{\cos\beta}{\sqrt{3}}=(w_i-w_3,\hat{z})=\cos\gamma.
\end{equation}
And this means that the reason why we impose \eqref{condizioneYb} as necessary for the minimality of the cone, is actually the same reason that makes it sufficient (the optimal profile angle between the sloping folds and $\Gamma$).

 Let us also remark that the cones of type $\Y_\beta$ satisfying that condition form a one parameter family of minimal cones, depending on the angle $\beta\in[0,\pi/2]$ or, equivalently, on the parameter $\alpha\in[0,1]$. In particular when $\alpha=0$ we have $\beta=\pi/2$, therefore $\Y_\beta$ becomes the union of a vertical half $\Y$ with the section of $\Gamma$ contained in between the two half-lines (cfr. \ref{retteq123})
\begin{equation}
\begin{aligned}
q_2 &=\left\{\left(t,-\sqrt{3}\sin\beta \;t,0\right): t\le0\right\}\\
q_3 &=\left\{\left(t,\sqrt{3}\sin\beta \;t,0\right): t\le0\right\}.
\end{aligned}
\end{equation}
However, since in this case the energy functional $J_0$ does not take into account any set laying on $\Gamma$, up to a $J_0$-negligible set $\Y_{\pi/2}$ is the same as a vertical half $\Y$. On the opposite, when $\alpha=1$ we have that $\beta=0$ and $\Y_\beta$ turns into $\mathbf{V}_{\pi/6}:=V_{\pi/6}\times\R$.

\subsection{$\overline{\Y}_\beta$}

\begin{figure}
\centering
\includegraphics[scale=0.4]{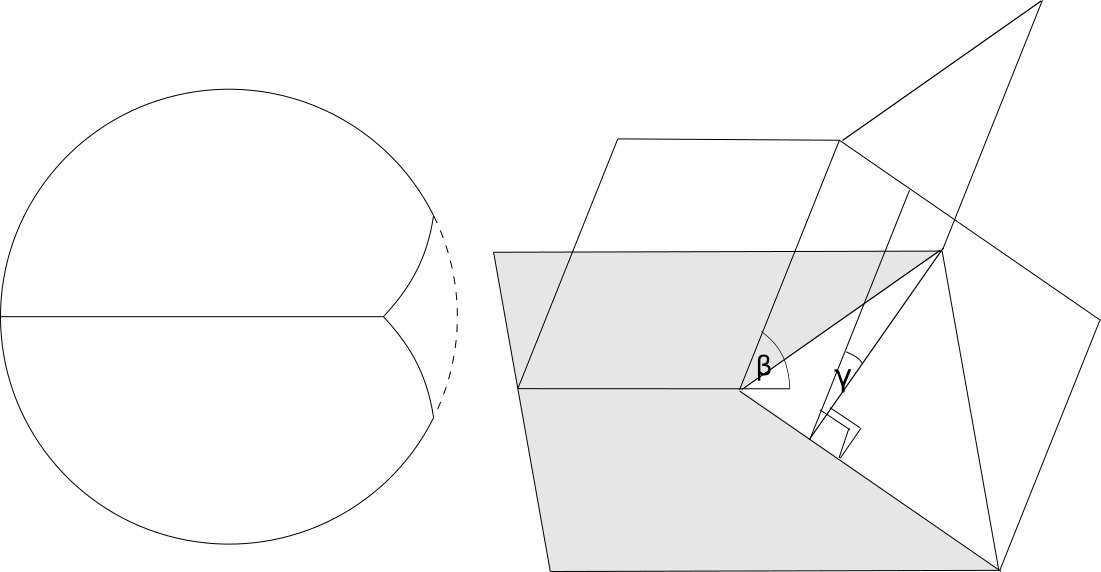}
\caption{The cone $\overline{\Y}_\beta$ (the grey region is the intersection between the cone and $\Gamma$) and on the left its intersection with the hemisphere.}
\label{altro_Y}
\end{figure}

We can use the previous construction to produce another cone that we call $\overline{\Y}_\beta$, for $\beta\in[0,\pi/2]$. It can be obtained with the same procedure as before: first take a cone $\Y\subset\mathbb{R}^3$ symmetric to the previous one, tilt it in a proper way, then intersect it with $\R^3_+$, and finally add to it a section of the horizontal plane $\Gamma$ (see Figure \ref{altro_Y}).

\begin{figure}[h]
\centering
\includegraphics[scale=0.35]{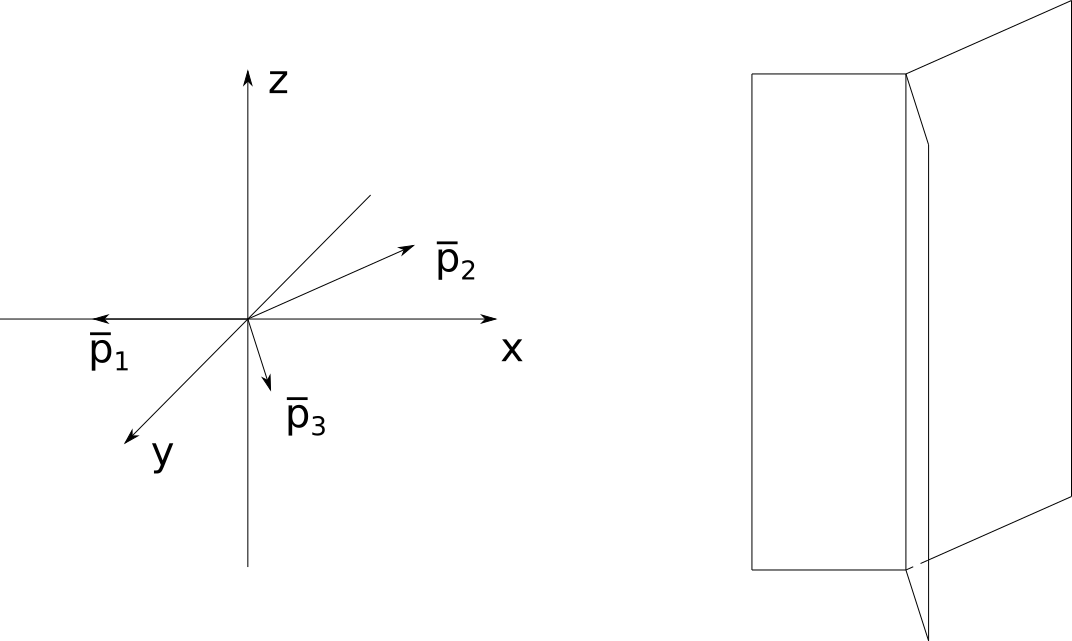}
\caption{The three points generating the cone $\overline{Y}$ on the left, and the cone $\overline{\Y}$ on the right.}
\label{yeY-.png}
\end{figure}

The cone $\overline{\Y}_\beta$ can be constructed as follows. Let $r$ be the straight line spanned by the vector $\hat{z}=(0,0,1)$ and $\overline{Y}\subset\R^3$ be the 1-dimensional cone over the three points
\begin{equation}
\begin{aligned}
\overline{p}_1 &=(-1,0,0)\\
\overline{p}_2 &=\left(\frac{1}{2},-\frac{\sqrt{3}}{2},0\right)\\
\overline{p}_3 &=\left(\frac{1}{2},\frac{\sqrt{3}}{2},0\right),
\end{aligned}
\end{equation}
in particular $\overline{p}_i=-p_i$ where the $p_i$ are the points defined in the previous section, and it follows that $\overline{Y}=-Y$ (where given $E\subset\R^3$ we denote $-E:=\{x\in\R^3:-x\in E\}$). Then we define the 2-dimensional cone $\overline{\Y}$ as the Cartesian product $\overline{Y}\times r$ (see Figure \ref{yeY-.png}), and again we have $\overline{\Y}=-\Y$.
Let us now rotate $\overline{\Y}$ around the $y$ axis with the rotation $R_\beta$. We obtain again a cone, $R_\beta(\overline{\Y})$, composed by one vertical fold and two sloping ones meeting $\Gamma$ with the same angle $\gamma$. The normal vectors to $\overline{p}_2$ and $\overline{p}_3$, respectively $\overline{n}_2$ and $\overline{n}_3$, after an appropriate rotation provide the normal vectors to the sloping folds of $R_\beta(\overline{\Y})$, respectively $\overline{m}_2$ and $\overline{m}_3$. That is
\begin{equation}
\begin{aligned}
\overline{n}_2 &:=\left(-\frac{\sqrt{3}}{2},-\frac{1}{2},0\right)\\
\overline{n}_3 &:=\left(-\frac{\sqrt{3}}{2},\frac{1}{2},0\right)\\
\overline{m}_2 &:=R_\beta(\overline{n}_2)=\left(-\frac{\sqrt{3}}{2}\sin\beta,-\frac{1}{2},\frac{\sqrt{3}}{2}\cos\beta\right)\\
\overline{m}_3 &:=R_\beta(\overline{n}_3)=\left(-\frac{\sqrt{3}}{2}\sin\beta,\frac{1}{2},\frac{\sqrt{3}}{2}\cos\beta\right).
\end{aligned}
\end{equation}
As before, using the fact that $\cos\gamma=(m_2,\hat{z})$, we find the following relation between $\beta$ and $\gamma$:
\begin{equation}\label{condizione_Yb-}
\frac{\sqrt{3}}{2}\cos\beta=\cos\gamma.
\end{equation}
The intersection between $R_\beta(\Y)$ and $\Gamma$ is the union of three half-lines meeting at the origin, each one being the intersection of one of the three folds with $\Gamma$. We name $\overline{q}_1$, $\overline{q}_2$ and $\overline{q}_3$ these half-lines, and using $\overline{m}_2$ and $\overline{m}_3$ we find that
\begin{equation}
\begin{aligned}
\overline{q}_1 &=\{(t,0,0): t\le0\}\\
\overline{q}_2 &=\left\{\left(t,-\sqrt{3}\sin\beta \;t,0\right): t\ge0\right\}\\
\overline{q}_3 &=\left\{\left(t,\sqrt{3}\sin\beta \;t,0\right): t\ge0\right\}.
\end{aligned}
\end{equation}
Let us call $\overline{S}$ the non convex subset of $\Gamma$ bounded by $\overline{q}_2\cup\overline{q}_3$; we can now define our cone as
\begin{equation}
\overline{\Y}_\beta:=(R_\beta(\overline{\Y})\cap\R^3_+)\cup \overline{S}.
\end{equation}

In the the rest of this section we will prove the following theorem
\begin{thm}\label{teoremaYbeta-}
The cone $\overline{\Y}_\beta$ is sliding minimal if and only if
\begin{equation}\label{condizioneYb-}
\alpha=\frac{\sqrt{3}}{2}\cos\beta.
\end{equation}
\end{thm}
First of all we have to check that the blow-up of $\overline{\Y}_\beta$ at any of its point $P$ (except for the origin) is a minimal cone. In case $P\in\overline{q}_1$ this condition is satisfied because the cone assumes the profile (iii). When $P\in \overline{q}_2\cup \overline{q}_3$ the profile taken by the cone is of type (iv), and in this case the blow-up is a minimal cone if and only if
\begin{equation}
\cos\gamma=\alpha,
\end{equation}
and we obtain that condition \eqref{condizioneYb-} is necessary for the minimality of $\overline{\Y}_\beta$.

Let us now provide a calibration for $\overline{\Y}_\beta$, this will show that the condition \eqref{condizioneYb-} is also sufficient for the minimality of $\overline{\Y}_\beta$. We will obtain a calibration for $\overline{\Y}_\beta$ by rotating with $R_\beta$ a calibration for the cone $\overline{\Y}$. Let
\begin{equation}
\begin{aligned}
v_1 &:=\left(\frac{1}{\sqrt{3}},0,0\right)\\
v_2 &:=\left(-\frac{1}{2\sqrt{3}},\frac{1}{2},0\right)\\
v_3 &:=\left(-\frac{1}{2\sqrt{3}},-\frac{1}{2},0\right)\\
\end{aligned}
\end{equation}
be a calibration for $\overline{\Y}$ in $\R^3$, we define the calibration for $\overline{\Y}_\beta$ as $\overline{w}_i:=R_\beta(v_i)$ and we get
\begin{equation}
\begin{aligned}
\overline{w}_1 &=\left(\frac{\sin\beta}{\sqrt{3}},0,-\frac{\cos\beta}{\sqrt{3}}\right)\\
\overline{w}_2 &=\left(-\frac{\sin\beta}{2\sqrt{3}},\frac{1}{2},\frac{\cos\beta}{2\sqrt{3}}\right)\\
\overline{w}_3 &=\left(-\frac{\sin\beta}{2\sqrt{3}},-\frac{1}{2},\frac{\cos\beta}{2\sqrt{3}}\right).
\end{aligned}
\end{equation}

\begin{figure}
\centering
\includegraphics[scale=0.35]{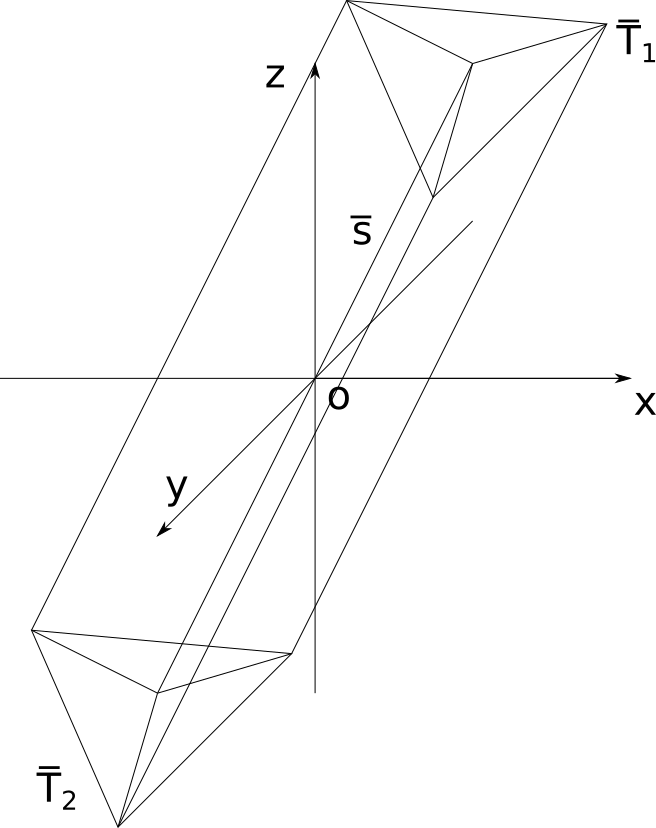}
\caption{The right prism $\overline{P}$ enclosing part of the cone $R_\beta(\overline{\Y})$.}
\label{prisma-.png}
\end{figure}

Let us call $\overline{s}$ the spine of $\overline{\Y}_\beta$, it is spanned by the vector $R_\beta(\hat{z})=(\cos\beta,0,\sin\beta)$. We fix the compact set in which the sliding deformation takes place as the right prism $\overline{P}$ whose bases are two equilateral triangles, $\overline{T}_1$ and $\overline{T}_2$, orthogonal to the spine  $\overline{s}$ and centred on it, such that their vertices lie in $\overline{\Y}_\beta$. We assume that the barycentre of $\overline{P}$ is the origin and its height is large enough such that the two triangular bases do not intersect $\Gamma$ (see Figure \ref{prisma-.png}).

For $i=1,2,3$ we call $F_i$ the lateral face of $\overline{P}$ which is orthogonal to $\overline{w}_i$. We set $\overline{P}_+:=\overline{P}\cap\R^3_+$ and $F_i^+:=F_i\cap\R^3_+$. Let $M\subset\R^3_+$ be a sliding competitor to $\overline{\Y}_\beta$ such that $M\triangle\overline{\Y}_\beta\subset \overline{P}_+$. Than $\R^3_+\setminus M$ has 3 unbounded connected components and we name $V_i$ the one containing $F_i^+$ for $i=1,2,3$, and we set $V_0:=\R^3\setminus\R^3_+$. In case $\R^3_+\setminus M$ also has bounded connected components we include them in $V_1$. For $i=0,1,2,3$ the sets $U_i:=V_i\cap \overline{P}$ are finite perimeter sets, and we introduce the following notation:
\begin{eqnarray}
M_{ij} &:=& \partial^* U_i\cap\partial^* U_j\\
M_i &:=& \bigcup_{j=0}^3M_{ij}\\
M_i^+ &:=& \bigcup_{j=1}^3M_{ij}, \textrm{ for } 1=1,2,3\\
\widetilde{M} &:=& M_1^+\cup M_2\cup M_3.
\end{eqnarray}
The set $\widetilde{M}$ defined above is contained in $M$ and $\H^2$-almost every point in it lies on the interface between exactly two of the regions $U_i$; moreover the sets $M_{ij}$ are essentially disjoint. Let us call $n_i$ the outer normal to $\partial U_i$, and $n_{ij}$ the unit normal to $M_{ij}$ pointing in direction of $U_j$. We can now compute as follows

\begin{equation}\label{calcoli_calibrazioneYb-}
\begin{aligned}
\frac{1}{\sqrt{3}}\sum_{i=1}^3\H^2(F_i^+) =& \sum_{i=1}^3\int_{F_i^+}\overline{w}_i\cdot n_id\mathcal{H}^2\\
=&-\sum_{i=1}^3\int_{M_i}\overline{w}_i\cdot n_id\mathcal{H}^2\\
=&-\sum_{i=1}^3\int_{M_i^+}\overline{w}_i\cdot n_id\mathcal{H}^2-\sum_{i=1}^3\int_{M_{i0}}\overline{w}_i\cdot n_id\mathcal{H}^2\\
=&\sum_{1\le i<j\le3}\int_{M_{ij}}(\overline{w}_j-\overline{w}_i)\cdot n_{ij}d\mathcal{H}^2+\sum_{i=1}^3\int_{M_{i0}}\overline{w}_i\cdot \hat{z}d\mathcal{H}^2.
\end{aligned}
\end{equation}
Let us now consider the second term in the last line

\begin{equation}\label{calcoli_interfacciaYb-}
\sum_{i=1}^3\int_{M_{i0}}\overline{w}_i\cdot \hat{z}d\mathcal{H}^2 =-\frac{\cos\beta}{\sqrt{3}}\mathcal{H}^2(M_{10})+\frac{\cos\beta}{2\sqrt{3}}\left(\mathcal{H}^2(M_{20})+\mathcal{H}^2(M_{30})\right).
\end{equation}
The two previous computation together with the fact that $\mathcal{H}^2(M_{10})=\mathcal{H}^2(M_0)-\mathcal{H}^2(M_{20})-\mathcal{H}^2(M_{30})$ lead to
\begin{equation}
\begin{aligned}
\frac{1}{\sqrt{3}}\mathcal{H}^2(\cup_iF_i^+)+\frac{\cos\beta}{\sqrt{3}}\mathcal{H}^2(M_0) &=\sum_{1\le i<j \le3}\int_{M_{ij}}(\overline{w}_j-\overline{w}_i)\cdot n_{ij}d\mathcal{H}^2\\
&+\frac{\sqrt{3}}{2}\cos\beta\left(\mathcal{H}^2(M_{20})+\mathcal{H}^2(M_{30})\right)\\
&\le\H^2\left(\widetilde{M}\setminus\Gamma\right)+\frac{\sqrt{3}}{2}\cos\beta\left(\widetilde{M}\cap\Gamma\right).
\end{aligned}
\end{equation}
Therefore, in case $\alpha=\frac{\sqrt{3}}{2}\cos\beta$, we get
\begin{equation}\label{calibrazione_quasi_fattoYb-}
C(\alpha,P)\le J_\alpha(\widetilde{M})\le J_\alpha(M)
\end{equation}
where $C(\alpha,K)$ is a constant only depending on the parameter $\alpha$ and on the compact set $K$ containing $M\triangle \overline{\Y}_\beta$, in our case $K=P$.
Since the left-hand side of \eqref{calibrazione_quasi_fattoYb-} is a constant, and the chain of inequalities turns into a chain of equalities for $M=\overline{\Y}_\beta$, we proved that this cone is minimal when $\alpha=\frac{\sqrt{3}}{2}\cos\beta$. Therefore condition \eqref{condizioneYb-} is both necessary and sufficient for the cone $\overline{\Y}_\beta$ to be minimal, and once again the explanation of this fact relies on the optimal angle profile.

The minimal cones of type $\overline{\Y}_\beta$ form a one parameter family depending on the angle $\beta\in[0,\pi/2]$ or, equivalently, on the parameter $\alpha\in[0,1]$. In particular when $\alpha=0$ we have $\beta=\pi/2$, therefore $\overline{\Y}_\beta$ becomes the union of a vertical half $\overline{\Y}$ with the section of $\Gamma$ not contained in between the two half-lines
\begin{equation}
\begin{aligned}
\overline{q}_2 &=\left\{\left(t,-\sqrt{3}\sin\beta \;t,0\right): t\ge0\right\}\\
\overline{q}_3 &=\left\{\left(t,\sqrt{3}\sin\beta \;t,0\right): t\ge0\right\}.
\end{aligned}
\end{equation}
However, since in this case the energy functional $J_0$ does not take into account any set laying on $\Gamma$, up to a $J_0$-negligible set $\overline{\Y}_{\pi/2}$ is the same as a vertical half $\overline{\Y}$. On the opposite, when $\alpha=1$ we have that $\beta=0$ and $\overline{\Y}_\beta$ turns into a cone composed by the union of $\Gamma$ with a vertical half-plane.

\subsection{Double $\Y$}

\begin{figure}
\centering
\includegraphics[scale=0.4]{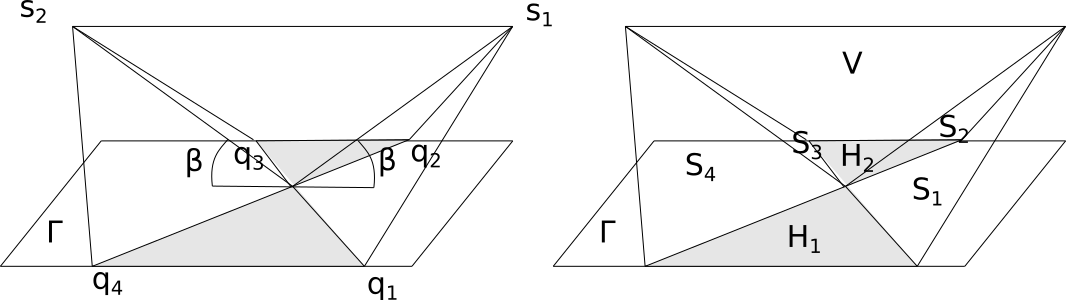}
\caption{The cone $\mathbf{W}_\beta$ with the names of the lines on the left and the names of the folds on the right (the grey region is the intersection between the cone and $\Gamma$).}
\label{doppioY}
\end{figure}

\begin{figure}
\centering
\includegraphics[scale=0.4]{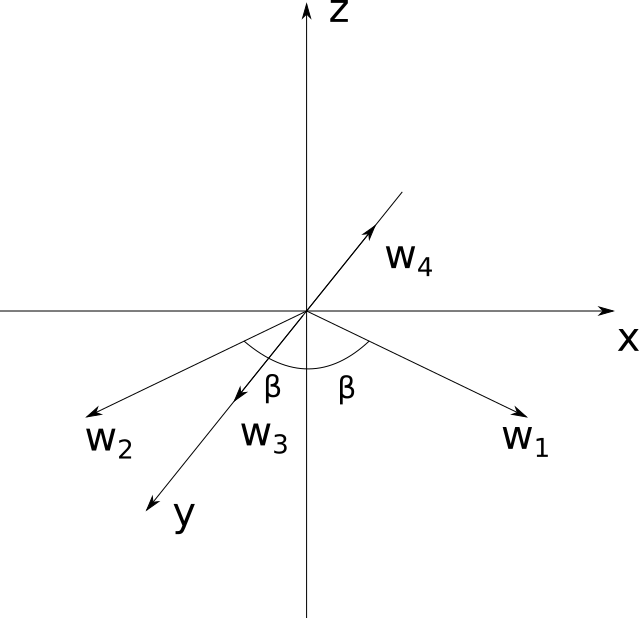}
\caption{Calibration for the cone $\mathbf{W}_\beta$.}
\label{doppioYcal}
\end{figure}

The next cone is called $\mathbf{W}_\beta$, and it is composed by two cones of type $\overline{\Y}_\beta$ symmetric to each other with respect to a vertical plane, and sharing the same vertical fold (see Figure \ref{doppioY}).
It can be constructed as follows. Let $H:=\{(x,y,z)\in\R^3:x\ge0\}$ be the half-space of positive $x$, and $\overline{\Y}_\beta^{x^+}:=\overline{\Y}_\beta\cap H$. Let $R_x$ be the reflection with respect to the $yz$ plane
\begin{equation}
R_x=\left(\begin{array}{ccc}
-1 & 0 & 0\\
0 & 1 & 0\\
0 & 0 & 1\\
\end{array}\right)
\end{equation}
then we can define $W_\beta:=\overline{\Y}_\beta^{x^+}\cup R_x(\overline{\Y}_\beta^{x^+})$ (see Figure \ref{mezzoYb-}).
\begin{figure}
\centering
\includegraphics[scale=0.35]{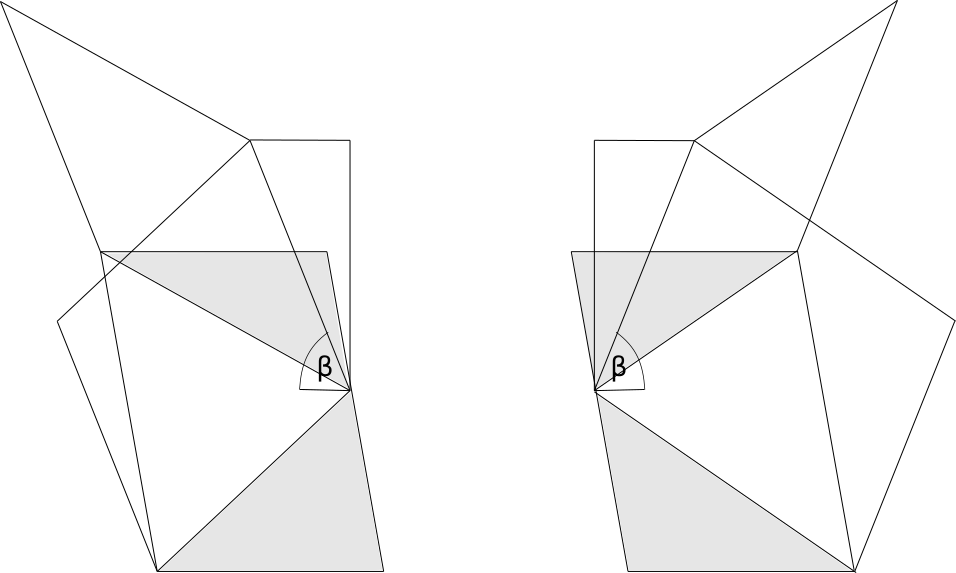}
\caption{The intersection of $\overline{\Y}_\beta$ with $H$ on the right, and its reflection with $R_x$ on the left.}
\label{mezzoYb-}
\end{figure}
\begin{figure}
\centering
\includegraphics[scale=0.35]{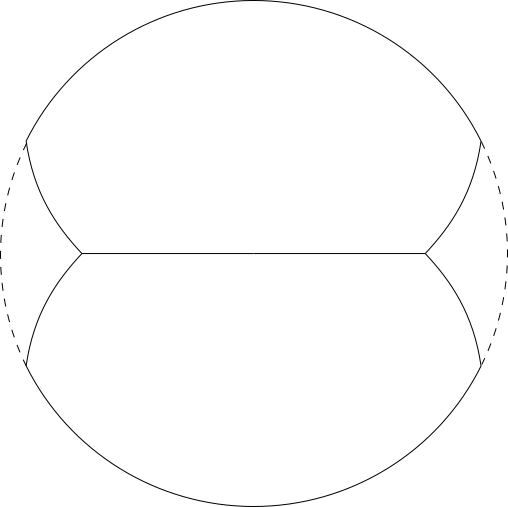}
\caption{The intersection of $\overline{\Y}_\beta$ with the hemisphere.}
\label{mezzoYb-}
\end{figure}
In this section we will prove the following theorem.
\begin{thm}\label{teoremaWbeta}
Let $\sin\beta\le1/\sqrt{3}$, then the cone $\mathbf{W}_\beta$ is sliding minimal if and only if
\begin{equation}\label{condizioneWb}
\alpha=\frac{\sqrt{3}}{2}\cos\beta.
\end{equation}
\end{thm}
As usual the necessity of condition \eqref{condizioneWb} is given by the minimality of the tangent cone of $\mathbf{W}_\beta$ at any of its point and the sufficiency is proved via calibration. By construction the two spines of the sloping $\Y$ cones are the two following half lines
\begin{equation}
\begin{aligned}
s_1 &=\{(\cos\beta\;t,0,\sin\beta\;t):t\ge0\}\\
s_2 &=\{(-\cos\beta\; t,0,\sin\beta\; t):t\ge0\},\\
\end{aligned}
\end{equation}
and the intersection of the four sloping folds with $\Gamma$ are the following four half lines
\begin{equation}
\begin{aligned}
q_1 &=\{(t,\sqrt{3}\sin\beta\;t,0):t\ge0\}\\
q_2 &=\{(t,-\sqrt{3}\sin\beta\;t,0):t\ge0\}\\
q_3 &=\{(t,\sqrt{3}\sin\beta\;t,0):t\le0\}\\
q_4 &=\{(t,-\sqrt{3}\sin\beta\;t,0):t\le0\}.
\end{aligned}
\end{equation}
We can name the folds of $\mathbf{W}_\beta$ as follows (see Figure \ref{doppioY}):
\begin{itemize}
\item[$V$:] the vertical planar face bounded by the two spines $s_1$ and $s_2$;
\item[$H_1$:] the horizontal planar face bounded by $q_1$ and $q_4$;
\item[$H_2$:] the horizontal planar face bounded by $q_2$ and $q_3$;
\item[$S_1$:] the sloping planar face bounded by $q_1$ and $s_1$;
\item[$S_2$:] the sloping planar face bounded by $q_2$ and $s_1$;
\item[$S_3$:] the sloping planar face bounded by $q_3$ and $s_2$;
\item[$S_4$:] the sloping planar face bounded by $q_4$ and $s_2$.
\end{itemize}
For $i=1,2,3,4$ let $\hat{s}_i$ be a unit vector orthogonal to $S_i$. Exploiting the computation of the previous section and the symmetries of $\mathbf{W}_\beta$ we get
\begin{equation}
\begin{aligned}
\hat{s}_1 &=\left(-\frac{\sqrt{3}}{2}\sin\beta, \frac{1}{2}, \frac{\sqrt{3}}{2}\cos\beta\right)\\
\hat{s}_2 &=\left(-\frac{\sqrt{3}}{2}\sin\beta, -\frac{1}{2}, \frac{\sqrt{3}}{2}\cos\beta\right)\\
\hat{s}_3 &=\left(\frac{\sqrt{3}}{2}\sin\beta, -\frac{1}{2}, \frac{\sqrt{3}}{2}\cos\beta\right)\\
\hat{s}_4 &=\left(\frac{\sqrt{3}}{2}\sin\beta, \frac{1}{2}, \frac{\sqrt{3}}{2}\cos\beta\right).\\
\end{aligned}
\end{equation}
Therefore we can choose the following vectors as our calibration (see Figure \ref{doppioYcal})
\begin{equation}
\begin{aligned}
w_1 &=\left(\frac{\sqrt{3}}{2}\sin\beta,0,-\frac{\sqrt{3}}{2}\cos\beta\right)\\
w_2 &=\left(-\frac{\sqrt{3}}{2}\sin\beta,0,-\frac{\sqrt{3}}{2}\cos\beta\right)\\
w_3 &=\left(0,\frac{1}{2},0\right)\\
w_4 &=\left(0,-\frac{1}{2},0\right),
\end{aligned}
\end{equation}
and it is easily seen that, except for $w_1-w_2$, the difference between any two vectors of the calibration is the unit normal to some of the folds of $\mathbf{W}_\beta$. In the following we will explain better the role played by this difference in the calibration argument.
Let us name the 4 connected components of $\R^3_+\setminus\W_\beta$ as follows
\begin{itemize}
\item[$V_1$:] the connected component bounded by $S_1$, $S_2$ and $\Gamma$;
\item[$V_2$:] the connected component bounded by $S_3$, $S_4$ and $\Gamma$;
\item[$V_3$:] the connected component bounded by $H_1$, $S_1$, $V$ and $S_4$;
\item[$V_4$:] the connected component bounded by $H_2$, $S_2$, $V$ and $S_3$.
\end{itemize}
Let $B:=B_1(0)$ be the ball with unitary radius centred at the origin. We choose $B$ as the compact set in which the deformation takes place, and we set $F_i:=V_i\cap\partial B$ for $i=1,2,3,4$ and $B_+:=B\cap\R^3_+$. Let $M$ be a sliding competitor to $\W_\beta$ such that $M\triangle\W_\beta\subset B_+$. It follows that $\R^3_+\setminus M$ has 4 unbounded connected components and, with an abuse of notation, we still call them $V_i$ for $i=1,2,3,4$ (in such a way that these connected components correspond to the previous ones when $M=\W_\beta$). We set $V_0:=\R^3\setminus\R^3_+$ and, in case $\R^3_+\setminus M$ also has bounded connected components we include them in $V_1$. The sets $V_i$ have locally finite perimeter, hence the sets $U_i:=V_i\cap B$ are finite perimeter sets. Let us introduce the following sets
\begin{eqnarray}
M_{ij} &:=& \partial^* U_i\cap\partial^* U_j\\
M_i &:=& \bigcup_{j=0, \, j\neq i}^4M_{ij} \\
M_i^+ &:=& \bigcup_{j=1}^4M_{ij}, \textrm{ for } 1=1,2,3\\
\widetilde{M} &:=& M_1^+\cup M_2^+\cup M_3\cup M_4.
\end{eqnarray}
It follows that $\widetilde{M}\subset M\cap B$, $\H^2$-almost every point in $\widetilde{M}$ lies on the interface between exactly two of the $U_i$, and the sets $M_{ij}$ are essentially disjoint. Finally we call $n_i$ the outer normal to $\partial U_i$, and $n_{ij}$ the unit normal to $M_{ij}$ pointing in direction of $U_j$. Thus
\begin{equation}\label{calcoli_calibrazioneW}
\begin{aligned}
\sum_{i=1}^4\int_{F_i}w_i\cdot n_i d\H^2=&-\sum_{i=1}^4\int_{M_i}w_i\cdot n_id\mathcal{H}^2\\
=&-\sum_{i=1}^4\int_{M_i^+}w_i\cdot n_id\mathcal{H}^2-\sum_{i=1}^4\int_{M_{i0}}w_i\cdot n_id\mathcal{H}^2\\
=&\sum_{1\le i<j\le4}\int_{M_{ij}}(w_j-w_i)\cdot n_{ij}d\mathcal{H}^2+\sum_{i=1}^4\int_{M_{i0}}w_i\cdot \hat{z}d\mathcal{H}^2.
\end{aligned}
\end{equation}
Let us consider the second term in the last line
\begin{equation}\label{calcoli_interfacciaW}
\begin{aligned}
\sum_{i=1}^4\int_{M_{i0}}w_i\cdot \hat{z}d\mathcal{H}^2 &=-\frac{\sqrt{3}}{2}\cos\beta\left(\H^2(M_{10})\H^2+(M_{20})\right)\\
&=\frac{\sqrt{3}}{2}\cos\beta\left(\H^2(M_{30})\H^2+(M_{40})\right)-\frac{\sqrt{3}}{2}\cos\beta\H^2(M_0).
\end{aligned}
\end{equation}
Putting together the two previous computation we get
\begin{equation}\label{calcoliW}
\begin{aligned}
&\sum_{i=1}^4\int_{F_i}w_i\cdot n_i d\H^2+\frac{\sqrt{3}}{2}\cos\beta\H^2(M_0)=\\
&=\sum_{1\le i<j\le4}\int_{M_{ij}}(w_j-w_i)\cdot n_{ij}d\mathcal{H}^2+\frac{\sqrt{3}}{2}\cos\beta\left(\H^2(M_{30})\H^2+(M_{40})\right)\\
&\le\H^2\left(\widetilde{M}\setminus\Gamma\right)+\frac{\sqrt{3}}{2}\cos\beta\H^2\left(\widetilde{M}\cap\Gamma\right).
\end{aligned}
\end{equation}
In order for the last inequality to be true we have to impose
\begin{equation}\label{condizionecalW}
|w_i-w_j|\le1\quad\forall1\le i<j\le4.
\end{equation}
As we remarked above, in order to satisfy \eqref{condizionecalW} we only need to check that $|w_1-w_2|\le1$, and this condition leads to
\begin{equation}
\sin\beta\le\frac{1}{\sqrt{3}}.
\end{equation}
Since the left-hand side in the first line of \eqref{calcoliW} is a constant depending only on the shape of the compact set chosen and on the calibration (which only depend on $\alpha$) in the following we will denote it with $C(\alpha,B)$. Therefore, assuming  $\alpha=\frac{\sqrt{3}}{2}\cos\beta$ and $\sin\beta\le\frac{1}{\sqrt{3}}$ \eqref{calcoliW} becomes
\begin{equation}
C(\alpha,B)\le J_\alpha(\widetilde{M})\le J_\alpha(M).
\end{equation}
Let us now remark that $\H^2(M_{12})=0$ when $M=\W_\beta$, then in this case the previous inequalities turn into a chain of equalities and we proved the minimality of $\W_\beta$.

\begin{figure}
\centering
\includegraphics[scale=0.39]{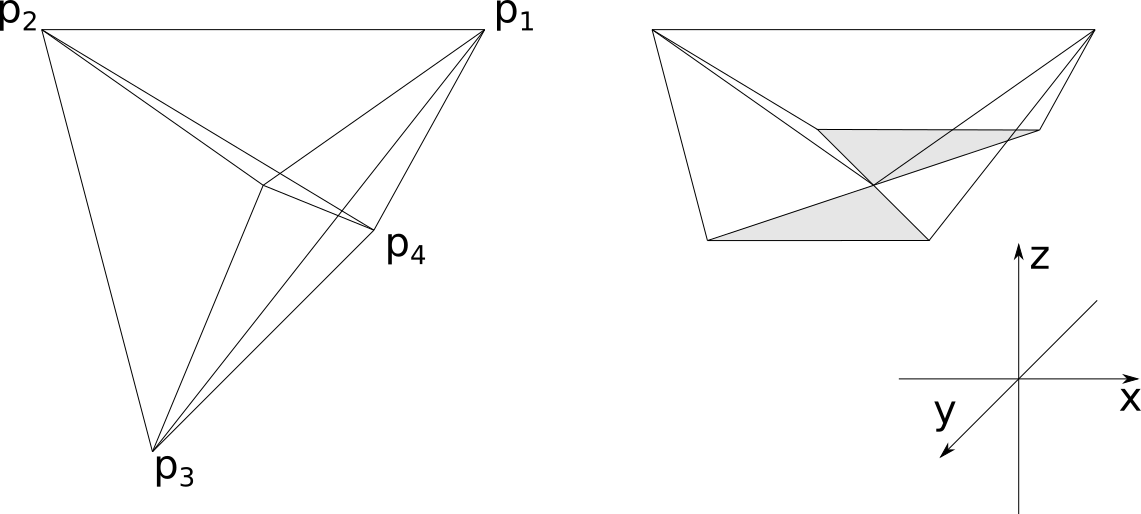}
\caption{On the left the tetrahedron $\Delta_3$ and the cone $\T$ over its skeleton. On the right the corresponding cone $\W_\beta$ and the relative Cartesian coordinate system (not centred at the origin).}
\label{Wtetraedro}
\end{figure}

The cones of type $\W_\beta$ satisfying the minimality condition form a one-parameter family which can be described in therm of the parameter $\alpha\in[1/\sqrt{2},1]$, or equivalently in term of the  angle $\beta \in [0,\arcsin(1/\sqrt{3})]$. In particular, when $\alpha=1$ the cone $\W_\beta$ turns into the union of $\Gamma$ with a vertical half-plane. On the other hand, when $\alpha=1/\sqrt{2}$ the two sloping $\Y$ cones of $\W_\beta$ actually belong to a cone of type $\T$ (see Figure \ref{Wtetraedro}). It can be obtained as the cone over the skeleton of the regular tetrahedron $\Delta_3$ whose vertices are
\begin{equation}
\begin{aligned}
p_1=& \left(\sqrt{\frac{2}{3}},0,\frac{1}{\sqrt{3}}\right) &p_2=& \left(-\sqrt{\frac{2}{3}},0,\frac{1}{\sqrt{3}}\right)\\
p_3=& \left(0,\sqrt{\frac{2}{3}},-\frac{1}{\sqrt{3}}\right) &p_4=& \left(0,-\sqrt{\frac{2}{3}},-\frac{1}{\sqrt{3}}\right).\\
\end{aligned}
\end{equation}

\section{Non-minimal cones}

So far we have seen cones for which the necessary condition for minimality turned out into being also a sufficient condition, at this point a natural question is whether it is always the case. As we are about to see the answer to the question is no. In the following we will give a non exhaustive overview on non minimal cones satisfying the necessary condition for minimality away from the origin. As in the previous sections, for each type of cone we will actually have a one-parameter family of cones, depending on the parameter $\alpha$. Therefore, in order to prove that a cone is not minimal, one has to provide a better competitor for each value of the parameter $\alpha$. In the following, for each family of cones we will show the existence of better competitors for some given values of $\alpha$. These better competitors have been obtained using Brakke's Surface Evolver (see \cite{brakke2013surface}) and, even if they do not provide right away a proof of the non minimality of the whole family, by testing different values of $\alpha$ we get a strong clue in this direction.

\subsection{$\Y+\Y$}
\begin{figure}[h]
\centering
\includegraphics[scale=0.4]{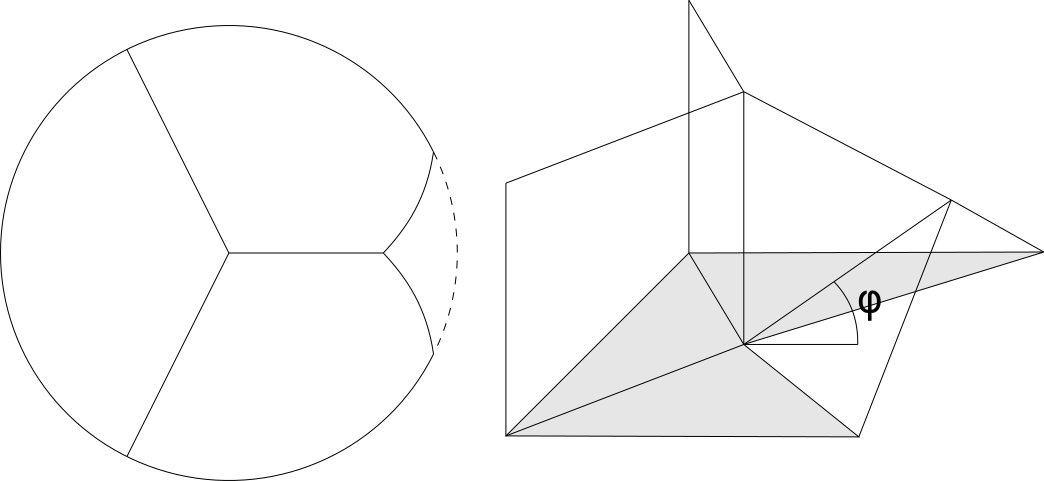}
\caption{On the right the cone $\Y+\Y$ (the grey region is the intersection between the cone and $\Gamma$), and on the left its intersection with the hemisphere.}
\label{Y+Y}
\end{figure}
\begin{figure}[h]
\centering
\includegraphics[scale=0.4]{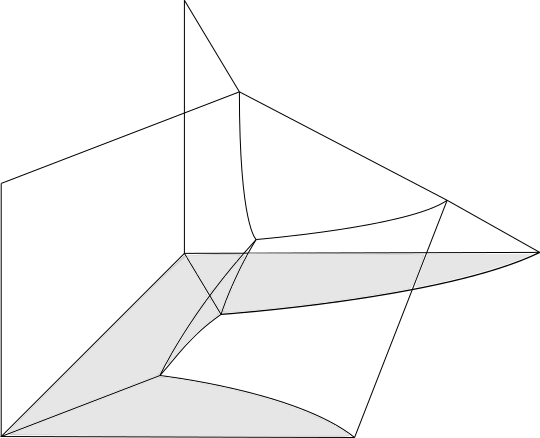}
\caption{A better competitor for the cone $\Y+\Y$ obtained with Brakke's Surface Evolver (the grey region is the intersection between the competitor and $\Gamma$). The curved triangle in the center represent the new interface produced by the pinching.}
\label{Y+Y-competitore}
\end{figure}
The first cone we are going to discuss can be obtained, loosely speaking, by adding a sloping $\Y$ to a vertical $\Y$ and then adding the proper region of $\Gamma$. More precisely the cone $\Y+\Y$ is formed of two cones of type $\Y$ having a common fold (see Figure \ref{Y+Y}). One of them is vertical and the other one is sloping in such a way as to produce an angle $\varphi$ between its spine and the horizontal plane. The cone also contains the part of $\Gamma$ not contained below the two sloping folds, and, in order to satisfy the necessary condition for minimality, $\alpha$ has to be chosen in function of $\varphi$ so that the angle between the sloping folds and the horizontal plane is the optimal one. In the following we will use $\sin\varphi\in[0,1]$ as the parameter of the family. Let us remark that the cone $\Y+\Y$ divides $\R^3_+$ in 4 different connected components, and the only two connected components not having a common interface are the one below the sloping $\Y$ (on the right in Figure \ref{Y+Y}) and the one bounded by $\Gamma$ and two folds of the vertical $\Y$ (on the left in the Figure). A better competitor for $\Y+\Y$ can be obtained by pinching together the folds in such a way as to create an interface in between the two regions that still don't have it (see Figure \ref{Y+Y-competitore}). Surface Evolver provided a better competitor of this type for different values of  $\sin\varphi$ well distributed on the interval $[0,1]$.

\subsection{$\Y+2\Y$}

\begin{figure}
\centering
\includegraphics[scale=0.4]{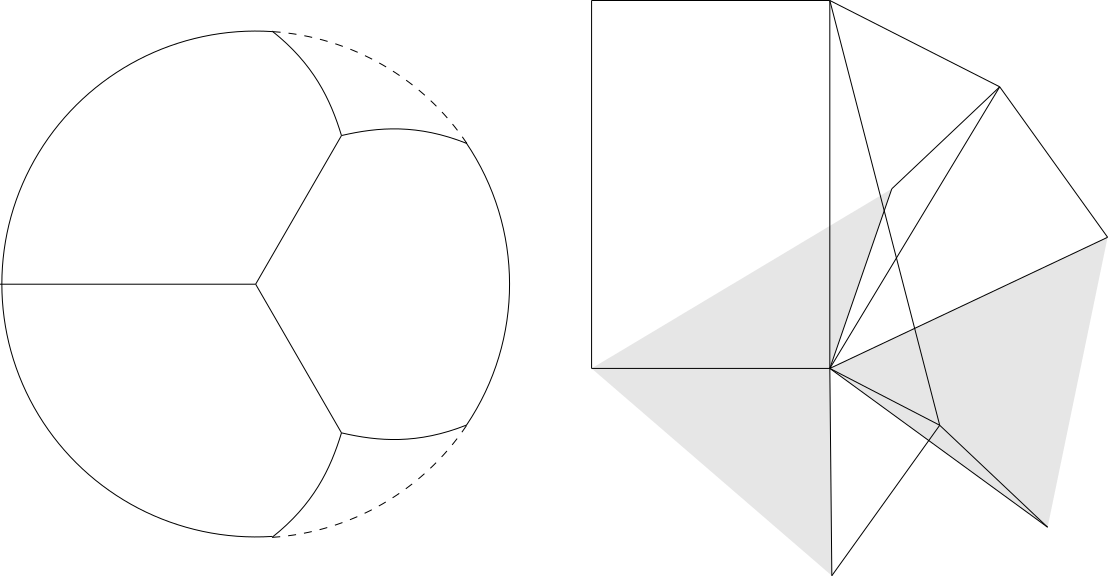}
\caption{On the right the cone $\Y+\Y$ corresponding to a value of $\sin\varphi$ lower than $1/2$ (the grey region is the intersection between the cone and $\Gamma$), and on the left its intersection with the hemisphere.}
\label{Y+2Y'}
\end{figure}

\begin{figure}
\centering
\includegraphics[scale=0.4]{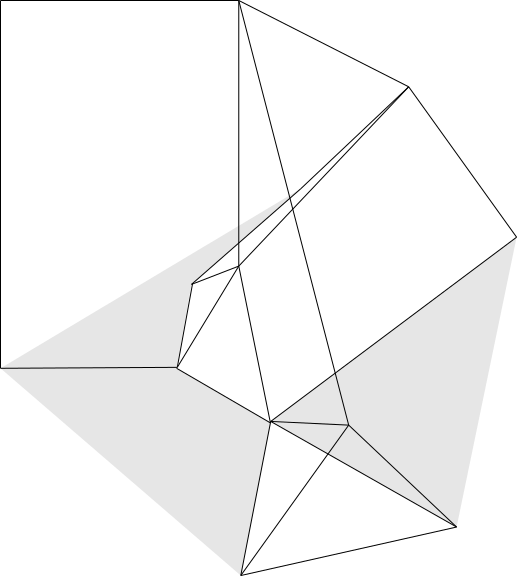}
\caption{A better competitor for the cone $\Y+\Y$ obtained with Brakke's Surface Evolver when $\sin\varphi\le\frac{1}{2}$ (the grey region is the intersection between the competitor and $\Gamma$). The little sloping triangle at the center of the picture represent the new interface obtained by the pinching.}
\label{Y+2Y'-competitore}
\end{figure}

\begin{figure}
\centering
\includegraphics[scale=0.4]{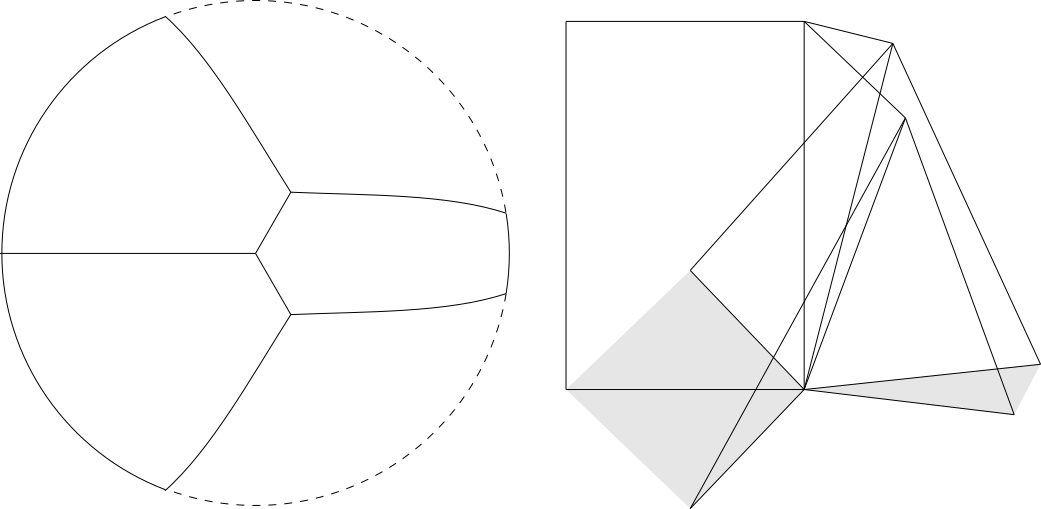}
\caption{On the right the cone $\Y+2\Y$ corresponding to a value of $\sin\varphi$ bigger than $1/2$ (the grey region is the intersection between the cone and $\Gamma$), and on the left its intersection with the hemisphere.}
\label{Y+2Y''}
\end{figure}

\begin{figure}
\centering
\includegraphics[scale=0.4]{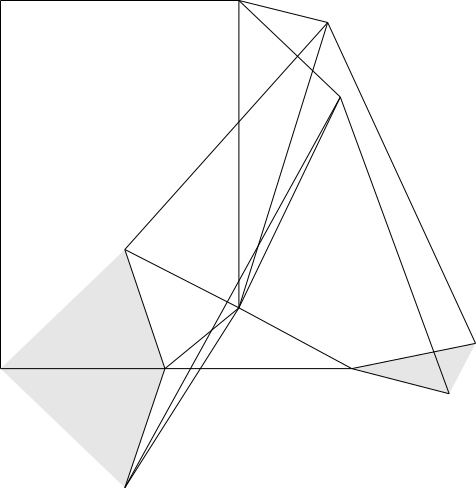}
\caption{A better competitor for the cone $\Y+2\Y$ obtained with Brakke's Surface Evolver when $\sin\varphi\ge\frac{1}{2}$ (the grey region is the intersection between the competitor and $\Gamma$). The little vertical triangle in the center of the picture represents the new interface produced by the pinching.}
\label{Y+2Y''-competitore}
\end{figure}

The second type of cone we will present can be obtained by adding one more sloping $\Y$ to the cone $\Y+\Y$ and then fixing the obtained cone with the proper regions of $\Gamma$ (see Figures \ref{Y+2Y'} and \ref{Y+2Y''}). In order to fulfil the necessary condition for minimality we have to impose the angle formed by the two spines of the sloping $\Y$ with $\Gamma$ to be the same. In this case there are various regions of the complement without a common interface and this means we have different ways to pinch the folds of the cone when looking for a better competitor. When $\sin\varphi\le\frac{1}{2}$ (as in Figure \ref{Y+2Y'}) it is convenient to produce an interface between the region below one of the sloping $\Y$ and the region in some sense opposed to it, and at the same time to pull the other sloping $\Y$ away from the origin by moving its vertex (see Figure \ref{Y+2Y'-competitore}). On the other hand, when $\sin\varphi\ge\frac{1}{2}$ (as in Figure \ref{Y+2Y''}), it is convenient to pinch together the two sloping folds next to each other in such a way as to create a common interface in between the two regions below the two sloping $\Y$ (see Figure \ref{Y+2Y''-competitore}).

\subsection{$\Y+3\Y$}

\begin{figure}
\centering
\includegraphics[scale=0.4]{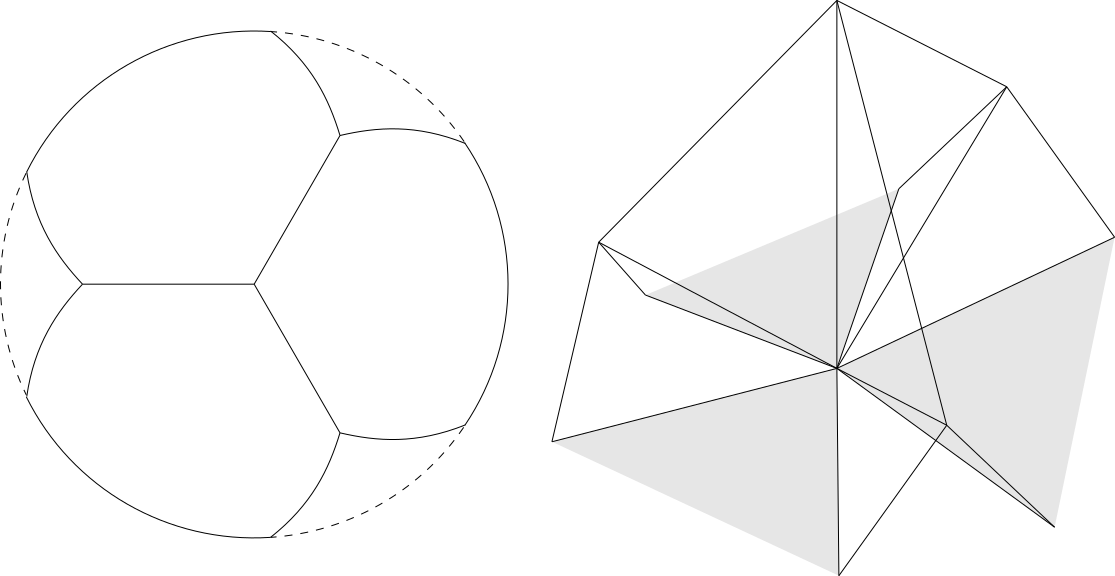}
\caption{On the right the cone $\Y+3\Y$ corresponding to a value of $\sin\varphi$ lower than $1/2$ (the grey region is the intersection between the cone and $\Gamma$), and on the left its intersection with the hemisphere.}
\label{Y+3Y}
\end{figure}

\begin{figure}
\centering
\includegraphics[scale=0.4]{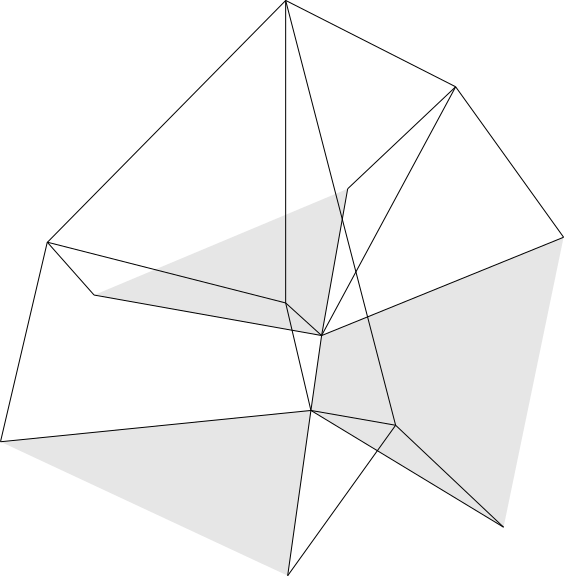}
\caption{A better competitor for the cone $\Y+3\Y$ obtained with Brakke's Surface Evolver when $\sin\varphi\le\frac{1}{2}$ (the grey region is the intersection between the competitor and $\Gamma$). The little sloping triangle in the center of the picture represents the new interface produced by the pinching.}
\label{Y+3Y-competitore}
\end{figure}

\begin{figure}
\centering
\includegraphics[scale=0.4]{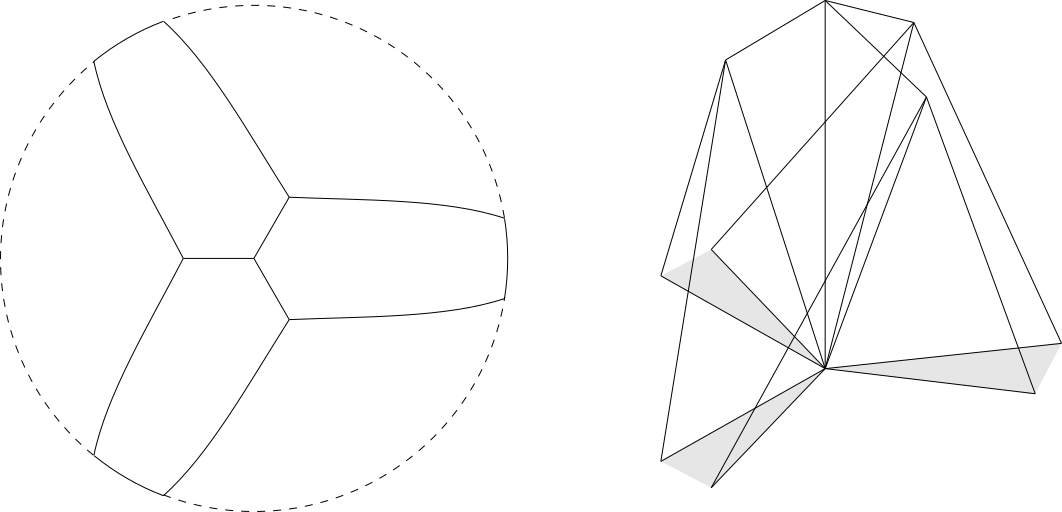}
\caption{On the right the cone $\Y+3\Y$ corresponding to a value of $\sin\varphi$ bigger than $1/2$ (the grey region is the intersection between the cone and $\Gamma$), and on the left its intersection with the hemisphere.}
\label{Y+3Y'}
\end{figure}

\begin{figure}
\centering
\includegraphics[scale=0.4]{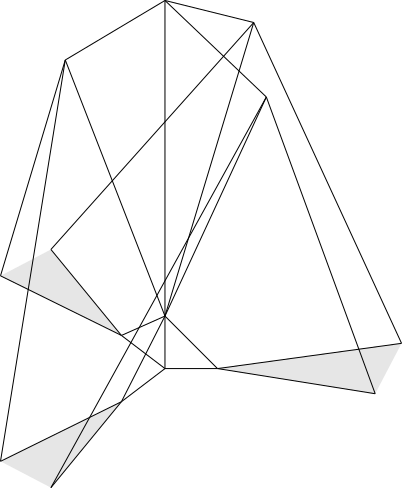}
\caption{A better competitor for the cone $\Y+3\Y$ obtained with Brakke's Surface Evolver when $\sin\varphi\ge\frac{1}{2}$ (the grey region is the intersection between the competitor and $\Gamma$). The three little vertical triangles meeting at the origin are the three new interfaces obtained by pinching the sloping folds with each others.}
\label{Y+3Y'-competitore}
\end{figure}

By iterating one more time the same construction (adding one more sloping $\Y$ and then fixing the cone with the proper regions of $\Gamma$) we obtain the cone $\Y+3\Y$ (see Figures \ref{Y+3Y} and \ref{Y+3Y'}). A better competitor for this type of cone can be found in a similar way as we did for $\Y+2\Y$. That is to say: when $\sin\varphi\le\frac{1}{2}$ (as in Figure \ref{Y+3Y}) it is convenient to produce an interface between the region below one of the sloping $\Y$ and the region opposed to it, and at the same time to pull the other two sloping $\Y$ away from the origin by moving their vertices (see Figure \ref{Y+3Y-competitore}); when $\sin\varphi\ge\frac{1}{2}$ (as in Figure \ref{Y+3Y'}), it is convenient to pinch together the three couples of sloping folds next to each other in such a way as to create the missing interfaces in between the regions below the sloping $\Y$, and then arrange them in such a way that they meet producing the shape of a $\Y$ in a neighbourhood of the origin (see Figure \ref{Y+3Y'-competitore}).



\subsection{$\T+\Y$}
\begin{figure}[h]
\centering
\includegraphics[scale=0.37]{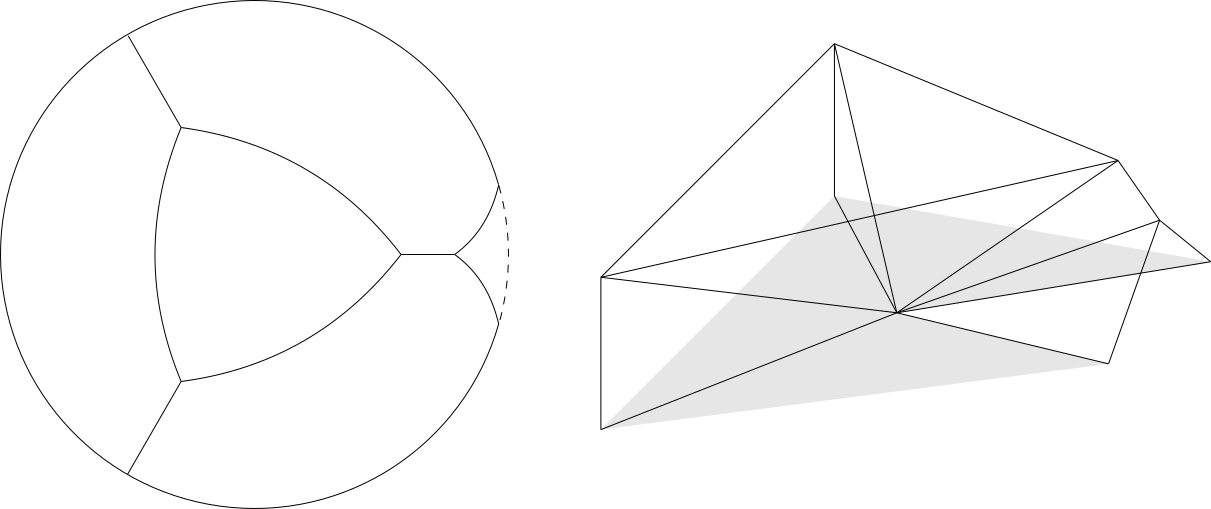}
\caption{On the right the cone $\T+\Y$ (the grey region is the intersection between the cone and $\Gamma$), and on the left its intersection with the hemisphere.}
\label{T+Y}
\end{figure}

We can apply the same technique of adding a sloping $\Y$ to the cone $\T_+$ (see Figure \ref{T+Y}). In this case we obtain a one-parameter family of cones depending on the parameter $\sin\varphi\in[0,1/3]$. A better competitor can be found by pinching the sloping folds of $\T_+$ down to $\Gamma$ (in the same fashion as we did in \ref{Half_T} for a competitor of the cone $\T_+$) and at the same time by pulling the sloping $\Y$ away from the origin (as we did in the previous cases). Because of this it looks inconvenient to add more sloping $\Y$ to the other vertical folds of $\T_+$.
\begin{figure}[h]
\centering
\includegraphics[scale=0.37]{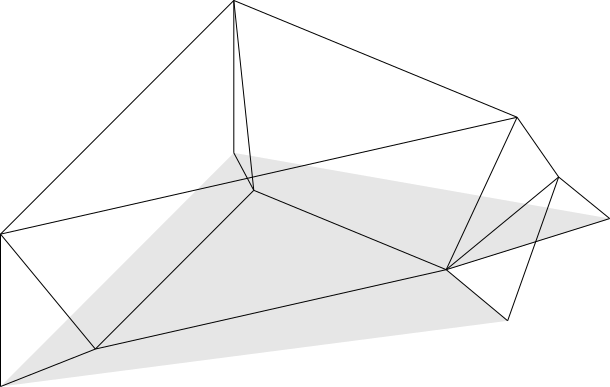}
\caption{A better competitor for the cone $\T+\Y$ obtained with Brakke's Surface Evolver (the grey region is the intersection between the competitor and $\Gamma$).}
\label{T+Y-competitore}
\end{figure}

\newpage
\subsection{Rectangle}
\begin{figure}
\centering
\includegraphics[scale=0.35]{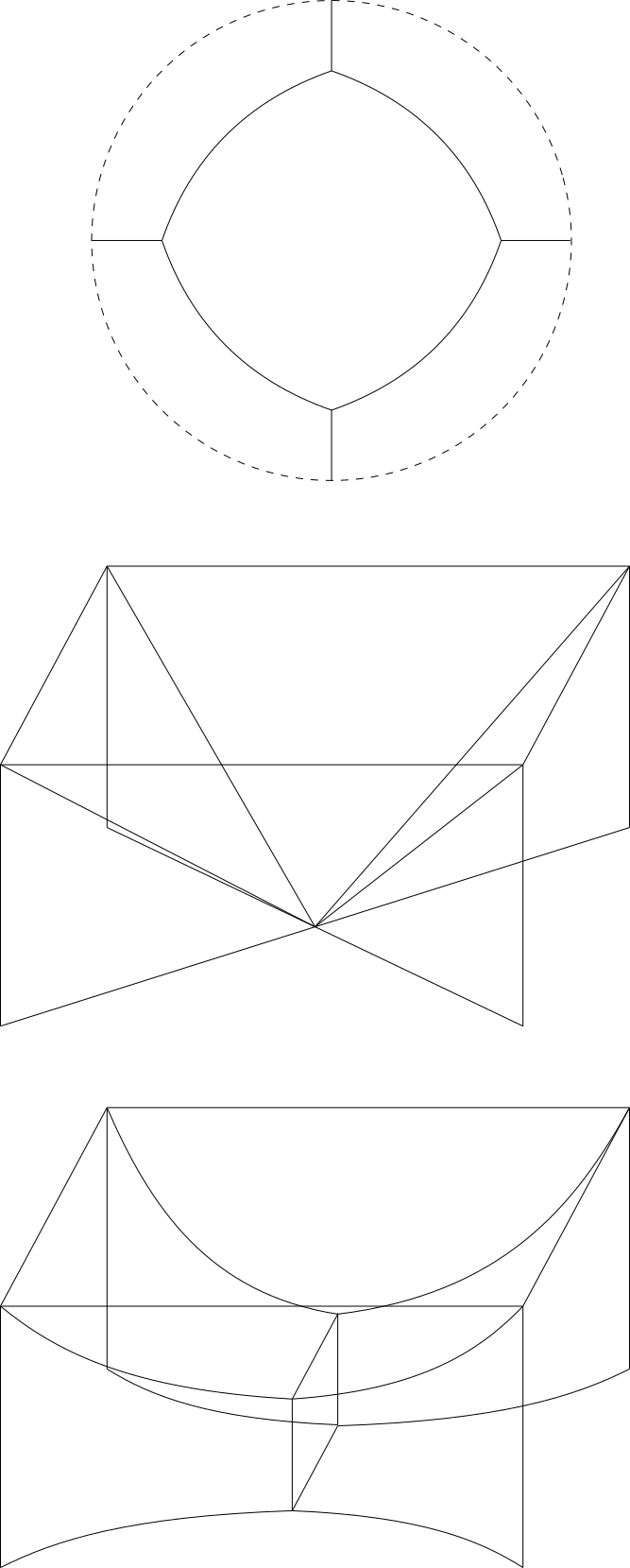}
\caption{In the center the cone $\mathbf{C}_+$, above its intersection with the hemisphere, below a better competitor (both the little square in the center and the bended folds meet $\Gamma$ orthogonally).}
\label{C+}
\end{figure}

We will discuss here some examples of cones such that in the associated network there is a rectangle whose edges are not contained in the equator, starting by the case when the rectangle is actually a square. Let us first consider $\mathbf{C}$, the cone over the edges of a cube that we have already met in Section \ref{Minimal_cones}, and let us call $M$ the better competitor for $\mathbf{C}$ provided by Brakke (see Figure \ref{conocubo}). We can now define $\mathbf{C}_+:=\mathbf{C}\cap\R^3_+$ (see Picture \ref{C+}) as the upper half of the aforementioned cone, than we can construct a better competitor for $\mathbf{C}_+$ starting by $M$. Let $P:=\{(x,y,z)\in\R^3:y=0\}$ be a vertical plane orthogonal to the $x$ axis and $P_+:=\{(x,y,z)\in\R^3:y\ge0\}$ be the half-space bounded by $P$. Let us also denote with $R:\R^3\to\R^3$ the rotation such that $R(P_+)=\R^3_+$. Then by the definition of $M$ it follows that $\H^2(M\cap P)=0$ and it means that $R(M\cap P_+)$ (see Figure \ref{C+}) is a better competitor than the cone $\mathbf{C}_+$ for every $\alpha\in[0,1]$.

\begin{figure}
\centering
\includegraphics[scale=0.35]{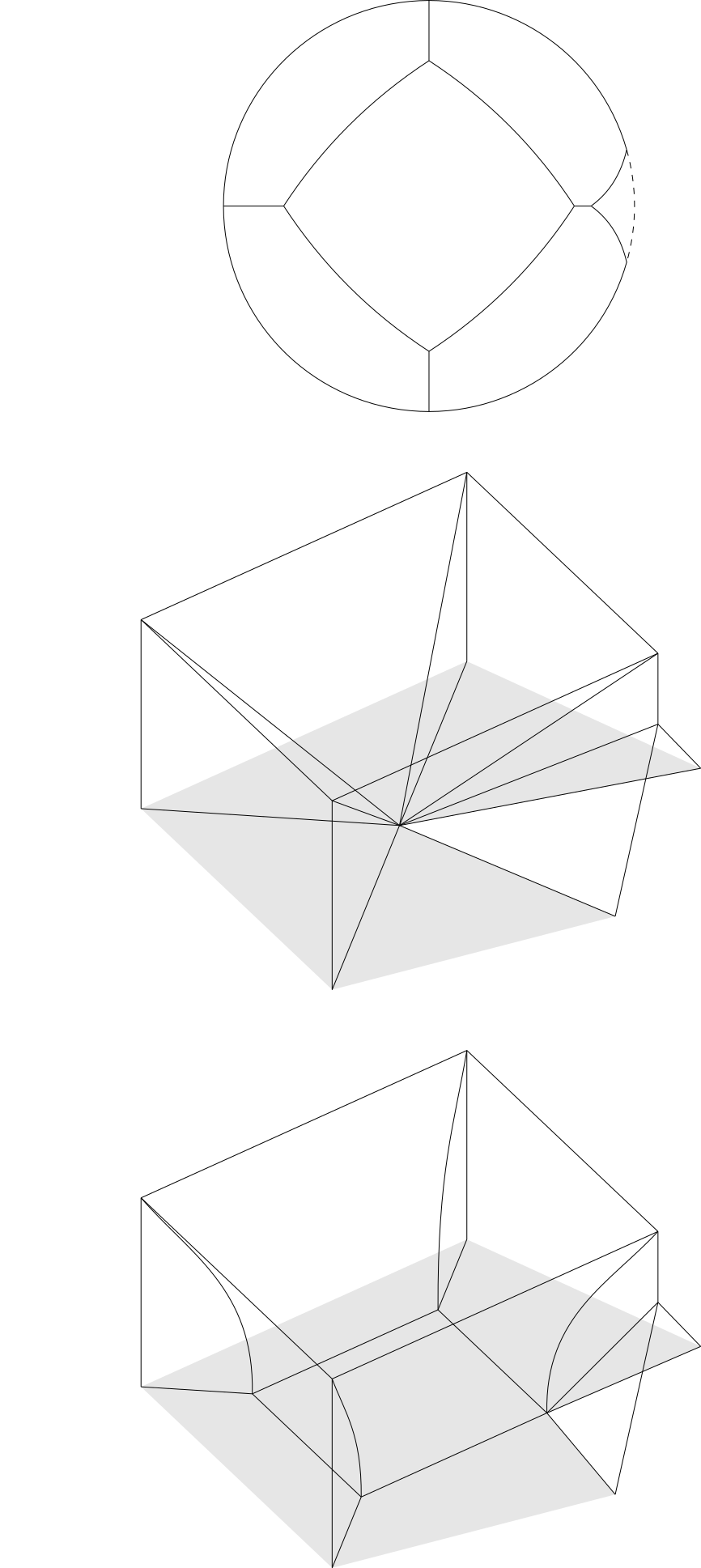}
\caption{In the center the cone $\mathbf{C}_+$ added with a sloping $\Y$, above its intersection with the hemisphere, below a better competitor obtained with Brakke's Surface Evolver (the grey region is the intersection between the competitor and $\Gamma$).}
\label{C+Y}
\end{figure}

By adding a sloping $\Y$ to the cone $\mathbf{C}_+$ (see Figure \ref{C+Y}) we obtain a one-parameter family of cones depending on the parameter $\sin\varphi\in[0,1/\sqrt{2}]$. A better competitor for this cone can be constructed in the same way we did for a better competitor for the cone $\T+\Y$. That is to say by pinching the sloping folds of $\mathbf{C}_+$ down to $\Gamma$ in such a way that they will end up producing an horizontal square, and at the same time by pulling the sloping $\Y$ away from the origin by moving its vertex. Again adding more sloping $\Y$ looks inconvenient.

\begin{figure}
\centering
\includegraphics[scale=0.5]{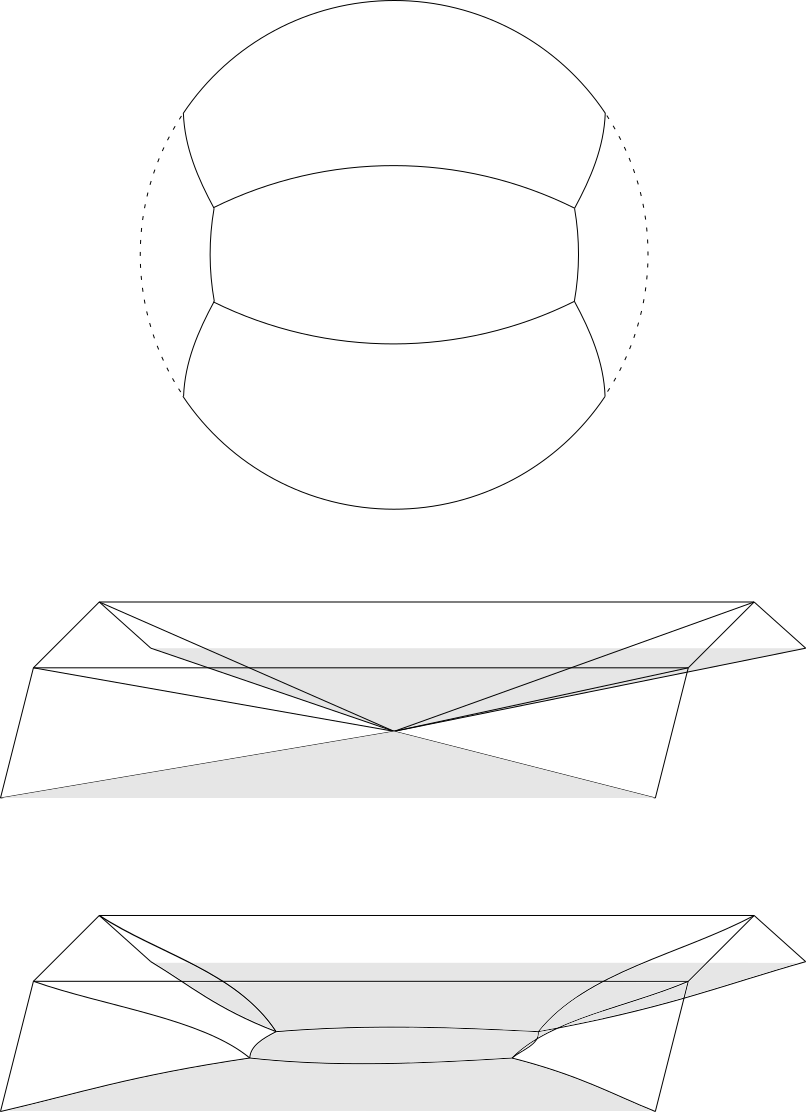}
\caption{In the center $\mathbf{R}$ (the grey region is the intersection between the competitor and $\Gamma$); above its intersection with the hemisphere; below a better competitor obtained with Brakke's Surface Evolver (the grey region is the intersection between the competitor and $\Gamma$).}
\label{rettangolo}
\end{figure}

Let us now consider the case where the sides of the rectangle have different length (which can be computed using \eqref{lati_rettangolo} or \eqref{lati_rettangolo/2}) and we assume it to be in symmetric position with respect to the coordinate axes. By this we mean that the folds generated by the longer edges, the folds generated by the shorter edges, and the four remaining sloping folds, respectively have the same slope (see Figure \ref{rettangolo}). We got a one parameter family of cones, and, by symmetry, we can chose as parameter the slope of either of the folds composing the rectangle. In case we chose the folds generated by the longer edges we have that the parameter is $\sin\varphi\in[1/\sqrt{2},1]$, otherwise, if we chose the folds generated by the shorter edges the parameter is $\sin\psi\in[0,1/\sqrt{2}]$. Once again a better competitor can be obtained by pushing the sloping folds down to $\Gamma$ in such a way as to produce a little ``fat'' rectangle (see Figure \ref{rettangolo}).

\begin{figure}
\centering
\includegraphics[scale=0.35]{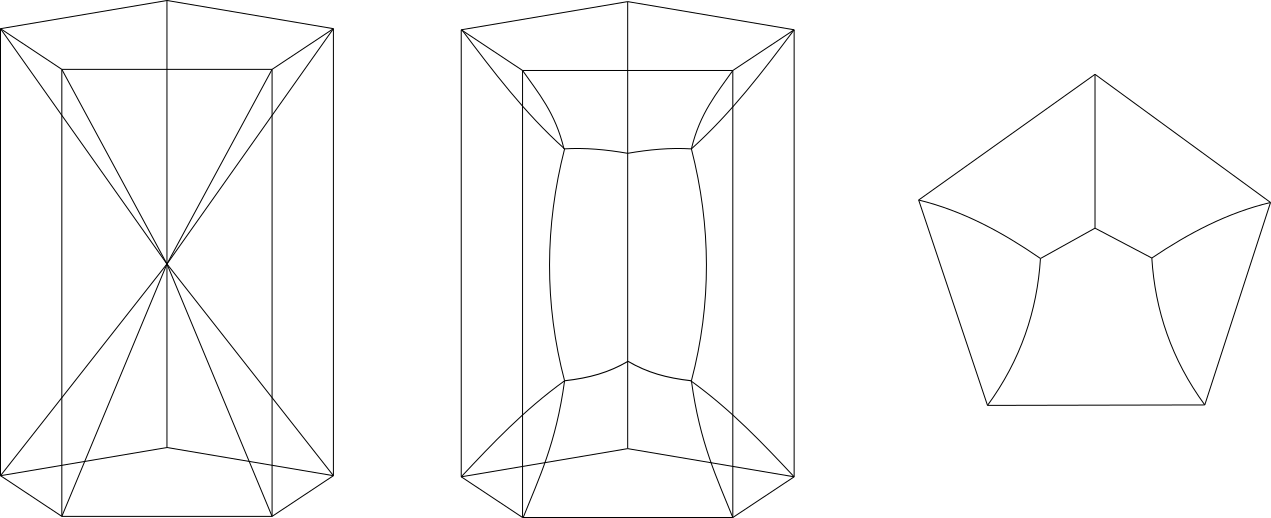}
\caption{On the left the cone $\mathbf{P}$ ; in the center a better competitor $M$; on the right a view of $M$ from above.}
\label{prisma-pentagono}
\includegraphics[scale=0.4]{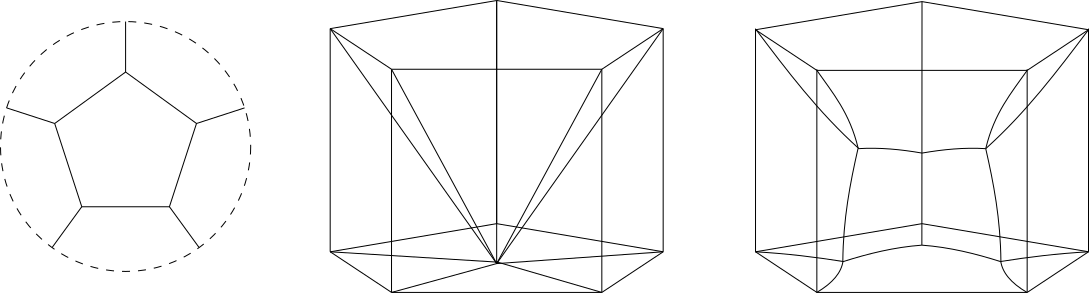}
\caption{In the center cone $\mathbf{P}\cap\R^3_+$, on the left its intersection with the unit hemisphere, and on the right a better competitor $M\cap\R^3_+$.}
\label{competitore-pentagono}
\includegraphics[scale=0.35]{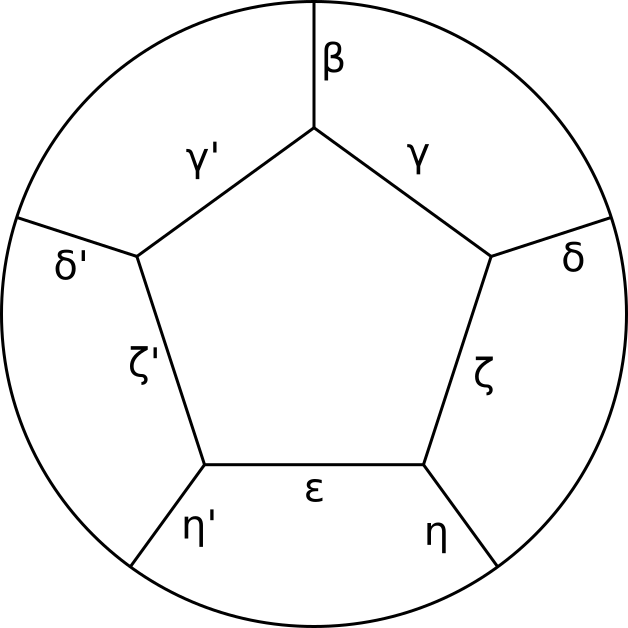}
\caption{The intersection with the unit sphere of a generic pentagonal cone (the equator is not meant to be contained in the network).}
\label{pentagono}
\end{figure}

\subsection{Pentagon}

Finally let us discuss here about pentagonal cones, that is to say cones whose intersection with the unit hemisphere produces a spherical pentagon. We already know that $\mathbf{P}$, the cone over the skeleton of a regular pentagonal prism, is not a minimal cone in $\R^3$ (see Figure \ref{prisma-pentagono}). Moreover can be found a better competitor $M$ such that $\H^2(M\cap\Gamma)=0$. It means that $M\cap\R^3_+$ is a better competitor with respect to the cone $\mathbf{P}\cap\R^3_+$ for every $\alpha\in[0,1]$.

A priori it is possible for more general spherical pentagons to appear in the intersection of a minimal cone with the hemisphere, since the only constraint they have to satisfy is given by condition \eqref{lati_pentagono}. Let us now assume only one spherical pentagon is generated by the network corresponding to a minimal cone. Let us also assume that, beside possible arcs contained in the equator, the network is only composed by the edges of the pentagon and by five more ``radial'' arcs, connecting the vertices of the pentagon to the equator. In this case, using the fact that the radial arcs can meet the equator only with an optimal profile, we can describe this family of pentagonal cones as a two-parameter family as follows.

Let us denote with Greek letters the arcs of the network, with an abuse of notation we will denote both an ark and its length with the same letter. First of all we can remark that, since the radial arcs are in odd number, at least one of them has to meet the equator orthogonally. Assume $\beta$ is a radial arc not meeting the equator orthogonally, then we have two cases. Either $\beta$ meets the equator with the optimal angle $\theta_\alpha$, or there exists another radial arc $\gamma$ such that the two arcs meet at the equator producing an optimal profile of type $V_\theta$. If we are in the first case than in the network there must be an arc $\delta$ adjacent to $\beta$ and contained in the equator. Since $\delta$ cannot be the whole equator it must have a second endpoint, and it means that it has to be adjacent to another radial arc $\eta$ meeting the equator with optimal angle $\theta_\alpha$. Therefore in both cases, given a non orthogonal radial arc we can find a second non orthogonal one. Since the radial arcs are in odd number, by iterating this argument we end up necessarily with an orthogonal radial arc.

\begin{figure}
\centering
\includegraphics[scale=0.35]{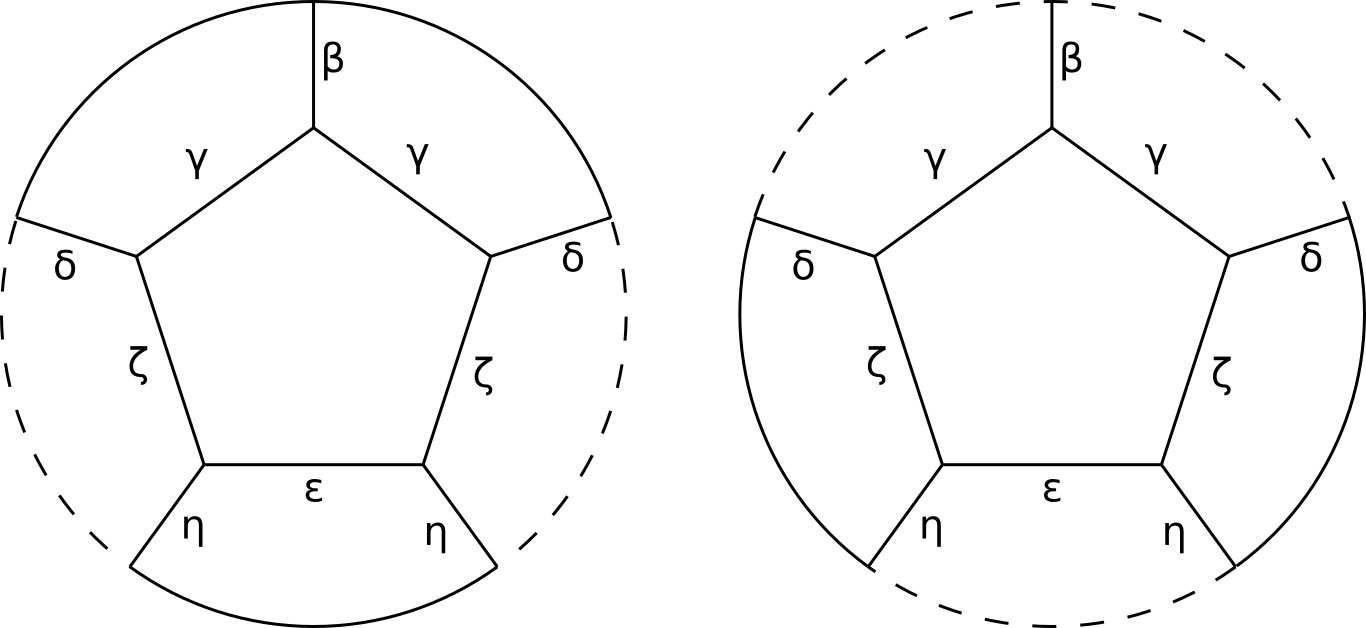}
\caption{The two possible configuration of a pentagonal minimal network with only one pentagon.}
\label{2pentagoni}
\end{figure}

Let us now name the arcs as in Figure \ref{pentagono}, where $\beta$ denotes the orthogonal radial arc found by the previous remark. We remark that, by the $120^\circ$ condition, the arc $\delta$ is completely determined by $\beta$ and $\gamma$, therefore $\delta=\delta(\beta,\gamma)$. The same holds for $\delta'=\delta'(\beta,\gamma')$ and, since every radial arc has to meet the equator with the same angle, $\delta'$ is forced to be symmetric with respect to $\delta$ and it implies $\gamma=\gamma'$. Using now \eqref{lati_pentagono} we can compute $\varepsilon$ in term of $\gamma$ and $\gamma'$, therefore $\varepsilon=\varepsilon(\gamma)$. In the same way we can find a relation between  $\varepsilon$, $\gamma$ and $\zeta$, therefore we can write $\zeta=\zeta(\varepsilon,\gamma)=\zeta(\gamma)$, and by the symmetry between $\gamma$ and $\gamma'$ it follows that also $\zeta$ and $\zeta'$ are symmetric. Therefore the cone must be symmetric with respect to the vertical plane containing $\beta$, and the network must be in one of the two configuration of Figure \ref{2pentagoni}.

Moreover all the edges of the pentagon can be written in term of $\gamma$, and all the radial arcs can be written in term of $\gamma$ and $\beta$. Therefore the pentagonal network satisfying the necessary condition for minimality and whose have only one pentagon can be described as a two-parameter family in term of $\gamma$ and $\beta
$.

\chapter{Higher dimension}\label{higher dimension}
In this chapter we are going to generalise the calibration argument used in Section \ref{Half_T} to the $n$-dimensional case. Let us start by describing the setting. The cone we will consider is contained in the half-space $\R^n_+:=\{(x_1,\cdots,x_n): x_1\ge0\}$ and the domain of the sliding boundary is the hyperplane $\Gamma:=\{x_1=0\}$ bounding $\R^n_+$. In this Chapter the coordinate $x_1$ will play the role of ``vertical direction'', which in the previous Chapter was played by $z$, while the plane $\{x_1=0\}$ plays the role of the ``horizontal hyperplane'', which in the previous chapter was played by $\{z=0\}$. A reference for the properties of simpices used here can be found in Appendix \ref{proprieta_simplessi}.

Let $\Delta^n:=[p_1,\cdots,p_{n+1}]$ be the $n$-dimensional regular simplex whose vertices are
\begin{equation}
\begin{array}{rcccccl}
p_1= \bigg( &-1, &0, &0, &\cdots, &0 &\bigg)\\
p_2= \bigg( &\frac{1}{n}, &-\frac{\sqrt{n^2-1}}{n}, &0, &\cdots, &0 &\bigg)\\
p_3= \bigg( &\frac{1}{n}, &\frac{1}{n}\sqrt{\frac{n+1}{n-1}}, &-\sqrt{\frac{(n+1)(n-2)}{n(n-1)}}, &\cdots, &0 &\bigg)\\
\vdots \phantom{ =\bigg( } &\vdots  &\vdots &\vdots &\ddots &\vdots& \phantom{\bigg) }\\
p_n= \bigg( &\frac{1}{n}, &\frac{1}{n}\sqrt{\frac{n+1}{n-1}}, &\sqrt{\frac{(n+1)}{n(n-1)(n-2)}}, &\cdots, &-p_{nn} &\bigg)\\
p_{n+1}= \bigg( &\frac{1}{n}, &\frac{1}{n}\sqrt{\frac{n+1}{n-1}}, &\sqrt{\frac{(n+1)}{n(n-1)(n-2)}}, &\cdots, &p_{nn} &\bigg).\\
\end{array}
\end{equation}
As we did in Section \ref{Minimal_cones} we define $\mathbf{\Delta}^n:=cone(\sk_{n-2}(\Delta^n))$. Let $\mathbf{\Delta}^n_+:=\mathbf{\Delta}^n\cap\R^n_+$, in the rest of this Chapter we are going to prove the following

\begin{thm}\label{teoremaDelta+}
Let $\alpha_n:=\sqrt{\frac{n+1}{2n}}$. If $\alpha\ge\alpha_n$ then the cone $\mathbf{\Delta}^n_+$ is an $\alpha$-sliding minimiser in the half-space $\R^n_+$ with respect to $\Gamma=\partial\R^n_+$.
\end{thm}

For $i=1,\cdots,n+1$ let $F_i:=\Delta^n_i$  be the $(n-1)$-dimensional face of $\Delta^n$ opposed to the vertex $p_i$, and $F_i^+:=F_i\cap\R^n_+$; in particular $F_1^+=F_1$. Let $M\subset\R^n_+$ be a sliding competitor for $\mathbf{\D}^n_+$ such that the symmetric difference between the two is contained in $\Delta^n_+$ and does not intersect any of the faces $F_i^+$ away from their boundaries. It follows that $\R^n\setminus M$ has 2 unbounded connected components: one of them contains $F_1$, and the other one contains the other $n$ faces $F_i^+$ for $i=2;\cdots,n+1$, as well as the lower half-space $\R^n\setminus\R^n_+$. However $\R^n_+\setminus M$ has $n+1$ unbounded connected components each one of them containing one of the faces $F_i^+$.

For $i=1,\cdots,n+1$ we name $V_i$ the connected component of $\R^n_+\setminus M$ containing $F_i^+$ and we set $V_0:=\R^n\setminus\R^n_+$. In case $\R^n_+\setminus M$ also has some bounded connected components we just include them in $V_1$. By the definition of sliding competitor we have that $M$ is a Lipschitz image of $\mathbf{\Delta}^n_+$, therefore $M$ has locally finite $(n-1)$-dimensional Hausdorff measure. This means that any of the $V_i$ is a set whose perimeter is locally finite. In particular the sets $U_i:=V_i\cap\Delta^n$ have finite perimeter. Let us now introduce the following notation:
\begin{eqnarray}
M_{ij} &:=& \partial^* U_i\cap\partial^* U_j\\
M_i &:=& \bigcup_{j=0}^{n+1}M_{ij}\\
M_i^+ &:=& \bigcup_{j=1}^{n+1}M_{ij}, \textrm{ for } i=1,\cdots,n+1\\
\widetilde{M} &:=& \left(\bigcup_{i=2}^{n+1} M_i^+\right)\cup M_1
\end{eqnarray}
where $i,j=0,\cdots,n+1$ unless otherwise specified. It follows that $\widetilde{M}\subset M\cap\Delta^n_+$ and $\H^{n-1}$-almost every point in $\widetilde{M}$ lies on the interface between exactly two regions of its complement (taking into account also $U_0$). Moreover the interfaces between different couples of regions are essentially disjoint with respect to $\H^{n-1}$.

 Let us now remark the following useful facts
\begin{eqnarray}
M_i &=& M_i^+\cup M_{i0}\\
\widetilde{M}\setminus\Gamma &=& \bigcup_{1\le i<j\le n+1}M_{ij}\quad \textrm{ (up to $\H^{n-1}$-negligible sets)}\\
\widetilde{M}\cap\Gamma &=& M_{10}\quad \textrm{ (up to $\H^{n-1}$-negligible sets)}\\
\partial^*U_i &=& F_i^+\cup M_i=F_i^+\cup M_i^+\cup M_{i0}\\
\partial^*U_0 &=& \partial^*(\Delta^n\setminus\R^n_+)=(\partial^*\Delta^n\setminus\R^n_+)\cup M_0.
\end{eqnarray}
Finally we denote with $n_i$ the exterior unit normal to $\partial U_i$, and $n_{ij}$ will denote the unit normal to $M_{ij}$ pointing in direction of $U_j$.

Let $\ell:=|p_i-p_j|$, the vectors of the calibration we are going to use are defined as follows
\begin{equation}
w_i:=-\frac{p_i}{\ell}=-\sqrt{\frac{n}{2(n+1)}}p_i \textrm{, for } i=1,\cdots, n+1.
\end{equation} 
In components we have
\begin{equation}
\begin{array}{rcccccl}\label{calibrazione_Delta+}
w_1= \bigg( &\sqrt{\frac{n}{2(n+1)}}, &0, &0, &\cdots, &0 &\bigg)\\
w_2= \bigg( &\frac{-1}{\sqrt{2n(n+1)}}, &\sqrt{\frac{n-1}{2n}}, &0, &\cdots, &0 &\bigg)\\
w_3= \bigg( &\frac{-1}{\sqrt{2n(n+1)}}, &\frac{-1}{\sqrt{2n(n-1)}}, &\sqrt{\frac{(n-2)}{2(n-1)}}, &\cdots, &0 &\bigg)\\
\vdots \phantom{ =\bigg( } &\vdots  &\vdots &\vdots &\ddots &\vdots& \phantom{\bigg) }\\
w_n= \bigg( &\frac{-1}{\sqrt{2n(n+1)}}, &\frac{-1}{\sqrt{2n(n-1)}}, &\frac{-1}{\sqrt{2(n-1)(n-2)}}, &\cdots, &\frac{p_{nn}}{\ell} &\bigg)\\
w_{n+1}= \bigg( &\frac{-1}{\sqrt{2n(n+1)}}, &\frac{-1}{\sqrt{2n(n-1)}}, &\frac{-1}{\sqrt{2(n-1)(n-2)}}, &\cdots, &-\frac{p_{nn}}{\ell} &\bigg).\\
\end{array}
\end{equation}
Let us remark that for $i=1,\cdots,n+1$ we have that $w_i\perp F_i$ and $|w_i|=1/\ell$; while in case $M=\mathbf{\Delta}^n_+$ we have that $w_j-w_i=n_{ij}$ for $1\le i<j\le n+1$.

We are now ready to start the calibration argument by applying the divergence theorem to the sets $U_i$ with the vectors $w_i$ as follows
\begin{equation}\label{calcoli_calibrazione_Delta+}
\begin{aligned}
&\frac{1}{\ell}\H^{n+1}\left(\cup_iF^+_i\right)=\sum_{i=1}^{n+1}\int_{F_i^+}w_i\cdot n_id\mathcal{H}^{n-1}=-\sum_{i=1}^{n+1}\int_{M_i}w_i\cdot n_id\mathcal{H}^{n-1}\\
&=-\sum_{i=1}^{n+1}\int_{M_i^+}w_i\cdot n_id\mathcal{H}^{n-1}-\sum_{i=1}^{n+1}\int_{M_{i0}}w_i\cdot n_{i0}d\mathcal{H}^{n-1}\\
&=\sum_{1\le i<j\le n+1}\int_{M_{ij}}(w_j-w_i)\cdot n_{ij}d\mathcal{H}^{n-1}+\sum_{i=1}^{n+1}\int_{M_{i0}}w_i\cdot \hat{x}_1d\mathcal{H}^{n-1}
\end{aligned}
\end{equation}
Let us now focus on the second of the two sums in the last line. Using \eqref{calibrazione_Delta+} we can see that
\begin{equation}
(w_1, \hat{x}_1)=\sqrt{\frac{n}{2(n+1)}},
\end{equation}
while for $i=2,\cdots,n+1$ we have
\begin{equation}
(w_i, \hat{x}_1)=(w_2, \hat{x}_1)=\frac{-1}{\sqrt{2n(n+1)}}.
\end{equation}
therefore we get
\begin{equation}
\begin{aligned}\label{calcoli_interfaccia_Delta+}
\sum_{i=1}^{n+1}\int_{M_{i0}}w_i\cdot \hat{x}_1d\mathcal{H}^{n-1} &=\sum_{i=2}^{n+1}(w_i, \hat{x}_1)\mathcal{H}^{n-1}(M_{i0})+(w_1, \hat{x}_1)\mathcal{H}^{n-1}(M_{40})\\
&=(w_2, \hat{x}_1)\mathcal{H}^{n-1}\left(\bigcup_{i=2}^{n+1}M_{i0}\right)+(w_1, \hat{x}_1)\mathcal{H}^{n-1}(M_{10}).
\end{aligned}
\end{equation}
Using the fact that
\begin{equation}
\H^{n-1}\left(\bigcup_{i=2}^{n+1}M_{i0}\right)=\H^{n-1}(M_0)-\H^{n-1}(M_{10})
\end{equation}
the equation \eqref{calcoli_interfaccia_Delta+} becomes
\begin{equation}
\begin{aligned}
\sum_{i=1}^{n+1}\int_{M_{i0}}w_i\cdot \hat{x}_1d\mathcal{H}^{n-1} &=(w_2, \hat{x}_1)\H^{n-1}(M_0)+(w_1-w_2, \hat{x}_1)\mathcal{H}^{n-1}(M_{10}).
\end{aligned}
\end{equation}
Therefore plugging the previous one in \eqref{calcoli_calibrazione_Delta+} we obtain
\begin{equation}
\begin{aligned}
&\frac{1}{\ell}\mathcal{H}^2(\cup_iF_i^+)-(w_2, \hat{x}_1)\H^{n-1}(M_0)=\\
&=\sum_{1\le i<j\le n+1}\int_{M_{ij}}(w_j-w_i)\cdot n_{ij}d\mathcal{H}^{n-1}+(w_1-w_2, \hat{x}_1)\mathcal{H}^{n-1}(M_{10})\\
&\le\H^{n-1}\left(\widetilde{M}\setminus\Gamma\right)+(w_1-w_2, \hat{x}_1)\H^{n-1}\left(\widetilde{M}\cap\Gamma\right).
\end{aligned}
\end{equation}
Let us now set
\begin{equation}
\alpha=(w_1-w_2, \hat{x}_1)=\sqrt{\frac{n}{2(n+1)}}+\frac{1}{\sqrt{2n(n+1)}}=\sqrt{\frac{n+1}{2n}}=\alpha_n,
\end{equation}
recalling that $\widetilde{M}\subset M\cap\Delta^n$ the previous inequality becomes
\begin{equation}\label{calibrazione_quasi_fatto_Delta+}
\frac{1}{\ell}\mathcal{H}^2(\cup_iF_i^+)-(w_2, \hat{x}_1)\H^{n-1}(M_0)\le J_\alpha(\widetilde{M})\le J_\alpha(M\cap\D^n_+).
\end{equation}
Since the left-hand side of \eqref{calibrazione_quasi_fatto_Delta+} is a constant, and the chain of inequalities turns into a chain of equalities when $M=\mathbf{\D}^n_+$, we proved that this cone is minimal when $\alpha=\alpha_n$. Moreover is is also minimal for $\alpha'\ge\alpha_n$ and it is due to the fact that $J_{\alpha'}(\mathbf{\D}^n_+\cap\D^n)=J_{\alpha}(\mathbf{\D}^n_+\cap\D^n)$ because $\H^{n-1}(\mathbf{\D}^n_+\cap\Gamma)=0$. To show the $\alpha'$-minimality of $\mathbf{\D}^n_+$ we can compute as follows:
\begin{equation}
J_{\alpha'}(\mathbf{\D}^n_+\cap\D^n)=J_{\alpha}(\mathbf{\D}^n_+\cap\D^n)\le J_{\alpha}(M\cap\D^n)\le J_{\alpha'}(M\cap\D^n)
\end{equation}
for every sliding competitor $M$.

\begin{appendices}
\chapter{Properties of simplices}\label{proprieta_simplessi}

In this Appendix we will describe some general properties about simplices that will be useful throughout the different chapters.

For $n\le m$, let us take $p_1,\cdots,p_{n+1}\in\R^m$, and let us assume that there is no affine space $H$ with dimension strictly lower than $n$ such that $p_i\in H$ for every $i=1,\cdots,n+1$. We denote with $[p_1,\cdots,p_{n+1}]$ the convex hull of such points and we call it the $n$-dimensional simplex with vertices $p_1,\cdots,p_{n+1}$. Let  $\Delta^n=[p_1,\cdots,p_{n+1}]$ be an $n$-dimensional simplex, and  $\varphi:\R^m\to\R^N$ be a linear function, then $\varphi(\Delta^n)$ is a simplex whose dimension is not bigger than $n$, indeed
\begin{equation}
\varphi(\Delta^n)=\varphi([p_1,\cdots,p_{n+1}])=[\varphi(p_1),\cdots,\varphi(p_{n+1})].
\end{equation}
Given and an index $i\in\{1,\cdots,n+1\}$ we can define
\begin{equation}
\Delta^n_i:=[p_1,\cdots,\hat{p}_i,\cdots,p_{n+1}].
\end{equation}
The simplex $\Delta^n_i$ has dimension $(n-1)$ and its vertices are the same as those of $\Delta^n$ except for $p_i$. In particular $\Delta^n_i$ is the $(n-1)$-dimensional face of $\Delta^n$ opposed to the vertex $p_i$. Analogously, for $1\le i<j\le n+1$ we can define
\begin{equation}
\Delta^n_{ij}:=[p_1,\cdots,\hat{p}_i,\cdots,\hat{p}_j,\cdots,p_{n+1}]
\end{equation}
as the $(n-2)$-dimensional simplex whose vertices are the same as those of $\Delta^n$ except for $p_i$ and $p_j$ and again $\Delta^n_{ij}$ is one of the $(n-2)$-dimensional faces of $\Delta^n$. For $k=0,\cdots,n$ this construction can be iterated $k$ times providing the $(n-k)$-dimensional faces of $\Delta^n$. Therefore we can define the $k$-dimensional skeleton of $\Delta^n$ as the union of all of its $(n-k)$-dimensional faces. That is to say
\begin{equation}
\sk_{(n-k)}(\Delta^n):=\bigcup_{0\le i_1<\cdots<i_k\le n}\Delta^n_{i_1,\cdots,i_k}.
\end{equation}

In the following of this section we will discuss the geometric properties of regular simplices. Let $\Delta^n=[p_1,\cdots,p_{n+1}]\subset\R^n$ be a regular simplex. Up to translations and dilations we can assume it to be \emph{centred at the origin}, that is to say its barycentre is the origin, and to be \emph{unitary}, that is to say its vertices lie on a unitary sphere centred at the origin. First of all one would like to know what do the vertices of such a simplex look like in coordinates. The requirements for $\Delta^n$ of being centred and unitary give us the following two condition on the vertices
\begin{equation}\label{condizioni_simplesso}
\frac{1}{n+1}\sum_{i=1}^{n+1}p_i=0 \; \textrm{ and } \; |p_i|=1 \; \textrm{ for } \; 1\le i\le n+1.
\end{equation}
Using the conditions in \eqref{condizioni_simplesso} we will now try to explicitly construct a regular unitary simplex centred at the origin. By the first condition we can write the vertices as follows 
\begin{equation}\label{simplesso_canonico1}
\begin{array}{rcccccccl}
p_1= \bigg( &p_{11}, &0, &0, &0, &\cdots &\cdots, &0 &\bigg)\\
p_2= \bigg( &-\frac{p_{11}}{n}, &p_{22}, &0, &0, &\cdots &\cdots, &0 &\bigg)\\
p_3= \bigg( &-\frac{p_{11}}{n}, &-\frac{p_{22}}{n-1}, &p_{33}, &0, &\cdots &\cdots, &0 &\bigg)\\
\vdots \phantom{ =\bigg( } &\\
p_k= \bigg( &-\frac{p_{11}}{n}, &\cdots, &\frac{p_{(k-1)(k-1)}}{n-k+1}, &p_{kk}, &0, &\cdots, &0 &\bigg)\\
\vdots \phantom{ =\bigg( } &\\
p_n= \bigg( &-\frac{p_{11}}{n}, &-\frac{p_{22}}{n-1}, &-\frac{p_{33}}{n-2}, &-\frac{p_{44}}{n-3}, &\cdots &\cdots, &p_{nn} &\bigg)\\
p_{n+1}= \bigg( &-\frac{p_{11}}{n}, &-\frac{p_{22}}{n-1}, &-\frac{p_{33}}{n-2}, &-\frac{p_{44}}{n-3}, &\cdots &\cdots, &-p_{nn} &\bigg).
\end{array}
\end{equation}
Now we can chose all the components $p_{ii}$ to be negative and imposing the second condition we get
\begin{equation}\label{simplesso_canonico2}
\begin{array}{rcccccl}
p_1= \bigg( &-1, &0, &0, &\cdots, &0 &\bigg)\\
p_2= \bigg( &\frac{1}{n}, &-\frac{\sqrt{n^2-1}}{n}, &0, &\cdots, &0 &\bigg)\\
p_3= \bigg( &\frac{1}{n}, &\frac{1}{n}\sqrt{\frac{n+1}{n-1}}, &-\sqrt{\frac{(n+1)(n-2)}{(n-1)}}, &\cdots, &0 &\bigg)\\
\vdots \phantom{ =\bigg( } &\vdots  &\vdots &\vdots &\ddots &\vdots& \phantom{\bigg) }\\
p_n= \bigg( &\frac{1}{n}, &\frac{1}{n}\sqrt{\frac{n+1}{n-1}}, &\sqrt{\frac{(n+1)}{n(n-1)(n-2)}}, &\cdots, &p_{nn} &\bigg)\\
p_{n+1}= \bigg( &\frac{1}{n}, &\frac{1}{n}\sqrt{\frac{n+1}{n-1}}, &\sqrt{\frac{(n+1)}{n(n-1)(n-2)}}, &\cdots, &-p_{nn} &\bigg).
\end{array}
\end{equation}
We will say that a simplex is \emph{canonical} or that it is in \emph{canonical position} if, up to orientations of $\R^n$, its vertices can be written as in \eqref{simplesso_canonico2}. Any other regular unitary simplex centred at the origin can be obtained as the rotation of a canonical simplex. Therefore in the following we will always assume $\Delta^n$ to be in canonical position. Let us now compute the length of a one-dimensional edge
\begin{equation}
\ell := |p_i-p_j|=\sqrt{\frac{2(n+1)}{n}}
\end{equation}
where $1\le i<j\le n+1$.

Given a simplex $\Delta^n$ in canonical position as in \eqref{simplesso_canonico2} let us now define the following cone 
\begin{equation}
\mathbf{\Delta}^n:=cone(\sk_{n-2}(\Delta^n)).
\end{equation}
First of all, since the notation might be misleading, let us remark that $\mathbf{\Delta}^n$ has dimension $n-1$. We also remark that, by definition
\begin{equation}
\mathbf{\Delta}^n=\bigcup_{1\le i<j\le n+1}cone(\Delta^n_{ij}),
\end{equation}
therefore the cone $\mathbf{\Delta}^n$ is composed of $(n+1)n/2$ folds, each one contained in a different hyperplane. It means that the intersection between two different folds is negligible with respect to the $(n-1)$-dimensional Hausdorff measure. In the following, abusing the notation, we might implicitly assume $\mathbf{\Delta}^n$ to be the intersection of the cone with the simplex itself.

Let us now prove an important property of regular simplices that plays a key role in the calibration argument. Let $1\le i<j\le n+1$ then the vector $p_i-p_j$ is orthogonal to the $(n-1)$-dimensional cone over the $(n-2)$-dimensional face $\Delta^n_{ij}$. Up to orientations of $\R^n$ we can assume the indices $i$ and $j$ to be $1$ and $2$. Since the cone over the face $\Delta^n_{12}$ is contained in the hyperplane spanned by the vectors $p_3,\cdots,p_{n+1}$ it will be enough to prove that $(p_1-p_2)\perp p_k$ for $k=3,\cdots,p_{n-1}$. Again up to orientation of $\R^n$ it will be sufficient to prove that $(p_1-p_2)\perp p_3$. Using \eqref{simplesso_canonico2} we have 
\begin{equation}
p_1-p_2=\left(\frac{1+n}{n},\frac{\sqrt{n^2-1}}{n},0,\cdots,0\right),
\end{equation}
we can compute the scalar product with $p_3$
\begin{equation}
((p_1-p_2),p_3)=-\frac{1+n}{n^2}+\frac{1}{n}\frac{\sqrt{n^2-1}}{n}\sqrt{\frac{n+1}{n-1}}=-\frac{1+n}{n^2}+\frac{1+n}{n^2}=0
\end{equation}
and the property is proved.

\end{appendices}


\bibliography{bibliografia}{}

\begin{thebibliography}{10}

\bibitem{almgren1968existence}
Frederick~J Almgren.
\newblock Existence and regularity almost everywhere of solutions to elliptic
  variational problems among surfaces of varying topological type and
  singularity structure.
\newblock {\em Annals of Mathematics}, pages 321--391, 1968.

\bibitem{ambrosio2000functions}
Luigi Ambrosio, Nicola Fusco, and Diego Pallara.
\newblock {\em Functions of bounded variation and free discontinuity problems},
  volume 254.
\newblock Clarendon Press Oxford, 2000.

\bibitem{brakke2013surface}
Kenneth~A Brakke.
\newblock Surface {E}volver, {V}ersion 2.70.
\newblock \url{http://facstaff.susqu.edu/brakke/evolver/evolver.html}.
\newblock Accessed: 01-03-2018.

\bibitem{brakke1991minimal}
Kenneth~A Brakke.
\newblock Minimal cones on hypercubes.
\newblock {\em The Journal of Geometric Analysis}, 1(4):329--338, 1991.

\bibitem{david2013regularity}
Guy David.
\newblock Regularity of minimal and almost minimal sets and cones: J.
  taylor’s theorem for beginners.
\newblock {\em Analysis and geometry of metric measure spaces}, 56:67--117,
  2013.

\bibitem{david2014local}
Guy David.
\newblock Local regularity properties of almost-and quasiminimal sets with a
  sliding boundary condition.
\newblock {\em arXiv preprint arXiv:1401.1179}, 2014.

\bibitem{david2014should}
Guy David.
\newblock Should we solve {P}lateau’s problem again.
\newblock {\em Advances in Analysis: The Legacy of Elias M. Stein. Edited by
  Charles Fefferman, Alexandru D. Ionescu, DH Phong, and Stephen Wainger,
  Princeton Mathematical Series}, 50:108--145, 2014.

\bibitem{david2000uniform}
Guy David and Stephen Semmes.
\newblock {\em Uniform rectifiability and quasiminimizing sets of arbitrary
  codimension}, volume 687.
\newblock American Mathematical Soc., 2000.

\bibitem{douglas1931solution}
Jesse Douglas.
\newblock Solution of the problem of {P}lateau.
\newblock {\em Transactions of the American Mathematical Society},
  33(1):263--321, 1931.

\bibitem{dugundji1970topology}
James Dugundji.
\newblock Topology boston, 1970.

\bibitem{evans1991measure}
Lawrence~Craig Evans and Ronald~F Gariepy.
\newblock {\em Measure Theory and Fine Properties of Functions}, volume~5.
\newblock CRC Press, 1991.

\bibitem{fang2016holder}
Yangqin Fang.
\newblock H{\"o}lder regularity at the boundary of two-dimensional sliding
  almost minimal sets.
\newblock {\em Advances in Calculus of Variations}, 2016.

\bibitem{fang2017local}
Yangqin Fang.
\newblock Local ${C}^{1,\beta}$-regularity at the boundary of two dimensional
  sliding almost minimal sets in $\mathbb{R}^3$.
\newblock arXiv:1611.01343, 2017.

\bibitem{federer2014geometric}
Herbert Federer.
\newblock {\em Geometric measure theory}.
\newblock Springer, 2014.

\bibitem{federer1960normal}
Herbert Federer and Wendell~H Fleming.
\newblock Normal and integral currents.
\newblock {\em Annals of Mathematics}, pages 458--520, 1960.

\bibitem{finn1974capillarity}
Robert Finn.
\newblock Capillarity phenomena.
\newblock {\em Russian Mathematical Surveys}, 29(4):133--153, 1974.

\bibitem{giusti1976boundary}
Enrico Giusti.
\newblock Boundary value problems for non-parametric surfaces of prescribed
  mean curvature.
\newblock {\em Annali della Scuola Normale Superiore di Pisa-Classe di
  Scienze}, 3(3):501--548, 1976.

\bibitem{lawlor1994paired}
Gary Lawlor and Frank Morgan.
\newblock Paired calibrations applied to soap films, immiscible fluids, and
  surfaces or networks minimizing other norms.
\newblock {\em Pacific Journal of Mathematics}, 166(1):55--83, 1994.

\bibitem{liang2012almgren}
Xiangyu Liang.
\newblock Almgren-minimality of unions of two almost orthogonal planes in
  $\mathbb{R}^4$.
\newblock {\em Proceedings of the London Mathematical Society},
  106(5):1005--1059, 2012.

\bibitem{liang2014almgren}
Xiangyu Liang.
\newblock Almgren and topological minimality for the set y$\times$ y.
\newblock {\em Journal of Functional Analysis}, 266(10):6007--6054, 2014.

\bibitem{liang2015topological}
Xiangyu Liang.
\newblock On the topological minimality of unions of planes of arbitrary
  dimension.
\newblock {\em International Mathematics Research Notices},
  2015(23):12490--12539, 2015.

\bibitem{marchese2016steiner}
Andrea Marchese and Annalisa Massaccesi.
\newblock The steiner tree problem revisited through rectifiable g-currents.
\newblock {\em Advances in Calculus of Variations}, 9(1):19--39, 2016.

\bibitem{massaccesi2014currents}
Annalisa Massaccesi.
\newblock {\em Currents with coefficients in groups, applications and other
  problems in Geometric Measure Theory}.
\newblock PhD thesis, Ph. D. thesis, Scuola Normale Superiore di Pisa, 2014.

\bibitem{mattila1999geometry}
Pertti Mattila.
\newblock {\em Geometry of sets and measures in Euclidean spaces: fractals and
  rectifiability}, volume~44.
\newblock Cambridge university press, 1999.

\bibitem{dephilippis2015regularity}
G~De Philippis and Francesco Maggi.
\newblock Regularity of free boundaries in anisotropic capillarity problems and
  the validity of young’s law.
\newblock {\em Archive for Rational Mechanics and Analysis}, 216(2):473--568,
  2015.

\bibitem{plateau1873statique}
Joseph Antoine~Ferdinand Plateau.
\newblock {\em Statique exp{\'e}rimentale et th{\'e}orique des liquides soumis
  aux seules forces mol{\'e}culaires}, volume~2.
\newblock Gauthier-Villars, 1873.

\bibitem{reifenberg1960solution}
Ernst~Robert Reifenberg.
\newblock Solution of the plateau problem form-dimensional surfaces of varying
  topological type.
\newblock {\em Acta Mathematica}, 104(1-2):1--92, 1960.

\bibitem{taylor1973regularity}
Jean~E Taylor.
\newblock Regularity of the singular sets of two-dimensional area-minimizing
  flat chains modulo 3 inr 3.
\newblock {\em Inventiones mathematicae}, 22(2):119--159, 1973.

\bibitem{taylor1976structure}
Jean~E Taylor.
\newblock The structure of singularities in soap-bubble-like and soap-film-like
  minimal surfaces.
\newblock {\em Annals of Mathematics}, pages 489--539, 1976.

\bibitem{taylor1977boundary}
Jean~E Taylor.
\newblock Boundary regularlty for solutions to various capillarity and free
  boundary problems.
\newblock {\em Communications in Partial Differential Equations},
  2(4):323--357, 1977.

\end{thebibliography}
\bibliographystyle{plain}

\newpage

\phantom{pain}

\newpage
\begin{figure}[h]
\vspace*{-5cm}
\hspace*{-4cm}
\includegraphics[scale=0.37]{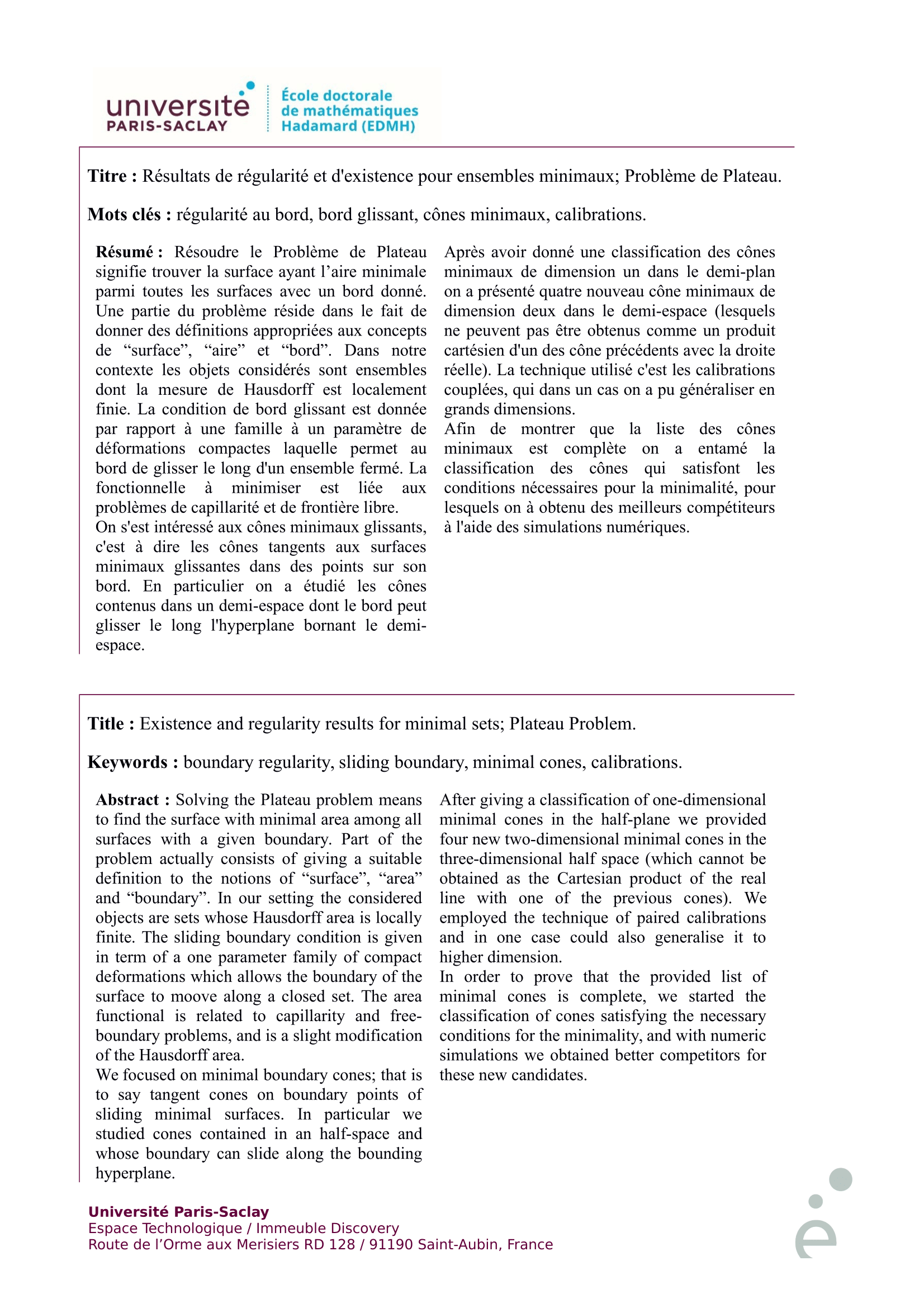}
\end{figure}

\end{document}